\documentclass[sort&compress,times,preprint]{elsarticle}

\usepackage{url}
\usepackage{xcolor}

\definecolor{gold}{rgb}{0.85, 0.65, 0.13}
\global\long\def\red#1{\textcolor{black}{#1}}%
\global\long\def\violet#1{\textcolor{black}{#1}}%
\global\long\def\blue#1{\textcolor{black}{#1}}%

\usepackage[hidelinks, colorlinks]{hyperref}
\usepackage{bookmark}
\bookmarksetup{
	numbered,
	addtohook={%
		\ifnum\bookmarkget{level}>1 %
		\bookmarksetup{numbered=false}%
		\fi
	},
}

\usepackage{framed,multirow}
\usepackage[margin=25mm]{geometry}
\global\bibsep=0pt

\usepackage{amssymb}
\usepackage{latexsym}

\usepackage{amsfonts}
\usepackage{mathrsfs}
\usepackage{upgreek}

\usepackage{graphicx}
\usepackage{epstopdf}
\usepackage{algorithmic}
\ifpdf
\DeclareGraphicsExtensions{.eps,.pdf,.png,.jpg}
\else
\DeclareGraphicsExtensions{.eps}
\fi
\usepackage[export]{adjustbox}

\usepackage{amsmath, amsthm}
\usepackage{xspace}
\usepackage{bm}
\usepackage{physics}
\usepackage{stmaryrd}
\usepackage{cleveref}
\usepackage{enumerate}
\usepackage[labelfont=bf,textfont=normalfont]{caption}
\usepackage{subcaption}
\usepackage{mathtools}

\usepackage{booktabs}
\usepackage{multirow, multicol, makecell}
\usepackage{longtable}
\usepackage{threeparttable}

\newtheorem{theorem}{Theorem}

\newtheorem{lemma}{Lemma}
\newtheorem{assumption}{Assumption}
\newtheorem{remark}{Remark}

\newtheorem{definition}{Definition}

\crefname{assumption}{Assumption}{Assumptions}
\crefname{definition}{Definition}{Definitions}
\crefname{lemma}{Lemma}{Lemmas}
\crefname{remark}{Remark}{Remarks}
\crefname{theorem}{Theorem}{Theorems}
\crefname{proposition}{Proposition}{Propositions}
\crefname{section}{Section}{Sections}
\crefname{figure}{Fig.}{Figs.}
\crefname{equation}{}{}
\crefname{table}{Table}{Tables}
\crefname{appendix}{}{}

\newcommand{\fnc}[1]{\ensuremath{\mathcal{#1}}}
\newcommand{\vecfnc}[1]{\ensuremath{\boldsymbol{\mathcal{#1}}}} 

\newcommand{\uk}[0]{\ensuremath{\bm{{u}}_{k}}}

\renewcommand{\H}[0]{\mathsf{H}}

\newcommand{\D}[0]{\mathsf{D}}
\newcommand{\Q}[0]{\mathsf{Q}}
\newcommand{\E}[0]{\mathsf{E}}

\newcommand{\T}[0]{\mathsf{T}}
\newcommand{\F}[0]{\mathsf{F}}
\newcommand{\C}[0]{\mathsf{C}}
\newcommand{\A}[0]{\mathsf{A}}
\newcommand{\I}[0]{\mathsf{I}}

\newcommand{\K}[0]{\mathsf{K}}
\renewcommand{\S}[0]{\mathsf{S}}

\newcommand{\X}[0]{\mathsf{X}}

\renewcommand{\P}[0]{\mathsf{P}}

\newcommand{\Dxik}[0]{{\mathsf{D}}_{\bm{x}_{i} k}}

\newcommand{\Sxik}[0]{{\mathsf{S}}_{\bm{x}_{i}k}}

\newcommand{\Qxik}[0]{{\mathsf{Q}}_{\bm{x}_{i}k}}

\newcommand{\Exi}[0]{{\mathsf{E}}_{\bm{x}_{i}}}
\newcommand{\Exik}[0]{{\mathsf{E}}_{\bm{x}_{i}k}}

\newcommand{\Nxig}[0]{{\mathsf{N}}_{\bm{x}_{i}\gamma}}

\newcommand{\B}[0]{\mathsf{B}}
\newcommand{\R}[0]{\mathsf{R}}
\newcommand{\N}[0]{\mathsf{N}}
\newcommand{\Bg}[0]{\mathsf{B}_{\gamma}}

\newcommand{\Rgk}[0]{\mathsf{R}_{\gamma k}}
\newcommand{\Rgv}[0]{\mathsf{R}_{\gamma v}}

\newcommand{\V}[0]{\mathsf{V}}

\DeclareMathOperator{\mydiag}{diag}

\newcommand{\etal}[0]{{et~al.\@}\xspace}
\newcommand{\eg}[0]{{e.g.\@}\xspace}
\newcommand{\ie}[0]{{i.e.\@}\xspace}

\newcommand{\ignore}[1]{} 

\newcommand{\poly}[1]{\ensuremath{\mathbb{P}^{#1}}({\Omega}_k)}

\newcommand{\cont}[1]{\ensuremath{\mathcal{C}^{#1}}({\Omega}_k)}

\newcommand{\vecLtwo}[1]{[\ensuremath{L^2}({\Omega}_k)]^{#1}}

\newcommand{\IR}[1]{\mathbb{R}^{#1}}%
\newcommand{\IRtwo}[2]{\mathbb{R}^{{#1}\times{#2}}}%

\newcommand{\avg} [1]{\left\{  {#1}\right\}  }%

\global\long\def\dx#1#2{\frac{\partial#1}{\partial#2}}%
\global\long\def\dxx#1#2{\frac{\partial^{2}#1}{\partial#2^{2}}}%
\newcommand{\pder[2]}[]{\dfrac{\partial#1}{\partial#2}}%
\newcommand{\der[2]}[]{\dfrac{{\rm d}#1}{{\rm d}#2}}%
\global\long\def\norm#1{\left\vert \left\vert #1\right\vert \right\vert }%
\global\long\def\fn#1{\mathcal{#1}}%
\global\long\def\fnb#1{\bm{\mathcal{#1}}}%
\global\long\def\mathds#1{\mathds{#1}}%

\newcommand*\overlinewithlessheight[1]{{\mathpalette\overline@aux{#1}}}

\begin{document}
	
	
	\begin{frontmatter}
		
		\title{Entropy-split multidimensional summation-by-parts discretization of the Euler and \red{compressible} Navier-Stokes equations} 
	
	\author{Zelalem Arega Worku \corref{cor1}}
	\cortext[cor1]{Corresponding author: } 
	\ead{zelalem.worku@mail.utoronto.ca}
	
	\author{David W. Zingg}
	\ead{dwz@oddjob.utias.utoronto.ca}
	
	\address{Institute for Aerospace Studies, University of Toronto, Toronto, Ontario, M3H 5T6, Canada}
	
	
	\addcontentsline{toc}{section}{Abstract}
	\begin{abstract}
		High-order Hadamard-form entropy stable multidimensional summation-by-parts discretizations of the Euler and \red{compressible} Navier-Stokes equations are considerably more expensive than the standard divergence-form discretization. In search of a more efficient entropy stable scheme, we extend the entropy-split method for implementation on unstructured grids and investigate its properties. The main ingredients of the scheme are Harten's entropy functions, diagonal-$ \mathsf{E} $ summation-by-parts operators with diagonal norm matrix, and entropy conservative simultaneous approximation terms (SATs). We show that the scheme is high-order accurate and entropy conservative on periodic curvilinear unstructured grids for the Euler equations. An entropy stable matrix-type \blue{interface} dissipation operator is constructed, which can be added to the SATs to obtain an entropy stable semi-discretization. Fully-discrete entropy conservation is achieved using a relaxation Runge-Kutta method. Entropy stable viscous SATs, applicable to both the Hadamard-form and entropy-split schemes, are developed for the \red{compressible} Navier-Stokes equations. In the absence of heat fluxes, the entropy-split scheme is entropy stable for the \red{compressible} Navier-Stokes equations. Local conservation in the vicinity of discontinuities is enforced using an entropy stable hybrid scheme. Several numerical problems involving both smooth and discontinuous solutions are investigated to support the theoretical results. Computational cost comparison studies suggest that the entropy-split scheme offers substantial efficiency benefits relative to Hadamard-form multidimensional SBP-SAT discretizations. 		
	\end{abstract}
	
\end{frontmatter}

\section{Introduction}
In recent years, research interest in nonlinearly stable skew-symmetric discretizations of the governing equations of fluid dynamics has been revived \cite{sjogreen2019entropy,yee2023recent,sjogreen2021construction,nordstrom2021nonlinear,nordstrom2022skew}. These studies are partly driven by the desire to improve upon some of the theoretical and practical aspects of the  high-order Hadamard-form or two-point flux entropy conservative or dissipative scheme \cite{fisher2013high, fisher2013discretely,carpenter2014entropy,chen2017entropy,crean2018entropy,parsani2015entropy}. Despite its demonstrated robustness \cite{rojas2021robustness,parsani2021high,chan2022robust,chan2022entropy,parsani2021simulation,bergmann2020assessment}, the Hadamard-form entropy conservative scheme is so far incomplete as there are remaining pieces that are actively being researched such as positivity preservation of thermodynamic quantities \cite{yamaleev2022positivity,lin2023positivity} and entropy stable boundary conditions \cite{parsani2015entropywall,svard2018entropy,dalcin2019conservative,svard2021entropy,chan2022entropy}. The efficiency of the scheme has also been given due attention \cite{ranocha2021efficient,chan2019efficient}, as its computational cost is larger than that of the divergence-form implementation. Nonlinearly stable skew-symmetric split schemes may potentially allow bypassing some of these issues although at a cost of compromising other coveted properties such as mass conservation, kinetic energy preservation, and robustness. \red{The entropy-split method is an example of nonlinearly stable skew-symmetric schemes, which has been applied in the traditional finite difference fashion \cite{gerritsen1996designing,yee2000entropy,sjogreen2019entropy,yee2023recent}. Entropy splitting refers to a skew-symmetric splitting of the flux derivatives of the Euler equations into conservative and nonconservative portions using the symmetrizability and homogeneity properties of the equations \cite{yee2000entropy}. So far, the entropy-split method has not been generalized for discontinuous element-type implementations, which approximate the solution on a mesh consisting of non-overlapping elements supporting solution discontinuities at element interfaces.} The principal objective of this paper is to extend \violet{the} existing skew-symmetric entropy-split scheme to element-type implementations and shed light on the trade-offs between a few of the numerical properties afforded by it in comparison to the Hadamard-form implementation of the Euler and \red{compressible} Navier-Stokes equations on unstructured meshes. In particular, the study aims to provide some insights to the following questions: how robust is the entropy-split method compared to the Hadamard-form scheme? how much more efficient is it? how can its loss of conservation be mitigated?

In his seminal paper of 1983, Harten \cite{harten1983symmetric} proposed a family of entropy functions that symmetrizes the Euler and \red{compressible} Navier-Stokes equations without the heat fluxes. Furthermore, he showed that the fluxes in the Euler equations are homogeneous in the entropy variables corresponding to the proposed entropy functions.  Olsson and Oliger \cite{olsson1994energy} used these symmetry and homogeneity results to write the Euler equations in a skew-symmetric entropy-split form and obtained a nonlinear energy estimate under some assumptions. Subsequently, Gerittsen and Olsson \cite{gerritsen1996designing} designed a split high-order entropy-conserving scheme by discretizing the entropy-split form of the Euler equations using the classical summation-by-parts (SBP) finite difference operators of Strand \cite{strand1994summation}. Yee \etal \cite{yee2000entropy} investigated various aspects of the entropy-split scheme including the choice of the arbitrary splitting parameter, effects of the splitting on nonlinear stability, and application of characteristic-based filters. Later, Sv{\"a}rd and Mishra \cite{svard2012entropy} constructed entropy stable boundary conditions which extend the entropy stability of the entropy-split scheme to bounded domains. More recently, Sj{\"o}green and Yee \cite{sjogreen2019entropy} showed that the entropy-split scheme is entropy conservative for the Euler equations and demonstrated its robustness in a series of papers \cite{sjogreen2021construction,yee2023recent,sjogreen2018high}. For a more detailed review and recent developments of the entropy-split scheme, see \cite{yee2023recent} and the references therein. 

In this paper, we extend the entropy-split scheme to element-type implementations, which \red{facilitate flexible mesh generation around complex geometries, especially when applied with simplicial elements}. The main ingredients of the new scheme are Harten's entropy functions, diagonal-norm SBP operators \cite{hicken2016multidimensional,chen2017entropy} with collocated volume and surface nodes, and entropy conservative or dissipative simultaneous approximation terms (SATs) \cite{carpenter1994time}, which are coupling terms at element interfaces and boundaries. We show that the scheme is high-order accurate and entropy conservative on periodic general curvilinear unstructured grids for the Euler equations. Generalization of the scheme to bounded domains is not pursued in this article, but it is possible to make this extension using the entropy stable boundary conditions proposed by Sv{\"a}rd and Mishra \cite{svard2012entropy}.  If heat fluxes are neglected, then the entropy stability result also applies to the \red{compressible} Navier-Stokes equations provided that appropriate viscous SATs are used. To this end, we extend the viscous SATs developed in \cite{yan2018interior,worku2021simultaneous} to systems of equations with a positive semidefinite diffusivity tensor, alleviating the requirement to have a positive definite diffusivity tensor for stability. These viscous SATs are applicable to the entopy-split as well as Hadamard-form discretizations and do not require setting the viscous terms as a system of first-order equations unlike the widely used viscous SATs in the literature \cite{carpenter2014entropy,parsani2015entropy,gassner2018br1,chen2020review,chan2022entropy}.  

The entropy-split discretization is not conservative in the sense of Lax-Wendroff \cite{lax1960systems}. Hence, for problems with discontinuous solutions, we introduce a strategy to mitigate the effects of the loss of conservation such as incorrect shock speed. The main idea is to enforce local conservation by using the Hadamard-form discretization near discontinuities while ensuring entropy stability via a seamless coupling between the Hadamard-form and entropy-split schemes. This strategy allows one to benefit from the conservative property of the Hadamard-form discretization in critical areas of the computational domain and the efficiency of the entropy-split scheme elsewhere without compromising entropy stability. A similar strategy of switching between different schemes, but without enforcing entropy stability, has been proposed by other authors, \eg, see \cite{hou1994nonconservative,sjogreen2021construction}. Hou and LeFloch \cite{hou1994nonconservative} showed that such a strategy of enforcing local conservation guarantees convergence to the weak solution for scalar conservation laws. 

High-order dissipation operators have been designed \cite{ismail2009affordable,derigs2017novel,winters2017uniquely} to augment the Hadamard-form entropy conservative scheme with the physical property needed at shock waves. The dissipation operators are usually integrated with the inviscid SATs, producing an entropy dissipative variant of the entropy conservative scheme. In this work, we design a matrix-type interface dissipation operator for the entropy-split and Sj{\"o}green-Yee Hadamard-form \cite{sjogreen2019entropy} discretizations. Finally, the extension from semi-discrete to full entropy conservation or dissipation is achieved by using the fourth-order explicit relaxation Runge-Kutta (RRK4) method of Ranocha \etal \cite{ranocha2020relaxation}. Several benchmark problems involving both  smooth and discontinuous solutions are investigated to support our theoretical findings and compare various numerical properties of the entropy-split and Hadamard-form discretizations.  We also present efficiency comparisons between these two types of discretizations as applied to problems governed by the Euler and \red{compressible} Navier-Stokes equations. 

In the next section, we begin by introducing SBP operators and the notation used, followed by a review of the homogeneity and symmetrization of the Euler and \red{compressible} Navier-Stokes equations. In \cref{sec:continuous analysis}, we present the continuous entropy conservation and stability analysis. \cref{sec:discrete analysis} contains the discrete entropy conservation analysis for the Euler equations, the conservation error analysis of the entropy-split scheme, a strategy to handle discontinuous solutions, and the construction of \blue{interface} dissipation operators. In \cref{sec:nse}, extension of the entropy-split scheme to the \red{compressible} Navier-Stokes equations is discussed and some implementation details are presented in \cref{sec:implementation details}. Numerical results are presented in \cref{sec:results} followed by conclusions in \cref{sec:conclusions}.

\section{Preliminaries}\label{sec:preliminaries}

\subsection{SBP operators and notation}
Spatial discretizations are handled using element-type SBP operators, which play a crucial role in the entropy stability analysis. The spatial domain, $ \Omega $, is assumed to be compact and connected, and it is tessellated into $ n_e $ non-overlapping elements, $ {\mathcal T}_h \coloneqq \{\{ \Omega_k\}_{k=1}^{n_e}: \Omega=\cup_{k=1}^{n_e} {\Omega}_k\}$. The boundaries of each element are assumed to be piecewise smooth and will be referred to as facets or interfaces, and the union of the facets of element $ \Omega_{k} $ is denoted by $ \Gamma_k \coloneqq \partial\Omega_{k} $. The set of all interior interfaces is denoted by $\Gamma ^I  \coloneqq \{\Gamma_k \cap \Gamma_v : k,v=1,\dots,n_e, k\neq v \}$. The set of facets of $ \Omega_k $ that are also interior facets is denoted by $ \Gamma^I_k  \coloneqq \Gamma^I\cap\Gamma_k $, and $ \Gamma  \coloneqq \{\partial\Omega \cap \Gamma_k : k = 1,\dots, n_e\}$ delineates the set of all boundary facets. The set of $ n_p $ volume nodes in element $ \Omega_k $ is represented by $ S_k=\{\bm{x}_{j}\}_{j=1}^{n_p} $, and $ n_f $ denotes the number of nodes on facet $ \gamma \in \Gamma_{k}$. Scalar functions are written in uppercase script type, \eg, $\fnc{U}_k \in \cont{\infty}$, and vector-valued functions of dimension $ n $ are represented by boldface uppercase script letters, \eg, $\vecfnc{W}_k \in \vecLtwo{n}$. The space of polynomials of total degree $ p $ is denoted by $\poly{p} $. Vectors containing function values at a set of grid points are denoted by bold letters, \eg, $ \uk \in \IR{n_p}$. We define $h \coloneqq \max_{a, b \in S_k} \norm{a - b}_2$ as the nominal element size. Matrices are denoted by sans-serif uppercase letters, \eg, $\V \in \IRtwo{n_p}{n_p}$; $ \bm{1} $ denotes a vector consisting of all ones, $ \bm {0} $ denotes a vector or matrix consisting of all zeros. The sizes of  $ \bm{1} $ and $ \bm {0} $ should be clear from context. Finally, $ \I_{n} $ represents the identity matrix of size $ n \times n $. 

SBP operators on a physical element, $ \Omega_{k} $, are defined as follows \cite{hicken2016multidimensional}.
\begin{definition} \label{def:sbp}
	The matrix $\Dxik\in \IRtwo{n_p}{n_p}$ is a degree $ p $ SBP operator in the $ i $-direction approximating the first derivative $ \pdv{\bm{x}_i} $ on the set of nodes $ {S}_{k}=\{\bm{x}_{j}\}_{j=1}^{n_p} $ if
	\begin{enumerate}
		\item $[ \Dxik \bm{p}]_j = {\pdv{\fnc{P}}{{x}_i}}({\bm{x}_j}) $  for all $\fnc{P} \in \poly{p} $
		\item $ \Dxik={\H}^{-1} _{k}\Qxik $ where $ {\H}_{k}$ is a symmetric positive definite (SPD) matrix, and 
		\item $ \Qxik = \Sxik + \frac{1}{2} \Exik$ where $ \Sxik= - \Sxik^T $, $ \Exik = \Exik^T $ and $ \Exik $ satisfies
		\[ \bm{p}^T \Exik  \bm {q} = \int_{{\Gamma}_{k}} \fnc{P}\fnc{Q} \;{n}_{{x}_i} \dd{\Gamma}\]
		$\forall \fnc{P},\fnc{Q} \in \poly{r} $, where $ r \ge p $, and $ {n}_{{x}_i} $ is the $ i$-component of the outward pointing unit normal vector on ${\Gamma}_{k} $.
	\end{enumerate}  
\end{definition} 

The third property in \cref{def:sbp} implies that $ \Q_{{x}_{i}k} + \Q_{{x}_{i}k}^T = \E_{{x}_{i}k}$, which is referred to as the SBP property. The SBP operators in \cref{def:sbp} contain grid metric terms as they are defined for the physical element $ \Omega_{k} $, which may be curved. We make the following assumption regarding the geometric mapping from the reference to physical elements. 
\begin{assumption}\label{ass:mapping}
	We assume that there is a bijective and time-independent polynomial mapping of degree $ p_{\text{geom}}$ from the reference to the physical elements. Furthermore, we assume that the degree of the mapping is limited to $ p_{\text{geom}} \le p+1 $ and $ p_{\text{geom}}\le \lfloor p/2 \rfloor +1 $ in two and three dimensions, respectively, where $ p $ is the degree of the SBP operator.
\end{assumption}
\cref{ass:mapping} ensures that the discrete metric identities are satisfied on curved elements \cite{crean2018entropy,shadpey2020entropy}. If a higher degree polynomial mapping is needed to represent a three-dimensional geometry, then one may employ the strategy outlined in  \cite{crean2018entropy}, which involves solving elementwise quadratic optimization problems. For construction of multidimensional SBP operators, we refer the reader to \cite{hicken2016multidimensional,fernandez2018simultaneous,chen2017entropy}, and see, \eg, \cite{crean2018entropy,shadpey2020entropy,worku2021simultaneous} for their implementation on curved elements.  A more general review of SBP-SAT discretizations can be found in \cite{fernandez2014review,svard2014review}. 

In this work, we exclusively consider SBP operators with diagonal $ \H $ and $ \E $ matrices. In the one-dimensional case, we use the Legendre-Gauss-Lobatto (LGL) operators \cite{gassner2013skew,fernandez2014generalized,carpenter2014entropy,carpenter2015entropy} and higher-dimensional cases use the SBP diagonal-$ \E $ multidimensional operators \cite{chen2017entropy,hicken2016multidimensional}. The diagonal of the norm matrix, $ \H $, contains strictly positive entries which are weights of a quadrature rule of at least degree $ 2p-1 $. Hence, it approximate an inner product,
\begin{equation}
	\bm{p}_k^T \H_k \bm{q}_k = \int_{\Omega_k} \fn{P}\fn{Q} \dd{\Omega}, \quad \forall \;\fn{P}\fn{Q} \in \poly{r}, \quad r\le 2p-1,
\end{equation}
and is used to define the $ \H $-norm, $ \norm{\bm{u}}_{\H}^2 = \bm{u}^T \H \bm{u}$. The facet nodes of diagonal-$ \E $ SBP operators are collocated with the volume nodes. This implies that an extrapolation operator, $ \Rgk \in \IRtwo{n_f}{n_p}$, that extrapolates the solution at the volume nodes to the facet nodes will simply produce the solution values at the facet nodes. Assuming that a quadrature rule of degree at least $ 2p $ with strictly positive weights exists on the nodes of facet $ \gamma \in \Gamma_{k} $ with its weights contained in the diagonal matrix $ \B_\gamma \in \IRtwo{n_f}{n_f} $, the matrix $ \Exik $ can be decomposed as \cite{fernandez2018simultaneous}
\begin{equation} \label{eq: Exik}
	\Exik = \sum_{\gamma\in\Gamma_{k}} \Exi^{kk}, \quad \text{where} \;\; \Exi^{k k}=\Rgk^T \Bg \Nxig \Rgk,
\end{equation}
where $ \Nxig \in \IRtwo{n_f}{n_f}$ is a diagonal matrix containing the $ i $-component of the outward unit normal vector on facet $ \gamma $. The facet nodes on adjoining facets of $ \Omega_{k} $ and its neighboring element, $ \Omega_{v} $, coincide, \ie, only conforming meshes are considered. The coupling terms between $ \Omega_{k} $ and $ \Omega_{v} $ make use of the following matrices and their relation, 
\begin{equation}
	\begin{aligned}
		\Exi^{kv}=\Rgk^T \Bg \Nxig \Rgv, \quad && 
		\Exi^{vk}=-\Rgv^T \Bg \Nxig \Rgk, \quad &&
		\Exi^{kv}= - (\Exi^{vk})^T.
	\end{aligned}
\end{equation}
For systems of equations with $ n_c $ components, the operators defined so far are applied using a Kronecker product with the identity matrix, \eg, $ \overline{\D}_{\bm{x}_i k} =\Dxik \otimes\I_{n_{c}} $, where the overbar indicates application to systems of equations.

\subsection{Governing equations}
The $ d $-dimensional \red{compressible} Navier-Stokes equations can be written as,
\begin{equation} \label{eq:nse}
	\begin{aligned}
		\pder[\fnb U]t+\sum_{i=1}^{d}\pder{{x}_{i}}\fnb F^{I(i)} & =\sum_{i=1}^{d}\pder{{x}_{i}}\fnb F^{V(i)},\quad \bm{x}=[\red{x_{1}, x_{2}, x_{3}}]^T\in\Omega,\quad t\in[0,t_{f}],\\
		\fn{B}\fnb U\left(\bm{x},t\right)  & =\fnb{U}_{\Gamma},\quad \bm{x}\in\Gamma,\\
		\fnb U\left(\bm{x},0\right) & =\fnb{U}_{0},\quad \bm{x}\in\Omega\red{,}
	\end{aligned}
\end{equation}
where $ \fnc{U}(\bm{x},t)$ contains the conservative variables, $ \fnb F^{I(i)}(\fnb U) $ and $  \fnb F^{V(i)}(\fnc{U}, \dx{\fnb U}{\bm{x}_j}) $ are the inviscid (Euler) and viscous fluxes in the $ i$-direction, respectively, and $ \fn{B} $ is a boundary operator. The Euler equations are obtained by nullifying the right-hand side (RHS) of the first line in \cref{eq:nse}. Assuming a Newtonian fluid, the conservative variables, the inviscid, and the viscous fluxes in three dimensions are given by
\begin{equation}
	\begin{aligned}
		\fnb U=\begin{bmatrix}\rho\\
			\rho \red{V_{1}}\\
			\rho \red{V_{2}}\\
			\rho \red{V_{3}}\\
			e
		\end{bmatrix},	
		\;\;
		\fnb F^{I(i)}=\begin{bmatrix}\rho V_{i}\\
			\rho V_{i}V_{1}+\delta_{i1}p\\
			\rho V_{i}V_{2}+\delta_{i2}p\\
			\rho V_{i}V_{3}+\delta_{i3}p\\
			(e+p)V_{i}
		\end{bmatrix},	
		\;\;
		\fnb F^{V(i)}&=\begin{bmatrix}0\\
			\tau_{{x}_{i}\red{x_{1}}}\\
			\tau_{{x}_{i}\red{x_{2}}}\\
			\tau_{{x}_{i}\red{x_{3}}}\\
			\Sigma_{j=1}^{d}V_{j}\tau_{{x}_{i}{x}_{j}}
		\end{bmatrix}+\begin{bmatrix}0\\
			0\\
			0\\
			0\\
			\kappa\dx{T}{{x}_{i}}
		\end{bmatrix},
		\;\;
		\tau_{{x}_{i}{x}_{j}}=\mu\left(\pder[V_{i}]{{x}_{j}}+\pder[V_{j}]{{x}_{i}}-\frac{2}{3}\delta_{ij}\sum_{k=1}^{d}\pder[V_{k}]{{x}_{k}}\right)\red{,}
	\end{aligned}
\end{equation}
where $ \rho $ is the density, $\red{ \bm{V}=[V_{1},V_{2},V_{3}]^{T}} $ contains the $ \bm{x}=[\red{x_{1},x_{2},x_{3}}] ^{T}$ velocity components, $ p $ is the pressure, $ e $ is the total energy, $ T $ is the temperature, $ \mu $ is the dynamic viscosity, $ \kappa $ is the heat conductivity, and $ \delta_{ij} $ is the Kronecker delta. Assuming thermally and calorically perfect gases, the system is closed by the relation $ e={p}/{(\gamma-1)}+{\rho}\bm{V}^{T}\bm{V}/2 $, where $ \gamma=7/5$ is the ratio of specific heats. The viscous fluxes are written as a sum of the viscous and heat fluxes, which is a suitable form for later discussion.
\subsection{Homogeneity}
\begin{definition}\label{def:homogeneous}
	$ \fnb F\left(\fnb U,\fnb V\right)  $ is a degree $ n $ homogeneous function in the argument $ \fnb U $ if $ \fnb F\left(\eta\fnb U,\fnb V\right)=\eta^{n}\fnb F\left(\fnb U,\fnb V\right) $, $ \eta\neq 0 $.
\end{definition}
\red{The definition of homogeneous functions can be found in several books, \eg, \cite{gokcen1996thermodynamics}, and the following lemma is a restatement of Euler's homogeneous function theorem, \eg, see \cite{lomax2013fundamentals}, in suitable forms for later use.} The proof of \red{the lemma} is omitted as it follows from elementary application of calculus and \cref{def:homogeneous}. 
\begin{lemma}\label{lem:homogeniety}
	Consider a function, $  \fnb{F} = \fnb F\left(\fnb U,\fnb V\right) $, that is degree $ n $ homogeneous in its first argument, and let $  \A=\dx{\fnb F}{\fnb U} $. Then, the the following identity holds,
	\begin{align}
		\A\fnb U&=n\fnb F, \label{eq:eq1 prop homogeniety}
	\end{align}
	Furthermore, if $  \fnb{F} = \fnb F\left(\fnb U\right)  $, then, in addition to \cref{eq:eq1 prop homogeniety}, it can be shown that
	\begin{align}
		\pder x\left(\A\fnb U\right)	&=n\pder[\fnb F]x=n\A\pder[\fnb U]x, \label{eq:eq2 prop homogeniety}
	\end{align}
	or equivalently,
	\begin{equation} \label{eq:eq3 prop homogeniety}
		\dx{\A}{x}\fnb U=(n-1)\A\dx{\fnb U}{x}.
	\end{equation}
\end{lemma}

It is well known that the inviscid fluxes are homogeneous functions of degree one in the conservative variables, \eg, see \cite{lomax2013fundamentals}. Similarly, the viscous fluxes are homogeneous functions of degree one in the argument $ \dx{\fnb U}{{x}_{j}} $ as well as $ \dx{\fnb Q}{{x}_{j}} $, where $ \fnb Q =  [\rho, \red{V_{1}, V_{2}, V_{3}}, T]^{T} $ is a vector of the primitive variables. It can also be shown that the viscous fluxes are homogeneous functions of degree $ -1 $ in the argument $ \fnb U $ but not homogeneous in $ \fnb Q $. Hence, \cref{lem:homogeniety} implies that 
\begin{align}
	\fnb F^{I(i)}	&=\A_{i}\fnb U,	 && \A_{i}	\coloneqq\pder[\fnb F^{I(i)}]{\fnb U}, &&\pder[\fnb F^{I(i)}]{{x}_{i}}=\A_{i}\pder[\fnb U]{{x}_{i}}, \label{eq:inviscid homo U}\\
	\fnb F^{V(i)} &=\sum_{j=1}^{d}\K_{ij}\pder[\fnb U]{x_{j}} =\sum_{j=1}^{d} \C_{ij}\pder[\fnb Q]{x_{j}},	&& \K_{ij}	\coloneqq\pder[\fnb F^{V(i)}]{\left(\dx{\fnb U}{{x}_{j}}\right)},  &&  \C_{ij}	\coloneqq\pder[\fnb  F^{V(i)}]{\left(\dx{\fnb Q}{{x}_{j}}\right)}. \label{eq:viscous homo dU}
\end{align}
The \red{compressible} Navier-Stokes equations, \cref{eq:nse}, can therefore be written as \cite{hughes1986new}
\begin{equation}\label{eq:nse nonconservative}
	\pder[\fnb U]t+\sum_{i=1}^{d}\A_{i}\pder[\fnb U]{{x}_{i}}=\sum_{i=1}^{d}\pder{{x}_{i}}\left(\sum_{j=1}^{d}\K_{ij}\pder[\fnb U]{{x}_{j}}\right).
\end{equation}
\subsection{Symmetrization}
A stability analysis via the energy method requires symmetrization of the inviscid flux Jacobians, $ \A_{i} $, and the diffusivity tensor, \ie, the block matrix containing the viscous flux Jacobians, $ \K=[\K_{ij}] $. One way to achieve this is through a change of variables. For hyperbolic systems, symmetrization of the flux Jacobians under a change of variables is intimately linked to the existence of an entropy function and the corresponding entropy fluxes; in fact, one implies the other \cite{harten1983symmetric}. However, not all entropy functions that symmetrize the inviscid flux Jacobians also symmetrize the diffusivity tensor. Entropy functions that enable symmetrization of both via a change of variables are referred to as generalized entropy functions \cite{hauke2006thermodynamics}. 

The Euler equations are accompanied by a scalar entropy balance equation that describes the second law of thermodynamics; \ie, the physical entropy, $ s = \ln(p\rho^{-\gamma})$, in a system does not decrease, or mathematically,
\begin{equation}\label{eq:entropy balance}
	\der[]t\int_{\Omega}\fn S \dd{\Omega}+\sum_{i=1}^{d}\int_{\Gamma}\fn G^{(i)}n_{\bm{x}_i}\dd{\Gamma}\le0,
\end{equation}
where the mathematical entropy function and entropy flux are given, respectively, by 
\begin{align}
	\fn S &=-\rho \fn{H}(s)/(\gamma-1) \label{eq:mathematical entropy},\\
	\fn G^{(i)}&=V_{i}\fn{S}, \label{eq:entropy flux}
\end{align}
and $ \fn{H} $ is any differentiable function such that $ \dxx{\fn{H}}{\fn S} \big/ \dx{\fn{H}}{\fn{S}} < 1/\gamma $ \cite{harten1983symmetric}. Note that due to the negative scaling in the definition of $ \fn{S} $, the inequality in \cref{eq:entropy balance} states that the mathematical entropy in a system does not increase. Furthermore, the inequality becomes an equality if the solution is smooth. The entropy inequality, \cref{eq:entropy balance}, also applies to the \red{compressible} Navier-Stokes equations, but with some restrictions on $ \fn{S} $. Hughes \etal \cite{hughes1986new} proved that the only family of generalized entropy functions for the \red{compressible} Navier-Stokes equations are those that are affine transformations of $ s $, \ie, $ \fn{H}(s) $ must be linear in $ s $. Usually, the generalized entropy function
\begin{equation}\label{eq:generalized entropy}
	\fn{S}_{G} = -\rho s/(\gamma-1)
\end{equation}
is used for the Hadamard-form entropy stable discretizations of the \red{compressible} Navier-Stokes equations. Entropy-split schemes, however, use a family of entropy functions that was introduced by Harten \cite{harten1983symmetric},
\begin{equation}\label{eq:harten entropy}
	\fn S_{H}=-\frac{\gamma+\alpha}{\gamma-1}\rho \exp({\frac{s}{\alpha+\gamma}})=\beta \rho\left({p}{\rho^{-\gamma}}\right)^{\frac{1}{\alpha+\gamma}},
\end{equation}
where $ \beta = -(\gamma+\alpha)/(\gamma-1) $, and $ \alpha $ is an arbitrary parameter satisfying the condition $ \alpha>0 $ or $ \alpha<-\gamma $. It is possible to apply Harten's entropy functions to construct an Hadamard-form discretization that is entropy conservative or stable for the Euler equations \cite{sjogreen2019entropy}.

For the symmetrization of \cref{eq:nse nonconservative}, we first introduce the entropy variables defined as
\begin{equation}\label{eq:entropy vars}
	\fnb W^{T}\coloneqq \pder[\fn S]{\fnb U}.
\end{equation}
Applying a change of variables, \cref{eq:nse nonconservative} can be rewritten as
\begin{equation}\label{eq:nse symmetric}
	\tilde{\A}_{0}\pder[\fnb W]t+\sum_{i=1}^{d}\tilde{\A}_{i}\pder[\fnb W]{{x}_{i}}=\sum_{i=1}^{d}\pder{{x}_{i}}\left(\sum_{j=1}^{d}\tilde{\K}_{ij}\pder[\fnb W]{{x}_{j}}\right),
\end{equation}
where 
\begin{equation} \label{eq:symmetric coefficients}
	\begin{aligned}
		\tilde{\A}_{0}=\dx{\fnb U}{\fnb W}, && \quad \tilde{\A}_{i}=\A_{i}\tilde{\A}_{0},  &&  \quad \tilde{\K}_{ij}=\K_{ij}\tilde{\A}_{0}. 
	\end{aligned}
\end{equation}
If $ \fn{S} = \fn{S}_{G}$, then $ \tilde{\A}_{0} $, $ \tilde{\A}_{i} $, and $ \tilde{\K} = [\tilde{\K}_{ij}]$ will be symmetric simultaneously, and we have $ \tilde{\A}_{0} \succ 0$ and $ \tilde{\K}\succeq 0 $, where $ \succ $ and $ \succeq $ denote symmetric positive definiteness and symmetric positive semidefiniteness, respectively. If, instead, $ \fn{S} = \fn{S}_{H} $, then $ \tilde{\A}_0 \succ 0$ and $ \tilde{\A}_{i} = \tilde{\A}_{i}^T $, but $ \tilde{\K} \not\succeq 0$ unless the heat fluxes in the \red{compressible} Navier-Stokes equations are neglected. The positive definiteness of $ \tilde{\A}_{0} $ implies that its inverse exists and that there is a bijective mapping between the conservative and entropy variables. Furthermore, the fact that $ \tilde{\A}_{0} $ and $ \tilde{\A}_{i} $ are symmetric implies that $ \fnb U $ and $ \fnb F^{I(i)} $ are Jacobians of scalar functions with respect to the entropy variables, \ie, 
\begin{align}
	\fnb U^{T}	=\pder[\varphi]{\fnb W},\\
	\left[\fnb F^{I(i)}\right]^{T}	=\pder[\psi^{(i)}]{\fnb W},
\end{align}
where the scalar functions $ \varphi $ and $ \psi^{(i)}$ are called the potential and potential flux functions, respectively, and they satisfy the following identities \cite{godunov1961interesting,mock1980systems,harten1983symmetric},
\begin{align}
	\varphi&=\fnb W^{T}\fnb U-\fn S, \label{eq:phi=wTu - s}\\
	\psi^{(i)}&=\fnb W^{T}\fnb F^{I(i)}-\fn G^{(i)}. \label{eq:psi=wTf - g}
\end{align}

\section{Continuous entropy stability analysis }\label{sec:continuous analysis}
Olsson and Oliger \cite{olsson1994energy} proposed a splitting of the symmetrized Euler fluxes that allows one to exploit the homogeneity of the fluxes to show entropy conservation. At the heart of the splitting is Harten's family of entropy functions, which, unlike the symmetrization using $ \fn{S}_{G}$, transfers the homogeneity of the Euler equations to the symmetrized system. \red{For smooth solutions,} the entropy-split Euler equations are written as 
\begin{equation}\label{eq:split Euler}
	\pder[\fnb U]t+\frac{\beta}{\beta+1}\sum_{i=1}^{d}\pder[\fnb F^{I(i)}]{{x}_{i}}+\frac{1}{\beta+1}\sum_{i=1}^{d}\tilde{\A}_{i}\pder[\fnb W]{{x}_{i}}=\bm{0},
\end{equation}
where the relative amount of conservative and nonconservative portions of the splitting is controlled by the arbitrary parameter $ \beta $, and the entropy variables are given by
\begin{equation}\label{eq:harten entropy variables}
	\fnb W^{T} =  \pder[\fn{S}_{H}]{\fnb U} =  \frac{\rho}{p}(p\rho^{-\gamma})^{\frac{1}{\alpha+\gamma}}\left[\begin{array}{ccccc}
		-\frac{\alpha}{\gamma-1}\frac{p}{\rho}-\frac{1}{2}\bm{V}^T\bm{V}, & \bm{V}^T, & -1\end{array}\right].
\end{equation}

\begin{lemma}\label{lem:homo in W}
	The Euler fluxes, $ \fnb F^{I(i)}\left(\fnb W\right)$, are homogeneous functions of degree $ \beta $ with respect to Harten's entropy variables, \cref{eq:harten entropy variables}, \ie, $\fnb F^{I(i)}(\eta \fnb{W}) = \eta^{\beta} \fnb{F}^{I(i)}(\fnb W)   $, $ \eta\neq 0 $.
\end{lemma}
\begin{proof}
	See, \eg, \cite{harten1983symmetric,sjogreen2019entropy}.
\end{proof}
The consequence of \cref{lem:homo in W} in combination with \cref{lem:homogeniety} is that 
\begin{equation}\label{eq:homo in W}
	\tilde{\A}_{i}\fnb W=\beta\fnb F^{I(i)},
\end{equation}
which is a crucial relation to show that the entropy-split Euler equations conserve entropy for smooth solutions. Multiplying \cref{eq:split Euler} by $ \fnb{W}^T $ from the left and integrating, we find
\begin{equation}
	\int_{\Omega}\fnb W^{T}\pder[\fnb U]t\dd{\Omega} + \frac{\beta}{\beta+1}\sum_{i=1}^{d}\int_{\Omega}\fnb W^{T}\pder[\fnb F^{I(i)}]{{x}_{i}}\dd{\Omega}+\frac{1}{\beta+1}\sum_{i=1}^{d}\int_{\Omega}\fnb W^{T}\tilde{\A}_{i}\pder[\fnb W]{{x}_{i}}\dd{\Omega} = 0.
\end{equation}
Using $ \fnb W^{T}\dx{\fnb U}{t}=\dx{\fn S}{t}$, integrating by parts the second term, and simplifying, we have 
\begin{equation}
	\int_{\Omega}\pder[\fn S]t \dd{\Omega}-\frac{\beta}{\beta+1}\sum_{i=1}^{d}\int_{\Omega}\pder[\fnb W^{T}]{{x}_{i}}\fnb F^{I(i)}\dd{\Omega}+\frac{1}{\beta+1}\sum_{i=1}^{d}\int_{\Omega}\fnb W^{T}\tilde{\A}_{i}\pder[\fnb W]{{x}_{i}}\dd{\Omega}+\frac{\beta}{\beta+1}\sum_{i=1}^{d}\int_{\Gamma}\fnb W^{T}\fnb F^{I(i)}n_{{x}_{i}}\dd{\Gamma}=0,
\end{equation}
which after making use of the symmetry of $ \tilde{\A}_{i} $ and \cref{eq:homo in W} yields
\begin{equation}\label{eq:entropy conservation 0}
	\der[]t\int_{\Omega}\fn S \dd{\Omega}+\frac{\beta}{\beta+1}\sum_{i=1}^{d}\int_{\Gamma}\fnb W^{T}\fnb F^{I(i)}n_{{x}_{i}}\dd{\Gamma}=0.
\end{equation}
In \cite{sjogreen2019entropy}, it was found, through a direct symbolic computation, that
\begin{equation}\label{eq:WTF=G}
	\frac{\beta}{\beta+1}\fnb W^{T}\fnb F^{I(i)}=\fn G^{(i)}.
\end{equation} 
Substituting \cref{eq:WTF=G} into \cref{eq:entropy conservation 0} gives the entropy balance stated in \cref{eq:entropy balance} \red{for smooth solutions}, 
\begin{equation}\label{eq:entropy conservation Euler}
	\der[]t \int_{\Omega} \fn S \dd{\Omega} + \sum_{i=1}^{d}\int_{\Gamma}\fn{G}^{(i)} n_{{x}_i}\dd{\Gamma} = 0,
\end{equation}
and in the case of periodic boundary conditions, we find $ \dv{}{t} \int_{\Omega}\fn{S}=0$.
\begin{lemma}\label{lem:potential and entropy relation}
	For Harten's family of entropy functions, the following relations hold
	\begin{align}
		\fn G^{(i)}&=\beta\psi^{(i)}, \label{eq:G=beta psi} \\
		\fn S&=\beta\varphi. \label{eq:S=beta phi}		
	\end{align}
\end{lemma}
\begin{proof}
	The first identity follows from substitution of \cref{eq:WTF=G} into \cref{eq:psi=wTf - g} and simplification. The second identity, \cref{eq:S=beta phi}, is obtained by direct computation of the vector product appearing in \cref{eq:phi=wTu - s}.
\end{proof}

\begin{remark}
	The identities in \cref{lem:potential and entropy relation} imply that for smooth problems governed by the Euler equations, the potential function obtained using Harten's entropy functions, $ \varphi $, is conserved. The identity in \cref{eq:G=beta psi} is crucial in the design of an entropy-split scheme for element-type discretizations.
\end{remark}

\violet{When Harten's entropy function, \cref{eq:harten entropy}, is used, the entropy balance given in \cref{eq:entropy balance} holds for the \red{compressible} Navier-Stokes equations if the heat fluxes are neglected.} In this case, $ \tilde{\K} $ is symmetric positive semidefinite \cite{harten1983symmetric}. We consider the RHS of \cref{eq:nse symmetric}; entropy conservation of the LHS is already shown for the entropy-split Euler equations in \cref{eq:entropy conservation Euler}. Multiplying the viscous terms in \cref{eq:nse symmetric} by the entropy variables, integrating, and applying integration by parts gives 
\begin{equation} \label{eq:entropy balance viscous}
	\int_{\Omega}\fnb{W}^T \sum_{i=1}^{d}\pder{{x}_{i}}\left(\sum_{j=1}^{d}\tilde{\K}_{ij}\pder[\fnb W]{{x}_{j}}\right) \dd{\Omega} = -\int_{\Omega}\pder[\fnb{W}]{\bm{x}}^T\tilde{\K}\pder[\fnb W]{\bm{x}}\dd{\Omega} + \sum_{i=1}^{d} \sum_{j=1}^{d} \int_{\Gamma}  \fnb{W}^T\tilde{\K}_{ij}\pder[\fnb W]{{x}_{j}} \dd{\Gamma},
\end{equation}
where $  \dx{\fnb{W}}{\bm{x}}^{T} = [\dx{\fnb{W}}{{x}_1}^{T},\dots,\dx{\fnb{W}}{{x}_d}^{T}]$. Since $ \tilde{\K}\succeq 0 $, the first term on the RHS of \cref{eq:entropy balance viscous} is negative, and if appropriate boundary conditions are used, \eg, periodic, then we have 
\begin{equation}\label{eq:entropy balance viscous 1}
	\int_{\Omega}\fnb{W}^T \sum_{i=1}^{d}\pder{{x}_{i}}\left(\sum_{j=1}^{d}\tilde{\K}_{ij}\pder[\fnb W]{{x}_{j}}\right) \dd{\Omega} \le 0,
\end{equation}
which combined with \cref{eq:entropy conservation Euler} yields the entropy balance in \cref{eq:entropy balance}.

\begin{remark}\label{rem:viscous flux with S_G}
	If the generalized entropy function \cref{eq:generalized entropy} is used for the symmetrization of the \red{compressible} Navier-Stokes equations, then $ \tilde{\K}\succeq 0 $ with the heat fluxes included, and the inviscid and viscous fluxes satisfy \cref{eq:entropy conservation Euler} and \cref{eq:entropy balance viscous 1}, respectively. The analysis for the viscous terms is the same as the above, while the inviscid analysis uses the identity \cite{harten1983symmetric} 
	\begin{equation}\label{eq:WTFu=Gu}
		\fnb W^{T}\pder[\fnb F^{I(i)}]{\fnb{U}}=\pder[\fn G^{(i)}]{\fnb{U}}
	\end{equation}
	to arrive at entropy conservation of the un-split form of the equations,
	\begin{align}
		\int_{\Omega}\fnb W^{T}\left(\pder[\fnb U]t+\sum_{i=1}^{d}\pder[\fnb F^{I(i)}]{{x}_{i}}\right)\dd{\Omega}&=\int_{\Omega}\pder[\fn S]t+\sum_{i=1}^{d}\fnb W^{T}\pder[\fnb F^{I(i)}]{\fnb U}\pder[\fnb U]{{x}_{i}}\dd{\Omega} \nonumber
		\\&=\der[]t\int_{\Omega}\fn S \dd{\Omega}+\sum_{i=1}^{d}\int_{\Gamma}\fn G^{(i)}n_{{x}_{i}}\dd{\Gamma}=0,
	\end{align}
	where the divergence theorem is applied to obtain the facet integral in the last equality \red{and smoothness of the solution is assumed}.
\end{remark}

If Harten's entropy functions are used, the explicit entries of $ \tilde{\A}_{0} $, $ \tilde{\A}_{i} $, and $ \tilde{\K} $ without the heat fluxes can be found in \cite{harten1983symmetric}. If the heat fluxes are considered, which are denoted by $ \fnb{F}^{V(i)}_{h} $, then the $ n_c\times n_c $ matrix in the following will have to be added to the diagonal blocks of the diffusivity tensor, $ \tilde{\K}_{ii} $, obtained without considering the heat fluxes,
\begin{equation}\label{eq:heat Kij}
	\pder[\fnb F_{h}^{V(i)}]{{x}_{i}}=\pder{{x}_{i}}\left(\pder[\fnb F_{T}^{V(i)}]{{x}_{i}}\right)=\pder{{x}_{i}}\pder[\fnb F_{T}^{V(i)}]{\fnb W}\pder[\fnb W]{{x}_{i}}=\pder{{x}_{i}}\left[\begin{array}{cccc}
		0 & \dots &  & 0\\
		\vdots & \ddots\\
		\\
		0 & \dots &  & 0\\
		\kappa\pder[T]{\fnb W_{1}} & \kappa\pder[T]{\fnb W_{2}} & \dots & \kappa\pder[T]{\fnb W_{n_{c}}}
	\end{array}\right]\pder[\fnb W]{{x}_{i}},
\end{equation}
where $ (\fnb F_{T}^{V(i)})^{T}=[0,\dots,\kappa T] $ is a vector of size $ 1\times n_{c} $, and $ T $, in terms of the entropy variables, is given by 
\begin{equation}\label{eq:T in W}
	T=\frac{\left(\gamma-1\right)\left[2\fnb W_{1}\fnb W_{n_{c}}-\Sigma_{i=2}^{n_{c}-1}\fnb W_{i}^{2}\right]}{2\alpha R\fn{\fnb W}_{n_{c}}^{2}}.
\end{equation}
Since $ T $ is a function of all the entries of $ \fnb{W} $, all entries in the last row of $ \dx{\fnb F_{T}^{V(i)}}{\fnb W} $ are nonzero, which destroys the symmetry and positive semidefiniteness of $ \tilde{\K}_{ii} $.The gradients of the temperature with respect to the entropy variables can be simplified to give 
\begin{align}\label{eq:grad T by W}
	\pder[T]{\fnb W_{1}}&=\frac{\gamma-1}{\alpha R\fnb W_{n_{c}}},&\pder[T]{\fnb W_{i}}&=-\frac{\left(\gamma-1\right)V_{i}}{\alpha R\fnb W_{n_{c}}},&\pder[T]{\fnb W_{n_{c}}}&=-\frac{T}{\fnb W_{n_{c}}}+\frac{\left(\gamma-1\right)\bm{V}^{T}\bm{V}}{2R\fnb W_{n_{c}}},\quad i=\left\{ 2,\dots,n_{c}-1\right\} .
\end{align}
Furthermore, $ \fnb{W}_{n_c} $ can be written as 
\begin{equation}\label{eq:W5}
	\fnb W_{n_{c}}=-\frac{\rho}{p}\left(p\rho^{-\gamma}\right)^{\frac{1}{\alpha+\gamma}}=-p^{\frac{1-\alpha-\gamma}{\alpha+\gamma}}\rho^{\frac{\alpha}{\alpha+\gamma}}=-\left(RT\right)^{\frac{1-\alpha-\gamma}{\alpha+\gamma}}\rho^{\frac{1-\gamma}{\alpha+\gamma}}.
\end{equation}
Assuming $ \alpha<-\gamma $ is used, which corresponds to the physically relevant branch of the split parameter $ \beta $ \cite{yee2000entropy}, and that the density and pressure remain positive, we can see from \cref{eq:W5} that the magnitude of $ \fnb{W}_{n_c} $ increases as $ \rho $ increases and when $ p $ or $ T $ decreases. Hence, \cref{eq:grad T by W} suggests that for high density, low pressure or temperature, and low velocity flows, the gradients of the temperature in $ \dx{\fnb F_{T}^{V(i)}}{\fnb W} $ are close to zero. Hence, the \red{impact of} the heat fluxes can be minimal in such cases. While it is unclear at what exact conditions the estimate from the viscous terms will bound the estimate from the heat fluxes, it is reasonable to expect that the entropy-split method will produce a robust discretization for many practical problems governed by the \red{compressible} Navier-Stokes equations since \red{heat fluxes are rarely a source of difficulty in numerical approximations of the equations.}

If the generalized entropy function, \cref{eq:generalized entropy}, is used, then only the last entry in the last row of $ \dx{\fnb F_{T}^{V(i)}}{\fnb W} $ is nonzero as shown in \cite{hughes1986new}, and $ \tilde{\K}_{ii} $ remains symmetric positive semidefinite. The explicit entries of $ \tilde{\A}_{0} $, $ \tilde{\A}_{i} $, and $ \tilde{\K} $ for this case can be found in, \eg, \cite{hughes1986new,fisher2013high}.

\section{Entropy-split discretization of the Euler equations}\label{sec:discrete analysis}
The spatial discretization of the entropy-split Euler equation on a periodic domain, \cref{eq:split Euler}, is handled using element-type SBP operators. Each element is coupled to its neighbors using SATs. On element $ \Omega_{k} $, we have the discretization
\begin{equation}\label{eq:discrete entropy split}
	\der[\bm{u}_{k}]t+\frac{\beta}{\beta+1}\sum_{i=1}^{d}\overline{\D}_{{x}_{i}k}\bm{f}_{k}^{I(i)}+\frac{1}{\beta+1}\sum_{i=1}^{d}\tilde{\A}_{i}\overline{\D}_{{x}_{i}k}\bm{w}_{k}=\text{SAT}_{k}^{I},
\end{equation}
where $ \text{SAT}_{k}^{I} $ denotes the inviscid SATs which are given by
\begin{equation}\label{eq:sats inv}
	\text{SAT}_{k}^{I}=\overline{\H}_{k}^{-1}\sum_{i=1}^{d}\overline{\E}_{{x}_{i}k}\bm{f}_{k}^{I(i)}-\frac{1}{2}\overline{\H}_{k}^{-1}\sum_{\gamma\in\Gamma_{k}^{I}}\sum_{i=1}^{d}\left[\overline{\E}_{{x}_{i}}^{kv}\circ\F_{i}\left(\bm{u}_{k},\bm{u}_{v}\right)\right]\bm{1},
\end{equation} 
where $ \circ $ denotes the Hadamard product, and $ \F_{i}\left(\bm{u}_{k},\bm{u}_{v}\right) \in \IRtwo{n_{pk} n_c}{n_{pv} n_c}$ is a block matrix containing evaluations of the two-point flux function of Sj{\"o}green and Yee \cite{sjogreen2019entropy} at the nodes of $ \Omega_{k} $ and $ \Omega_{v} $, \ie,
\begin{equation}\label{eq:F 2-point}
	\F_{i}\left(\bm{u}_{k},\bm{u}_{v}\right)=2\left[\begin{array}{ccc}
		\mydiag\left[\fnb F_{i}^{*}\left(\bm{u}_{k,1},\bm{u}_{v,1}\right)\right] & \dots & \mydiag\left[\fnb F_{i}^{*}\left(\bm{u}_{k,1},\bm{u}_{v,n_{pv}}\right)\right]\\
		\vdots & \ddots & \vdots\\
		\mydiag\left[\fnb F_{i}^{*}\left(\bm{u}_{k,n_{pk}},\bm{u}_{v,1}\right)\right] & \dots & \mydiag\left[\fnb F_{i}^{*}\left(\bm{u}_{k,n_{pk}},\bm{u}_{v,n_{pv}}\right)\right]
	\end{array}\right],
\end{equation}
where $ \fnb{F}^{*}_{i}$ is the $ i $-direction Sj{\"o}green -Yee two point flux, and $ n_{pk} $ and $ n_{pv} $ are the number of volume nodes in $ \Omega_{k} $ and $ \Omega_{v} $, respectively.  \blue{To define $ \fnb{F}^{*}_{i}$, we first introduce the parameter vector \cite{sjogreen2019entropy}
\begin{equation}
	\bm{z}=\left[\begin{array}{ccc}
		\frac{\rho}{p}\left(p\rho^{-\gamma}\right)^{\frac{1}{\alpha+\gamma}}, & \bm{V}^T, & p\end{array}\right]^{T},\label{eq:param vector}
\end{equation}
and the arithmetic and exponential averages defined by
\begin{align} 
	\avg{{u}}&=\frac{1}{2}\left({u}_{k}+{u}_{v}\right),\label{eq:avg arithmetic} \\
	\avg{{u}^{r-1}}_{e}&=\frac{1}{r}\frac{{u}_{k}^{r}-{u}_{v}^{r}}{{u}_{k}-{u}_{v}}, \quad \text{for } r\in\IR{}\label{eq:avg exponential}. 
\end{align}
The entries of the Sj{\"o}green-Yee two-point flux function in the $ i $-direction are given by \cite{sjogreen2019entropy}
\begin{equation}
	\mathtt{\begin{aligned}\left[\fnb F_{i}^{*}\right]_{1} & =\frac{\avg{\bm{z}_{1}\bm{z}_{(i+1)}}}{\avg{\bm{z}_{1}^{-\gamma/\alpha}}\avg{\bm{z}_{n_{c}}^{(1-\gamma-\alpha)/\alpha}}_{e}}\\
			\left[\fnb F_{i}^{*}\right]_{j} & =\left[\fnb F_{i}^{*}\right]_{1}\avg{z_{j}}+\delta_{(i+1)j}\avg{z_{n_{c}}},\quad j={2,\dots,n_{c}-1}\\
			\left[\fnb F_{i}^{*}\right]_{n_{c}} & =\left[\fnb F_{i}^{*}\right]_{1}\left(\frac{\gamma}{\gamma-1}\avg{\bm{z}_{n_{c}}^{-(\gamma-1)/\alpha}}\avg{\bm{z}_{1}^{-(\gamma+\alpha)/\alpha}}_{e}+\sum_{j=2}^{n_{c}-1}\avg{z_{j}}^{2}-\frac{1}{2}\avg{\sum_{j=2}^{n_{c}-1}z_{j}^{2}}\right).
		\end{aligned}
	}
\end{equation}} The two-point flux function satisfies three conditions: consistency, $ \fnb F_{i}^{*}\left(\fnb U,\fnb U\right)=\fnb F^{I(i)}\left(\fnb U\right) $, symmetry, $ \fnb F_{i}^{*}\left(\fnb U^{\prime},\fnb U^{\prime\prime}\right)=\fnb F_{i}^{*}\left(\fnb U^{\prime\prime},\fnb U^{\prime}\right) $, and Tadmor's entropy consistency condition \cite{tadmor1987numerical}, 
\begin{equation}\label{eq:Tadmor condition}
	\left(\fnb W^{\prime\prime}-\fnb W^{\prime}\right)^{T}\fnb F_{i}^{*}\left(\fnb U^{\prime},\fnb U^{\prime\prime}\right)=\left(\psi^{(i)}\right)^{\prime\prime}-\left(\psi^{(i)}\right)^{\prime},
\end{equation}
where $ (\cdot)^{\prime} $ and $ (\cdot)^{\prime\prime} $ indicate two states of a variable. For further discussion on the form of $ \F_{i}\left(\bm{u}_{k},\bm{u}_{v}\right) $ and its properties, see \cite{crean2018entropy}, and a general procedure to construct two-point fluxes can be found in \cite{ranocha2018comparison}.
\begin{remark}
	Although we have used two point fluxes in the inviscid SATs, they need only be evaluated where $ \overline{\E}_{x_{i}}^{kv} $ is nonzero. Since we are using diagonal-$ \E  $ SBP operators, the number of two-point flux evaluations per facet is equal to the number of single-point flux evaluations. Hence, the implication of using a two-point flux in the SATs on the computational cost is not significant.
\end{remark}

\subsection{High-order accuracy and consistency}
In this section, we show the high-order accuracy of the entropy-split discretization of the Euler equations \cref{eq:discrete entropy split}, by virtue of which we show the consistency of the scheme. We first present intermediate results through a series of lemmas. 
\begin{lemma}\label{lem:sat diag E}
	For SBP diagonal-$ \E $ operators, the following identity holds,
	\begin{equation}
		\overline{\E}_{{x}_{i}k}\bm{f}_{k}^{I(i)}-\frac{1}{2}\left[\overline{\E}_{{x}_{i}k}\circ\F_{i}\left(\bm{u}_{k},\bm{u}_{k}\right)\right]\bm{1}=\bm{0}.
	\end{equation}
\end{lemma}
\begin{proof}
	Since $ \overline{\E}_{{x}_{i}k} $ is diagonal, the $ m^{th} $ entry of the second vector on the LHS can be written as
	\begin{equation}
		\frac{1}{2}\left\{ \left[\overline{\E}_{{x}_{i}k}\circ\F_{i}\left(\bm{u}_{k},\bm{u}_{k}\right)\right]\bm{1}\right\} _{m}	=\frac{1}{2}\left[\overline{\E}_{{x}_{i}k}\right]_{mm}\F_{i}\left(\bm{u}_{k,m},\bm{u}_{k,m}\right)\\
		=\left[\overline{\E}_{{x}_{i}k}\right]_{mm}\fnb F_{i}^{*}\left(\bm{u}_{k,m},\bm{u}_{k,m}\right)\\
		=\left[\overline{\E}_{{x}_{i}k}\right]_{mm}\bm{f}^{(i)}_{k}\left(\bm{u}_{k,m}\right),
	\end{equation}
	where we have used the consistency of $ \fnb{F}^{*}_{i} $ in the last equality, and $ m=\{1, \dots, n_c n_{p}\} $. But we also have that
	\begin{equation}
		\left[\overline{\E}_{{x}_{i}k}\bm{f}_{k}^{I(i)}\right]_{m}=\left[\overline{\E}_{{x}_{i}k}\right]_{mm}\bm{f}^{(i)}_{k}\left(\bm{u}_{k,m}\right);
	\end{equation}
	thus, the result follows since the difference of the vectors at every entry is zero.
\end{proof}

\cref{lem:sat diag E} implies that the inviscid SATs, \cref{eq:sats inv}, can also be written in the form
\begin{equation}\label{eq:sats inv 2nd form}
	\text{SAT}_{k}^{I}=\frac{1}{2}\overline{\H}_{k}^{-1}\sum_{i=1}^{d}\left[\overline{\E}_{{x}_{i}k}\circ\F_{i}\left(\bm{u}_{k},\bm{u}_{k}\right)\right]\bm{1} - \frac{1}{2}\overline{\H}_{k}^{-1}\sum_{\gamma\in\Gamma_{k}^{I}}\sum_{i=1}^{d}\left[\overline{\E}_{{x}_{i}}^{kv}\circ\F_{i}\left(\bm{u}_{k},\bm{u}_{v}\right)\right]\bm{1},
\end{equation}
which has an identical structure as the inviscid SATs used in the Hadamard-form discretization of the Euler equations in \cite{crean2018entropy}. As will be shown later, this connection enables us to switch between the entropy-split and the Sj{\"o}green-Yee Hadamard-form discretizations.
\begin{lemma}\label{lem:sats vanish}
	For sufficiently smooth solutions, the SATs in \cref{eq:sats inv} or \cref{eq:sats inv 2nd form} vanish, \ie, $ \text{SAT}_{k}^{I}=\bm{0} $.
\end{lemma}
A proof of \cref{lem:sats vanish} for general types of SBP operators can be found in \cite[Theorem 2]{crean2018entropy}. For completeness, we have provided a proof applicable to the SBP diagonal-E operators. 
\begin{proof}
	Since SBP diagonal-$ \E $ operators are used, $ \R_{\gamma k} $ picks out the function values at the facet nodes. We thus have
	\begin{align*}
		\left\{ \frac{1}{2}\left[\overline{\E}_{{x}_{i}k}^{kk}\circ\F_{i}\left(\bm{u}_{k},\bm{u}_{k}\right)\right]\bm{1}\right\} _{m}&=\frac{1}{2}\sum_{j=1}^{n_{pk}}\left(\R_{\gamma k}^{T}\B_{\gamma}\N_{{x}_{i}\gamma}\R_{\gamma k}\right)_{mj}\fnb F_{i}^{*}\left(\bm{u}_{k,m},\bm{u}_{k,j}\right) \\&=\frac{1}{2}\sum_{n=1}^{n_f}\left(\R_{\gamma k}^{T}\B_{\gamma}\N_{{x}_{i}\gamma}\right)_{mn}\sum_{j=1}^{n_{pk}}\left(\R_{\gamma k}\right)_{nj}\fnb F_{i}^{*}\left(\bm{u}_{k,m},\bm{u}_{k,j}\right)\\&=\frac{1}{2}\sum_{n=1}^{n_f}\left(\R_{\gamma k}^{T}\B_{\gamma}\N_{{x}_{i}\gamma}\right)_{mn}\fnb F_{i}^{*}\left(\bm{u}_{k,m},\bm{u}_{\gamma k,n}\right),
	\end{align*}
	where $ m=\left\{ 1,\dots,n_{pk}\right\}  $. Using similar steps, we arrive at
	\begin{equation*}
		\left\{ \frac{1}{2}\left[\overline{\E}_{x_{i}}^{kv}\circ\F_{i}\left(\bm{u}_{k},\bm{u}_{v}\right)\right]\bm{1}\right\} _{m}=\frac{1}{2}\sum_{n=1}^{n_f}\left(\R_{\gamma k}^{T}\B_{\gamma}\N_{{x}_{i}\gamma}\right)_{mn}\fnb F_{i}^{*}\left(\bm{u}_{k,m},\bm{u}_{\gamma v,n}\right),
	\end{equation*}
	but due to the smoothness of the solution, we have $ \bm{u}_{\gamma k}=\bm{u}_{\gamma v} $; hence, the difference of the two vectors vanishes, \ie, 
	\begin{equation*}
		\frac{1}{2}\overline{\H}_{k}^{-1}\left[\overline{\E}_{{x}_{i}k}^{kk}\circ\F_{i}\left(\bm{u}_{k},\bm{u}_{k}\right)\right]\bm{1}- \frac{1}{2}\overline{\H}_{k}^{-1}\left[\overline{\E}_{x_{i}}^{kv}\circ\F_{i}\left(\bm{u}_{k},\bm{u}_{v}\right)\right]\bm{1}=\bm{0},
	\end{equation*}
	which, after summing over all facets of $ \Omega_{k} $, gives the desired result.
\end{proof}

\begin{lemma}\label{lem:Dx accuracy}
	Assuming $ \fnb{U} $ and $ \fnb{F}^{I(i)} $ are sufficiently smooth and that \cref{ass:mapping} holds, the entropy-split approximation of the derivatives of the Euler fluxes in \cref{eq:discrete entropy split} is design-order accurate, \ie, 
	\begin{equation}
		\frac{\beta}{\beta+1}\overline{\D}_{{x}_{i}k}\bm{f}_{k}^{I(i)}+\frac{1}{\beta+1}\tilde{\A}_{i}\overline{\D}_{{x}_{i}k}\bm{w}_{k}=\pder[\fnb F^{I(i)}]{{x}_{i}}+\fn O\left(h^{p}\right).
	\end{equation}
\end{lemma}
\begin{proof}
	The accuracy of the derivative operator for a vector containing a restriction of a sufficiently smooth function is shown in \cite[Theorem 9]{crean2018entropy}, \ie, $ \overline{\D}_{\bm{x}_i k}\bm{u}_k = \dx{\fnb{U}_k}{\bm{x}_i} +\fn{O}(h^p)$ for all $ \fnb{U}_k  \in \cont{p+1} $. This implies that 
	\begin{align*}
		\frac{\beta}{\beta+1}\overline{\D}_{{x}_{i}k}\bm{f}_{k}^{I(i)}+\frac{1}{\beta+1}\tilde{\A}_{i}\overline{\D}_{{x}_{i}k}\bm{w}_{k}&=\frac{\beta}{\beta+1}\pder[\fnb F^{I(i)}_k]{{x}_{i}}+\frac{1}{\beta+1}\pder[\fnb F^{I(i)}_k]{\fnb W}\pder[\fnb W]{{x}_{i}}+\fn O\left(h^{p}\right)
		=\pder[\fnb F^{I(i)}]{{x}_{i}}+\fn O\left(h^{p}\right),
	\end{align*}
	as desired.
\end{proof}

\begin{theorem}\label{thm:high-order}
	The entropy-split discretization of the Euler equations, \cref{eq:discrete entropy split}, is consistent and high-order accurate, \ie, for sufficiently smooth solutions, the discretization approximates the Euler equations to the design order of accuracy. 
\end{theorem}
\begin{proof}
	The proof follows directly from \cref{lem:sat diag E,lem:sats vanish,lem:Dx accuracy}.
\end{proof}

\subsection{Entropy conservation}
For the discrete entropy conservation analysis, we follow similar steps employed in the corresponding continuous entropy conservation analysis in \cref{sec:continuous analysis}. We premultiply the entropy-split discretization of the Euler equations, \cref{eq:discrete entropy split}, by $ \bm{w}_{k}^{T}\overline{\H}_{k}$ to find 
\begin{equation}\label{eq:discrete entropy conservation 0}
	\begin{aligned}
		\bm{w}_{k}^{T}\overline{\H}_{k}\der[\bm{u}_{k}]t+\frac{\beta}{\beta+1}\sum_{i=1}^{d}\bm{w}_{k}^{T}\overline{\H}_{k}\overline{\D}_{{x}_{i}k}\bm{f}_{k}^{I(i)}&+\frac{1}{\beta+1}\sum_{i=1}^{d}\bm{w}_{k}^{T}\overline{\H}_{k}\tilde{\A}_{i}\overline{\D}_{{x}_{i}k}\bm{w}_{k}
		\\& =\sum_{\gamma\in\Gamma_{k}^{I}}\sum_{i=1}^{d}\bm{w}_{k}^{T}\overline{\E}_{{x}_{i}k}^{kk}\bm{f}_{k}^{I(i)}-\frac{1}{2}\sum_{\gamma\in\Gamma_{k}^{I}}\sum_{i=1}^{d}\bm{w}_{k}^{T}\left[\overline{\E}_{{x}_{i}}^{kv}\circ\F_{i}\left(\bm{u}_{k},\bm{u}_{v}\right)\right]\bm{1}.
	\end{aligned}
\end{equation}
We first note that since $ \overline{\H}_{k} $ is diagonal, we have $ \bm{w}_{k}^{T}\overline{\H}_{k}=\bm{1}^{T}\mydiag(\bm{w}_{k})\overline{\H}_{k}=\bm{1}^{T}\overline{\H}_{k}\mydiag(\bm{w}_{k})  $ and thus
\begin{equation} \label{eq:entropy conv time}
	\bm{w}_{k}^{T}\overline{\H}_{k}\der[\bm{u}_{k}]t=\bm{1}^{T}\overline{\H}_{k}\der[\bm{s}_{k}]t.
\end{equation}
In a similar manner, we analyze each term in \cref{eq:discrete entropy conservation 0} in a series of lemmas to facilitate the proof of entropy conservation.
\begin{lemma}\label{lem:sat1 to potential}
	The discrete facet integrals of the dot product of the entropy variables and the Euler fluxes can be expressed in terms of the potential fluxes, \ie, 
	\begin{align}
		\bm{w}_{k}^{T}\overline{\E}_{{x}_{i}k}\bm{f}_{k}^{I(i)}&=\left(\beta+1\right)\bm{1}^{T}{\E}_{{x}_{i}k}\bm{\psi}_{k}^{(i)},\\\bm{w}_{v}^{T}\overline{\E}_{{x}_{i}v}\bm{f}_{v}^{I(i)}&=\left(\beta+1\right)\bm{1}^{T}{\E}_{{x}_{i}v}\bm{\psi}_{v}^{(i)}.
	\end{align}
\end{lemma}
\begin{proof}
	Using the fact that $ \overline{\E}_{{x}_{i}k} $ is diagonal, we have that 
	\begin{equation*}
		\bm{w}_{k}^{T}\overline{\E}_{{x}_{i}k}\bm{f}_{k}^{I(i)}=\bm{1}^{T}\mydiag\left(\bm{w}_{k}\right)\overline{\E}_{{x}_{i}k}\bm{f}_{k}^{I(i)}=\bm{1}^{T}\overline{\E}_{{x}_{i}k}\mydiag\left(\bm{w}_{k}\right)\bm{f}_{k}^{I(i)}.
	\end{equation*}
	Identity \cref{eq:WTF=G} then gives 
	\begin{equation*}
		\bm{w}_{k}^{T}\overline{\E}_{{x}_{i}k}\bm{f}_{k}^{I(i)}=\frac{\beta+1}{\beta}\bm{1}^{T}{\E}_{{x}_{i}k}\bm{g}_{k}^{(i)},
	\end{equation*}
	which in turn, after application of \cref{eq:G=beta psi}, yields
	\begin{equation*}
		\bm{w}_{k}^{T}\overline{\E}_{{x}_{i}k}\bm{f}_{k}^{I(i)}=\left(\beta+1\right)\bm{1}^{T}{\E}_{{x}_{i}k}\bm{\psi}_{k}^{(i)}.
	\end{equation*}
	Following similar steps, we find
	\begin{equation*}
		\bm{w}_{v}^{T}\overline{\E}_{{x}_{i}v}\bm{f}_{v}^{I(i)}=\left(\beta+1\right)\bm{1}^{T}{\E}_{{x}_{i}v}\bm{\psi}_{v}^{(i)},
	\end{equation*}
	which concludes the proof.	
\end{proof}
\begin{lemma}\label{lem:entropy conv volume}
	The discrete volume integral of the dot product of the entropy variables and entropy-split flux derivatives can be written as a discrete facet integral of the potential fluxes, \ie, 
	\begin{equation} \label{eq:entropy conv volume}
		\frac{\beta}{\beta+1}\bm{w}_{k}^{T}\overline{\H}_{k}\overline{\D}_{{x}_{i}k}\bm{f}_{k}^{I(i)}+\frac{1}{\beta+1}\bm{w}_{k}^{T}\overline{\H}_{k}\tilde{\A}_{i}\overline{\D}_{{x}_{i}k}\bm{w}_{k}=\beta\bm{1}^{T}{\E}_{{x}_{i}k}\bm{\psi}_{k}^{(i)}.
	\end{equation} 
\end{lemma}
\begin{proof}
	Since the norm matrix is assumed to be diagonal and  $ \tilde{\A}_{i} $ is block diagonal, we have that
	\begin{equation*}
		\overline{\H}_{k}\tilde{\A}_{i}=\left[\H_{k}\otimes\I_{n_{c}}\right]\tilde{\A}_{i}=\tilde{\A}_{i}\left[\H_{k}\otimes\I_{n_{c}}\right].
	\end{equation*}
	Moreover, the homogeneity property \cref{eq:homo in W} and the symmetry of $ \tilde{\A}_{i} $ yield
	\begin{equation*}
		\bm{w}_{k}^{T}\tilde{\A}_{i}=\bm{w}_{k}^{T}\tilde{\A}_{i}^{T}=\beta\left[\bm{f}_{k}^{I(i)}\right]^{T}.
	\end{equation*}
	Hence, the second term on the LHS of  \cref{eq:entropy conv volume} can be written as
	\begin{equation*}
		\frac{1}{\beta+1}\bm{w}_{k}^{T}\overline{\H}_{k}\tilde{\A}_{i}\overline{\D}_{{x}_{i}k}\bm{w}_{k}=\frac{\beta}{\beta+1}\left[\bm{f}_{k}^{I(i)}\right]^{T}\overline{\H}_{k}\overline{\D}_{{x}_{i}k}\bm{w}_{k}=\frac{\beta}{\beta+1}\bm{w}_{k}^{T}\overline{\D}_{{x}_{i}k}^{T}\overline{\H}_{k}\bm{f}_{k}^{I(i)}.
	\end{equation*}
	Applying the SBP property, the first term on the LHS of \cref{eq:entropy conv volume}  becomes
	\begin{align*}
		\frac{\beta}{\beta+1}\bm{w}_{k}^{T}\overline{\H}_{k}\overline{\D}_{{x}_{i}k}\bm{f}_{k}^{I(i)}&=\frac{\beta}{\beta+1}\bm{w}_{k}^{T}\overline{\Q}_{{x}_{i}k}\bm{f}_{k}^{I(i)}=\frac{\beta}{\beta+1}\bm{w}_{k}^{T}\left[-\overline{\Q}_{{x}_{i}k}^{T}+\overline{\E}_{{x}_{i}k}\right]\bm{f}_{k}^{I(i)}\\
		&=-\frac{\beta}{\beta+1}\bm{w}_{k}^{T}\overline{\D}_{{x}_{i}k}^{T}\overline{\H}_{k}\bm{f}_{k}^{I(i)}+\frac{\beta}{\beta+1}\bm{w}_{k}^{T}\overline{\E}_{{x}_{i}k}\bm{f}_{k}^{I(i)}.
	\end{align*}
	Therefore, we obtain 
	\begin{equation*}
		\frac{\beta}{\beta+1}\bm{w}_{k}^{T}\overline{\H}_{k}\overline{\D}_{{x}_{i}k}\bm{f}_{k}^{I(i)}+\frac{1}{\beta+1}\bm{w}_{k}^{T}\overline{\H}_{k}\tilde{\A}_{i}\overline{\D}_{{x}_{i}k}\bm{w}_{k}=\frac{\beta}{\beta+1}\bm{w}_{k}^{T}\overline{\E}_{{x}_{i}k}\bm{f}_{k}^{I(i)} =\beta\bm{1}^{T}{\E}_{{x}_{i}k}\bm{\psi}_{k}^{(i)},
	\end{equation*}
	where we have used \cref{lem:sat1 to potential} in last equality.
\end{proof}

The next lemma and its proof can be found in \cite[Lemma 3]{crean2018entropy}. We reproduce the proof here for completeness.
\begin{lemma}\label{lem:sat2 to potential}
	At an interface shared by $ \Omega_{k} $ and $ \Omega_{v} $, we have that 
	\begin{equation}\label{eq:sat2 to potential}
		\frac{1}{2}\bm{w}_{k}^{T}\left[\overline{\E}_{{x}_{i}}^{kv}\circ\F_{i}\left(\bm{u}_{k},\bm{u}_{v}\right)\right]\bm{1}+\frac{1}{2}\bm{w}_{v}^{T}\left[\overline{\E}_{{x}_{i}}^{vk}\circ\F_{i}\left(\bm{u}_{v},\bm{u}_{k}\right)\right]\bm{1}=\bm{1}^{T}\E_{{x}_{i}}^{kk}\bm{\psi}_{k}^{(i)}+\bm{1}^{T}\E_{{x}_{i}}^{vv}\bm{\psi}_{v}^{(i)}.
	\end{equation}
\end{lemma}
\begin{proof}
	The LHS of \cref{eq:sat2 to potential} can be written as 
	\begin{align*}
		\frac{1}{2}\bm{w}_{k}^{T}\left[\overline{\E}_{{x}_{i}}^{kv}\circ\F_{i}\left(\bm{u}_{k},\bm{u}_{v}\right)\right]\bm{1}+\frac{1}{2}\bm{w}_{v}^{T}\left[\overline{\E}_{{x}_{i}}^{vk}\circ\F_{i}\left(\bm{u}_{v},\bm{u}_{k}\right)\right]\bm{1} &=\frac{1}{2}\bm{w}_{k}^{T}\left[\overline{\E}_{{x}_{i}}^{kv}\circ\F_{i}\left(\bm{u}_{k},\bm{u}_{v}\right)\right]\bm{1}-\frac{1}{2}\bm{w}_{v}^{T}\left[\left(\overline{\E}_{{x}_{i}}^{kv}\right)^{T}\circ\F_{i}\left(\bm{u}_{v},\bm{u}_{k}\right)\right]\bm{1}
		\\&=\sum_{m=1}^{n_{pk}}\sum_{n=1}^{n_{pv}}\left[\bm{w}_{k,m}^{T}\left({\E}_{{x}_{i}}^{kv}\right)_{mn}\fnb F_{i}^{*}\left(\bm{u}_{k,m},\bm{u}_{v,n}\right)-\bm{w}_{v,n}^{T}\left({\E}_{{x}_{i}}^{kv}\right)_{mn}\fnb F_{i}^{*}\left(\bm{u}_{v,n},\bm{u}_{k,m}\right)\right]\\&=\sum_{m=1}^{n_{pk}}\sum_{n=1}^{n_{pv}}\left({\E}_{{x}_{i}}^{kv}\right)_{mn}\left(\bm{w}_{k,m}-\bm{w}_{v,n}\right)^{T}\fnb F_{i}^{*}\left(\bm{u}_{k,m},\bm{u}_{v,n}\right)\\&=\sum_{m=1}^{n_{pk}}\sum_{n=1}^{n_{pv}}\left({\E}_{{x}_{i}}^{kv}\right)_{mn}\left(\bm{\psi}_{k,m}^{(i)}-\bm{\psi}_{v,n}^{(i)}\right)\\&=\bm{1}^T\R_{\gamma v}^{T}\B_{\gamma}\N_{\gamma{x}_{i}}\R_{\gamma k}\bm{\psi}_{k}^{(i)}-\bm{1}^{T}\R_{\gamma k}^{T}\B_{\gamma}\N_{\gamma{x}_{i}}\R_{\gamma v}\bm{\psi}_{v}^{(i)}\\&=\bm{1}^{T}\R_{\gamma k}^{T}\B_{\gamma}\N_{\gamma{x}_{i}}\R_{\gamma k}\bm{\psi}_{k}^{(i)}-\bm{1}^{T}\R_{\gamma v}^{T}\B_{\gamma}\N_{\gamma{x}_{i}}\R_{\gamma v}\bm{\psi}_{v}^{(i)}\\&=\bm{1}^{T}{\E}_{\bm{x}_{i}}^{kk}\bm{\psi}_{k}^{(i)}+\bm{1}^{T}{\E}_{\bm{x}_{i}}^{vv}\bm{\psi}_{v}^{(i)},
	\end{align*}
	where we have used ${\E}_{{x}_{i}}^{kv}=-({\E}_{{x}_{i}}^{vk})^{T}  $ in the first equality, the symmetry of two-point fluxes in the third equality, Tadmor's consistency condition, \cref{eq:Tadmor condition}, in the fourth equality, and the fact that $ \R_{\gamma v}\bm{1}=\R_{\gamma k}\bm{1} $ in the penultimate equality. 
\end{proof}

We are now ready to present the main result of the paper.
\begin{theorem}\label{thm:entropy conservation}
	On a periodic domain, $ \Omega $, the entropy-split SBP-SAT discretization of the Euler equations, \cref{eq:discrete entropy split}, with smooth solutions is entropy conservative, \ie, 
	\begin{equation}
		\sum_{\Omega_{k}\in\fn T_{h}}\bm{1}^{T}\H_{k}\der[\bm{s}]t=0.
	\end{equation}
\end{theorem}
\begin{proof}
	We start from \cref{eq:discrete entropy conservation 0}, which is obtained by premultiplying \cref{eq:discrete entropy split} by $ \bm{w}_k^{T}\overline{\H}_k $. We then note that \cref{lem:entropy conv volume,lem:sat1 to potential} imply that 
	\begin{align*}
		-\frac{\beta}{\beta+1}\bm{w}_{k}^{T}\overline{\H}_{k}\overline{\D}_{{x}_{i}k}\bm{f}_{k}^{I(i)}-\frac{1}{\beta+1}\bm{w}_{k}^{T}\overline{\H}_{k}\tilde{\A}_{i}\overline{\D}_{{x}_{i}k}\bm{w}_{k}+\bm{w}_{k}^{T}\overline{\E}_{{x}_{i}k}\bm{f}_{k}^{I(i)}&=-\beta\bm{1}^{T}{\E}_{{x}_{i}k}\psi_{k}^{(i)}+\left(\beta+1\right)\bm{1}^{T}{\E}_{{x}_{i}k}\bm{\psi}_{k}^{(i)} \\
		&=\bm{1}^{T}{\E}_{{x}_{i}k}\psi_{k}^{(i)},
	\end{align*}
	which along with \cref{eq:entropy conv time}, allows us to write \cref{eq:discrete entropy conservation 0} as
	\begin{equation}\label{eq:entropy conv proof 0}
		\bm{1}^{T}\overline{\H}_{k}\der[\bm{s}_{k}]t=\sum_{i=1}^{d}\bm{1}^{T}{\E}_{{x}_{i}k}\psi_{k}^{(i)}-\frac{1}{2}\sum_{\gamma\in\Gamma_{k}^{I}}\sum_{i=1}^{d}\bm{w}_{k}^{T}\left[\overline{\E}_{{x}_{i}}^{kv}\circ\F_{i}\left(\bm{u}_{k},\bm{u}_{v}\right)\right]\bm{1}.
	\end{equation}
	Summing over all elements, and collecting terms by facet, we have 
	\begin{equation*}
		\sum_{\Omega_{k}\in\fn T_{h}}\bm{1}^{T}\overline{\H}_{k}\der[\bm{s}_{k}]t=\sum_{\gamma\in\Gamma^{I}}\sum_{i=1}^{d}\left(\bm{1}^{T}{\E}_{{x}_{i}k}^{kk}\bm{\psi}_{k}^{(i)}+\bm{1}^{T}{\E}_{{x}_{i}v}^{vv}\bm{\psi}_{v}^{(i)}\right)-\sum_{\gamma\in\Gamma^{I}}\sum_{i=1}^{d}\left(\frac{1}{2}\bm{w}_{k}^{T}\left[\overline{\E}_{{x}_{i}}^{kv}\circ\F_{i}\left(\bm{u}_{k},\bm{u}_{v}\right)\right]\bm{1}+\frac{1}{2}\bm{w}_{v}^{T}\left[\overline{\E}_{{x}_{i}}^{vk}\circ\F_{i}\left(\bm{u}_{v},\bm{u}_{k}\right)\right]\bm{1}\right).
	\end{equation*}
	Finally, applying \cref{lem:sat2 to potential}, we obtain
	\begin{align*}
		\sum_{\Omega_{k}\in\fn T_{h}}\bm{1}^{T}\overline{\H}_{k}\der[\bm{s}_{k}]t&=\sum_{\gamma\in\Gamma^{I}}\sum_{i=1}^{d}\left(\bm{1}^{T}{\E}_{{x}_{i}k}^{kk}\bm{\psi}_{k}^{(i)}+\bm{1}^{T}{\E}_{{x}_{i}v}^{vv}\bm{\psi}_{v}^{(i)}\right)-\sum_{\gamma\in\Gamma^{I}}\sum_{i=1}^{d}\left(\bm{1}^{T}\E_{{x}_{i}}^{kk}\bm{\psi}_{k}^{(i)}+\bm{1}^{T}\E_{{x}_{i}}^{vv}\bm{\psi}_{v}^{(i)}\right) =0,
	\end{align*}
	as desired.
\end{proof}

\subsection{Primary conservation error} \label{subsec:conservation}
The entropy-split scheme does not conserve mass, momentum, and energy due to the nonconservative portion of the split Euler flux derivatives. The goal of this subsection is to find an estimate of the primary conservation error for smooth solutions and understand how it behaves with mesh resolution and order of the spatial discretization.

To study the conservation of all components of the vector of conservative variables simultaneously, we first introduce the block identity vector, $  \overline{\I}\in\IR{n_{p}n_{c}\times n_{c}} $ defined as 
\begin{align}
	\overline{\I}&\coloneqq\bm{1}\otimes\I_{n_{c}},
\end{align}
where $ \bm{1} $ is of size $ n_{p}\times1 $ in this case. We will now present a lemma that will be used to derive the conservation error.
\begin{lemma}\cite[Theorem 3]{crean2018entropy}\label{lem:conv sat cancel}
	At a shared facet $ \gamma \in \{\Gamma_{k}\cap \Gamma_{v}\} $, the following identity holds,
	\begin{equation}\label{eq:convserv sat}
		\overline{\I}^{T}\left[\overline{\E}_{{x}_{i}}^{kv}\circ\F_{i}\left(\bm{u}_{k},\bm{u}_{v}\right)\right]\bm{1}+\overline{\I}^{T}\left[\overline{\E}_{{x}_{i}}^{vk}\circ\F_{i}\left(\bm{u}_{v},\bm{u}_{k}\right)\right]\bm{1}=\bm{0}.
	\end{equation}
\end{lemma}
\begin{proof}
	Using $ \E_{{x}_{i}}^{kv}=-\left(\E_{{x}_{i}}^{vk}\right)^{T} $ and considering the $  j $-th entry of the LHS of \cref{eq:convserv sat} at the interface $ \gamma $, we have
	\begin{align*}
		\left\{ \bar{\I}^{T}\left[\overline{\E}_{x_{i}}^{kv}\circ\F_{i}\left(\bm{u}_{k},\bm{u}_{v}\right)\right]\bm{1}+\bar{\I}^{T}\left[\overline{\E}_{x_{i}}^{vk}\circ\F_{i}\left(\bm{u}_{v},\bm{u}_{k}\right)\right]\bm{1}\right\} _{j}	&=\left\{ \bar{\I}^{T}\left[\overline{\E}_{x_{i}}^{kv}\circ\F_{i}\left(\bm{u}_{k},\bm{u}_{v}\right)\right]\bm{1}-\bar{\I}^{T}\left[\left(\overline{\E}_{{x}_{i}}^{kv}\right)^{T}\circ\F_{i}\left(\bm{u}_{v},\bm{u}_{k}\right)\right]\bm{1}\right\} _{j}\\
		&=\bm{1}^{T}\left[\E_{x_{i}}^{kv}\circ\F_{i}^{(j)}\left(\bm{u}_{k},\bm{u}_{v}\right)\right]\bm{1}-\bm{1}^{T}\left[\left(\E_{x_{i}}^{kv}\right)^{T}\circ\F_{i}^{(j)}\left(\bm{u}_{k},\bm{u}_{v}\right)^{T}\right]\bm{1}\\
		&=\bm{1}^{T}\left[\E_{x_{i}}^{kv}\circ\F_{i}^{(j)}\left(\bm{u}_{k},\bm{u}_{v}\right)\right]\bm{1}-\bm{1}^{T}\left[\E_{x_{i}}^{kv}\circ\F_{i}^{(j)}\left(\bm{u}_{k},\bm{u}_{v}\right)\right]\bm{1}\\
		&=0,
	\end{align*}
	where $  j=\left\{ 1,\dots,n_{c}\right\}  $ and $  \F_{i}^{(j)} $ is a matrix containing the two-point flux evaluations in $  \F_{i} $ corresponding to the $ j $-th component of the conservative variables. The symmetry of $  \F_{i}^{(j)} $ is used in the second equality, \ie, $  \F_{i}^{(j)} (\bm{u}_v, \bm{u}_k) = \F_{i}^{(j)} (\bm{u}_k, \bm{u}_v)^T$, and the property of Hadamard products, $( \A \circ \B)^T = \A^T \circ \B^T $, is used in the penultimate equality.
\end{proof}

\begin{theorem}\label{thm:conservation}
	For element-wise sufficiently smooth solution and fluxes on a periodic domain, the entropy-split scheme, \cref{eq:discrete entropy split}, is conservative up to the design order of the discretization, \ie, 
	\begin{equation} \label{eq:conservation}
		\sum_{\Omega_{k}\in\fn T_{h}}\overline{\I}^{T}\overline{\H}_{k}\der[\bm{u}_{k}]t= \overline{\fn O}\left(h^{p+1}\right).
	\end{equation}
\end{theorem}
\begin{proof}
	We start by premultiplying \cref{eq:discrete entropy split} by $ \overline{\I}^{T}\overline{\H}_{k} $ to find
	\begin{equation}\label{eq:conv proof 0}
		\overline{\I}^{T}\overline{\H}_{k}\der[\bm{u}_{k}]t+\frac{\beta}{\beta+1}\sum_{i=1}^{d}\overline{\I}^{T}\overline{\H}_{k}\overline{\D}_{{x}_{i}k}\bm{f}_{k}^{I(i)}+\frac{1}{\beta+1}\sum_{i=1}^{d}\overline{\I}^{T}\overline{\H}_{k}\tilde{\A}_{i}\overline{\D}_{{x}_{i}k}\bm{w}_{k}=\sum_{i=1}^{d}\overline{\I}^{T}\overline{\E}_{{x}_{i}k}\bm{f}_{k}^{I(i)}-\frac{1}{2}\sum_{\gamma\in\Gamma_{k}^{I}}\sum_{i=1}^{d}\overline{\I}^{T}\left[\overline{\E}_{{x}_{i}}^{kv}\circ\F_{i}\left(\bm{u}_{k},\bm{u}_{v}\right)\right]\bm{1}.
	\end{equation}
	Using the accuracy of the SBP operators, we can write
	\begin{align*}
		\frac{\beta}{\beta+1}\sum_{i=1}^{d}\overline{\I}^{T}\overline{\H}_{k}\overline{\D}_{{x}_{i}k}\bm{f}_{k}^{I(i)}+\frac{1}{\beta+1}\sum_{i=1}^{d}\overline{\I}^{T}\overline{\H}_{k}\tilde{\A}_{i}\overline{\D}_{{x}_{i}k}\bm{w}_{k}
		&=\frac{\beta}{\beta+1}\sum_{i=1}^{d}\int_{\Omega_{k}}\pder[\fnb F_{k}^{I(i)}]{{x}_{i}}\dd{\Omega}+\frac{1}{\beta+1}\sum_{i=1}^{d}\int_{\Omega_{k}}\tilde{\A}_{i}\pder[\fnb W]{{x}_{i}}\dd\Omega+\fn O\left(h^{p+1}\right)\\
		&=\sum_{i=1}^{d}\int_{\Omega_{k}}\pder[\fnb F_{k}^{I(i)}]{{x}_{i}}\dd{\Omega}+\fn O\left(h^{p+1}\right)
		\\&=\sum_{i=1}^{d}\int_{\Gamma_{k}}\fnb F_{k}^{I(i)}\bm{n}_{x_{i}}\dd{\Gamma}+\fn O\left(h^{p+1}\right)\\
		&=\sum_{i=1}^{d}\overline{\I}^{T}\overline{\E}_{{x}_{i}k}\bm{f}_{k}^{I(i)}+\fn O\left(h^{p+1}\right).
	\end{align*}
	Hence, \cref{eq:conv proof 0} simplifies to
	\begin{align*}
		\overline{\I}^{T}\overline{\H}_{k}\der[\bm{u}_{k}]t=-\frac{1}{2}\sum_{\gamma\in\Gamma_{k}^{I}}\sum_{i=1}^{d}\overline{\I}^{T}\left[\overline{\E}_{{x}_{i}}^{kv}\circ\F_{i}\left(\bm{u}_{k},\bm{u}_{v}\right)\right]\bm{1}+\fn O\left(h^{p+1}\right).
	\end{align*}
	Summing over all elements and collecting terms by facet, we find 
	\begin{align*}
		\sum_{\Omega_{k}\in\fn T_{h}}\bm{1}^{T}\H_{k}\der[\bm{u}_{k}]t&= -\frac{1}{2}\sum_{\gamma\in\Gamma^{I}}\sum_{i=1}^{d}\left(\overline{\I}^{T}\left[\overline{\E}_{{x}_{i}}^{kv}\circ\F_{i}\left(\bm{u}_{k},\bm{u}_{v}\right)\right]\bm{1}+\overline{\I}^{T}\left[\overline{\E}_{{x}_{i}}^{vk}\circ\F_{i}\left(\bm{u}_{v},\bm{u}_{k}\right)\right]\bm{1}\right)+\fn O\left(h^{p+1}\right),
	\end{align*}
	which, after applying  \cref{lem:conv sat cancel}, gives the desired result.
\end{proof}

\begin{remark}\label{rem:conservation}
	The term on the LHS of \cref{eq:conservation} is a volume functional. It is well known that, for smooth enough problems, if the discretization is adjoint (dual) consistent, then certain types of functionals superconverge, \eg, see \cite{hicken2011superconvergent,hicken2014dual,hicken2012output,hartmann2013higher,penner2020superconvergent,penner2022accurate,worku2021simultaneous}. 
\end{remark}
The implication of \cref{thm:conservation} is that if the solution and fluxes are sufficiently smooth on each element, then the primary conservation error vanishes with mesh refinement or with increasing order of the discretization. However, for problems with solution discontinuities, a different strategy must be sought, as the theorem no longer applies unless discontinuities are aligned with element interfaces. The next section discusses a strategy to handle precisely this issue.

\subsection{Entropy conservative hybrid scheme}
According to the Lax-Wendroff theorem \cite{lax1960systems}, conservation is sufficient to ensure convergence to a weak solution, provided that the numerical scheme for hyperbolic systems of equations is convergent. The entropy-split scheme is not conservative; hence, convergence to a weak solution is not guaranteed. However, Hou and LeFloch \cite{hou1994nonconservative} showed that for scalar conservation laws, local conservation is sufficient to ensure convergence to the right weak solutions. Accordingly, they proposed a hybrid scheme that applies a nonconservative method for smooth portions and a conservative method near discontinuities as a remedy. A similar approach has been applied to systems of conservation laws in \cite{sjogreen2021construction} where the entropy-split scheme is replaced by a conservative method near shocks, which enables convergence to the correct weak solution but at the cost of losing entropy stability. 

In this work, we enforce local conservation in critical portions of the computational domain by switching between the entropy-split and the Sj{\"o}green-Yee Hadamard-form discretizations while maintaining entropy conservation. This strategy allows to exploit both the conservation property of the Hadamard-form scheme and efficiency of the entropy-split scheme. The next theorem stipulates that the resulting hybrid scheme is entropy conservative. 

\begin{theorem}\label{thm:switching hybrid}
	The Hadamard-form entropy conservative discretization based on the Sj{\"o}green-Yee two-point fluxes \cite{sjogreen2019entropy} can be applied on a single element or a group of elements without destroying the entropy conservation property of the entropy-split scheme, \cref{eq:discrete entropy split}. Local conservation is achieved on the elements where the Hadamard-form discretization is applied. Furthermore, the discretization on a single or group of elements may be switched between the Hadamard-form and entropy-split discretizations at different timesteps without loss of entropy conservation of the resulting hybrid scheme.
\end{theorem}

\begin{proof}
	The result is a consequence of the shared structure of the inviscid SATs of the Hadamard-form and entropy-split schemes as well as the shared family of entropy functions used to construct the Sj{\"o}green-Yee two-point fluxes and the entropy-split schemes. The Hadamard-form discretization of the Euler equations can be written as \cite{crean2018entropy,shadpey2020entropy}
	\begin{equation}\label{eq:2pt disc}
		\der[\bm{u}_{k}]t+\sum_{i=1}^{d}\left[\overline{\D}_{{x}_{i}k}\circ\F_{i}\left(\bm{u}_{k},\bm{u}_{k}\right)\right]\bm{1}=\frac{1}{2}\overline{\H}_{k}^{-1}\sum_{i=1}^{d}\left[\overline{\E}_{{x}_{i}k}\circ\F_{i}\left(\bm{u}_{k},\bm{u}_{k}\right)\right]\bm{1}-\frac{1}{2}\overline{\H}_{k}^{-1}\sum_{\gamma\in\Gamma_{k}^{I}}\sum_{i=1}^{d}\left[\overline{\E}_{{x}_{i}}^{kv}\circ\F_{i}\left(\bm{u}_{k},\bm{u}_{v}\right)\right]\bm{1}.
	\end{equation}
	Note that the first term on the RHS of \cref{eq:2pt disc} is equal to the first SAT term in \cref{eq:discrete entropy split} due to \cref{lem:sat diag E}. Premultiplying \cref{eq:2pt disc} by $ \bm{w}_{k}^{T}\overline{\H}_{k} $ and using the identity $\overline{\S}_{{x}_{i}}=\overline{\Q}_{{x}_{i}}-\frac{1}{2}\overline{\E}_{{x}_{i}} $, we find 
	\begin{equation} \label{eq:hybrid proof 0}
		\bm{1}^{T}\overline{\H}_{k}\der[\bm{s}_{k}]t+\sum_{i=1}^{d}\bm{w}_{k}^{T}\left[\overline{\S}_{{x}_{i}k}\circ\F_{i}\left(\bm{u}_{k},\bm{u}_{k}\right)\right]\bm{1}=-\frac{1}{2}\sum_{\gamma\in\Gamma_{k}^{I}}\sum_{i=1}^{d}\bm{w}_{k}^{T}\left[\overline{\E}_{{x}_{i}}^{kv}\circ\F_{i}\left(\bm{u}_{k},\bm{u}_{v}\right)\right]\bm{1}.
	\end{equation}
	Using properties of the two-point flux function, it is possible to show that (e.g., see \cite{crean2018entropy}) 
	\begin{equation*}
		\bm{w}_{k}^{T}\left[\overline{\S}_{\bm{x}_i k}\circ\F_{i}\left(\bm{u}_{k},\bm{u}_{k}\right)\right]\bm{1}=-\bm{1}^{T}\E_{{x}_{i}k}\bm{\psi}_{k}^{(i)}.
	\end{equation*}
	Substituting the above in \cref{eq:hybrid proof 0} and rearranging gives
	\begin{equation*}
		\bm{1}^{T}\overline{\H}_{k}\der[\bm{s}_{k}]t=\sum_{i=1}^{d}\bm{1}^{T}\E_{{x}_{i}k}\bm{\psi}_{k}^{(i)}-\frac{1}{2}\sum_{\gamma\in\Gamma_{k}^{I}}\sum_{i=1}^{d}\bm{w}_{k}^{T}\left[\overline{\E}_{{x}_{i}}^{kv}\circ\F_{i}\left(\bm{u}_{k},\bm{u}_{v}\right)\right]\bm{1}.
	\end{equation*}
	Comparing the last equality with \cref{eq:entropy conv proof 0}, we conclude that the term 
	\begin{equation*}
		\left[\overline{\D}_{{x}_{i}k}\circ\F_{i}\left(\bm{u}_{k},\bm{u}_{k}\right)\right]\bm{1}
	\end{equation*}
	in \cref{eq:2pt disc} and the split flux derivative,
	\begin{equation*}
		\frac{\beta}{\beta+1}\overline{\D}_{{x}_{i}k}\bm{f}_{k}^{I(i)}+\frac{1}{\beta+1}\tilde{\A}_{i}\overline{\D}_{{x}_{i}k}\bm{w}_{k},
	\end{equation*}
	in \cref{eq:discrete entropy split} can be used interchangeably on any element in the domain at any timestep without compromising the entropy conservation property of the scheme. The local conservation of the Hadamard-form discretization follows from the elementwise conservation result in \cite[Theorem 3]{crean2018entropy}.
\end{proof}

\blue{Although not pursued in this work, in practice, the hybrid discretization should be coupled with a shock wave detection mechanism as in \cite{sjogreen2021construction,yee2023recent}, where a shock detector is used to switch between a primary conservative scheme and the entropy-split method. Automatic switching between the Sj{\"o}green-Yee Hadamard-form and entropy-split discretizations using a shock detector is essential to efficiently handle discontinuous solutions and enhance flexibility of the hybrid scheme.}

\subsection{\blue{Interface} dissipation}\label{subsec:interface dissipation}
Thermodynamic entropy increases for irreversible processes; thus, the mathematical entropy should be dissipated for flows with shock waves. \blue{Interface} dissipation is added to the entropy-split, \cref{eq:discrete entropy split}, Hadamard-form, \cref{eq:2pt disc}, and hybrid discretizations to produce entropy stable variants of the schemes. Following \cite{ismail2009affordable}, we consider the matrix-form interface dissipation term,
\begin{equation}\label{eq:dissipation operator}
	-\frac{1}{2}\overline{\H}_{k}^{-1}\overline{\R}_{\gamma k}^{T}\overline{\B}_{\gamma} \hat{\X}\left|\Lambda_{n}\right|\hat{\S}\hat{\X}^{T}\left(\bm{w}_{\gamma k}-\bm{w}_{\gamma v}\right),
\end{equation}
which is added to the inviscid SAT terms in  \cref{eq:discrete entropy split} and \cref{eq:2pt disc}. The hat $ \left(\hat{\cdot}\right)  $ indicates that a quantity is computed at a certain intermediate state between the states $ \bm{w}_{\gamma k}=\overline{\R}_{\gamma k}\bm{w}_{k} $ and $ \bm{w}_{\gamma v}=\overline{\R}_{\gamma v}\bm{w}_{v} $ at a shared interface $ \gamma $ of $ \Omega_{k} $ and $ \Omega_{v} $. The dissipation term, \cref{eq:dissipation operator}, is modeled after a dissipation operator involving the jump in conservative variables, 
\begin{equation}\label{eq:dissipation in conservative}
	-\frac{1}{2}\overline{\H}_{k}^{-1}\overline{\R}_{\gamma k}^{T}\overline{\B}_{\gamma}\hat{\X}\left|\Lambda_{n}\right|\hat{\X}^{-1}\left(\bm{u}_{\gamma k}-\bm{u}_{\gamma v}\right).
\end{equation}
Computing $[ \tilde{\A}_{0}]_{jj} =[\dx{\fnb{U}}{\fnb{W}}]_j $, for $  j\in \{1,\dots,n_f\} $, at the same intermediate state as $ \hat{\X} $, we have 
\begin{equation}
	-\frac{1}{2}\overline{\H}_{k}^{-1}\overline{\R}_{\gamma k}^{T}\overline{\B}_{\gamma}\hat{\X}\left|\Lambda_{n}\right|\hat{\X}^{-1}\left(\bm{u}_{\gamma k}-\bm{u}_{\gamma v}\right)_j \approx-\frac{1}{2}\overline{\H}_{k}^{-1}\overline{\R}_{\gamma k}^{T}\overline{\B}_{\gamma}\hat{\X}\left|\Lambda_{n}\right|\hat{\X}^{-1}\hat{\tilde{\A}}_{0}\left(\bm{w}_{\gamma k}-\bm{w}_{\gamma v}\right).
\end{equation}
Applying the eigenvector scaling theorem of Barth \cite{barth1999numerical}, there exists a diagonal positive definite scaling matrix, $ \hat{\S} $, such that $ \hat{\tilde{\A}}_{0}=\hat{\X}\hat{\S}\hat{\X}^{T} $. It follows that 
\begin{equation}
	-\frac{1}{2}\overline{\H}_{k}^{-1}\overline{\R}_{\gamma k}^{T}\overline{\B}_{\gamma}\hat{\X}\left|\Lambda_{n}\right|\hat{\X}^{-1}\left(\bm{u}_{\gamma k}-\bm{u}_{\gamma v}\right)_j \approx-\frac{1}{2}\overline{\H}_{k}^{-1}\overline{\R}_{\gamma k}^{T}\overline{\B}_{\gamma}\hat{\X}\left|\Lambda_{n}\right|\hat{\S}\hat{\X}^{T}\left(\bm{w}_{\gamma k}-\bm{w}_{\gamma v}\right).
\end{equation}

If an appropriate intermediate state is used, \cref{eq:dissipation operator} produces an entropy stable scheme when used with an entropy conservative discretization. We note that since $ \overline{\B}_{\gamma} = \B_{\gamma} \otimes \I_{n_c}$ is diagonal and $ \hat{\X} \in \IRtwo{n_cn_f}{n_c n_f}$ is a block diagonal matrix with blocks of size $ n_c\times n_c $, we have the commutativity property $  \overline{\B}_{\gamma}\hat{\X} = \hat{\X} \overline{\B}_{\gamma}$, which is crucial to show that \cref{eq:dissipation operator} is entropy dissipative. In contrast, the dissipation based on the jump of the conservative variables, \cref{eq:dissipation in conservative}, does not guarantee entropy stability. The intermediate states in this work are chosen such that the \blue{interface} dissipation term, \cref{eq:dissipation operator}, is consistent in the sense that in the limit of infinite mesh resolution we obtain $ \dd{\fnb U}=\X\S\X^{T}\dd{\fnb W} $ for smooth solutions. We used the symbolic manipulation software Maxima \cite{maxima} to find intermediate states that satisfy this consistency condition. 

The intermediate states of the primitive variables are computed \blue{in terms of the parameter vector, \cref{eq:param vector},} as 
\begin{align}
	\hat{\rho}&=\frac{\avg{\bm{z}_{1}}}{\avg{\bm{z}_{1}^{-\gamma/\alpha}}\avg{\bm{z}_{n_c}^{(1-\gamma-\alpha)/\alpha}}_{e}}, \\
	\hat{V}_{i}&=\frac{\avg{\bm{z}_{1}\bm{z}_{(i+1)}}}{\avg{\bm{z}_{1}}},\quad i=\left\{ 1,\dots,d\right\}, \\
	\hat{p}&=\avg{\bm{z}_{n_c}},
\end{align}
\blue{where the arithmetic and exponential averages are defined in \cref{eq:avg arithmetic} and \cref{eq:avg exponential}, respectively.} The scaling matrix $  \hat{\S} \in \IRtwo{n_cn_f}{n_cn_f}$ is a block diagonal matrix with diagonal blocks of size $ n_c\times n_c $ given by 
\begin{equation}
	[\hat{\S}]_{jj}=\mydiag\left(\frac{\hat{\theta}_j}{2\gamma},\frac{\left(\gamma-1\right)\left(\gamma+\alpha\right)\hat{\theta}_j}{\gamma\alpha},\hat{\eta}_j,\dots,\hat{\eta}_j,\frac{\hat{\theta}_j}{2\gamma}\right),\quad j\in\{1,\dots, n_f\},
\end{equation}
where if $ d\ge 2 $, then $[ \hat{\S}]_{jj} $ contains $ \hat{\eta}_j $ from the third up to the $ (d+1)$-st diagonal entry, and  
\begin{align}
	\hat{\theta}_j&=\left[\frac{\hat{\rho}}{\hat{p}\avg{\bm{z}_{1}^{-\gamma/\alpha}}\avg{\bm{z}_{n_c}^{(1-\gamma-\alpha)/\alpha}}_{e}}\right]_{j},\\
	\hat{\eta}_j&=\left[\frac{1}{\avg{\bm{z}_{1}}\avg{\bm{z}_{1}^{-(\gamma+\alpha)/\alpha}}_{e}\avg{\bm{z}_{n_c}^{(1-\gamma-\alpha)/\alpha}}_{e}}\right]_{j}.
\end{align}
Dropping the facet node indices from the fluid dynamics variables for convenience, the diagonal blocks of the absolute values of the eigenvalues in \cref{eq:dissipation operator} are given by 
\begin{equation}
	[\left|\Lambda_{n}\right|]_{jj}=\mydiag\left(\left|\hat{\bm{V}}\cdot\bm{n}-\hat{a}\right|,\left|\hat{\bm{V}}\cdot\bm{n}\right|,\left|\hat{\bm{V}}\cdot\bm{n}\right|,\left|\hat{\bm{V}}\cdot\bm{n}\right|,\left|\hat{\bm{V}}\cdot\bm{n}+\hat{a}\right|\right), \quad j\in\{1,\dots,n_f\},
\end{equation}
where $ \hat{\bm{V}}=\left[\hat{V}_1,\dots,\hat{V}_d\right]^{T} $, $ \bm{n}=\left[n_{\bm{x}_1},\dots,n_{\bm{x}_d}\right]^{T} $, and $ \hat{a}=\sqrt{\gamma\hat{p}/\hat{\rho}} $ .  Finally, $ \hat{\X} $ contains the right eigenvectors of the normal flux Jacobian computed at the intermediate states, and for the $ d=3 $ case, its diagonal blocks are given by \cite{rohde2001eigenvalues}
\begin{equation}
	[\hat{\X}]_{jj}=\left[\begin{array}{ccccc}
		1 & 1 & 0 & 0 & 1\\
		\hat{\red{V}}_{1}-\hat{a}n_{\red{x_{1}}} & \hat{\red{V}}_{1} & n_{\red{x_{2}}} & -n_{\red{x_{3}}} & \hat{\red{V}}_{1}+\hat{a}n_{\red{x_{1}}}\\
		\hat{\red{V}}_{2}-\hat{a}n_{\red{x_{2}}} & \hat{\red{V}}_{2} & -n_{\red{x_{1}}} & 0 & \hat{\red{V}}_{2}+\hat{a}n_{\red{x_{2}}}\\
		\hat{\red{V}}_{3}-\hat{a}n_{\red{x_{3}}} & \hat{\red{V}}_{3} & 0 & n_{\red{x_{1}}} & \hat{\red{V}}_{3}+\hat{a}n_{\red{x_{3}}}\\
		\hat{H}-\hat{a}\hat{\bm{V}}\cdot\bm{n} & \frac{1}{2}\hat{\bm{V}}\cdot\hat{\bm{V}} & \hat{\red{V}}_{1}n_{\red{x_{2}}}-\hat{\red{V}}_{2}n_{\red{x_{1}}} & \hat{\red{V}}_{3}n_{\red{x_{1}}}-\hat{\red{V}}_{1}n_{\red{x_{3}}} & \hat{H}+\hat{a}\hat{\bm{V}}\cdot\bm{n}
	\end{array}\right], \quad j\in\{1,\dots,n_f\},
\end{equation}
where $ \hat{H}=\frac{\hat{a}^{2}}{\gamma-1}+\frac{1}{2}\hat{\bm{V}}\cdot\hat{\bm{V}} $. For $ d<3 $, $[ \hat{\X}]_{jj} $ is obtained by eliminating the rows and columns corresponding to the extra dimension, \eg, if $ d=2 $, then the fourth row and column are eliminated.

The dissipation operator \cref{eq:dissipation operator} can be applied with the entropy-split, Sj{\"o}green-Yee Hadamard-form, or the hybrid entropy conservative schemes to obtain an entropy dissipative discretization. We note that in \cite{derigs2017novel, winters2017uniquely}, it was shown that it is important to satisfy a discrete version of the dissipation term consistency condition, \ie, 
\begin{equation}
	\left(\bm{u}_{\gamma k}-\bm{u}_{\gamma v}\right)=\hat{\X}\hat{\S}\hat{\X}^{T}\left(\bm{w}_{\gamma k}-\bm{w}_{\gamma v}\right),
\end{equation}
to avoid numerical problems that may arise in certain flow configurations \red{that involve very large jumps or Mach numbers}. The intermediate states presented in this work are not designed to satisfy this discrete consistency condition; \red{hence, they may cause numerical issues in some special cases}. A discretely consistent dissipation operator may improve robustness of the entropy dissipative schemes, but this line of research is out of the scope of this paper. 

\section{Extension to the \red{compressible} Navier-Stokes equations}\label{sec:nse}
Harten's family of entropy functions do not symmetrize the diffusivity tensor in the presence of heat fluxes, as discussed in \cref{sec:continuous analysis}. Nonetheless, we will extend the entropy-split scheme to the \red{compressible} Navier-Stokes equations and implement it for various problems to understand some of its properties. We assume that $ \tilde{\K} $ is symmetric positive semidefinite, which is true for Hadamard-form discretizations that are based on the entropy function \cref{eq:generalized entropy} and the entropy-split scheme without heat fluxes. When Harten's family of entropy functions, \cref{eq:harten entropy}, is used and the heat fluxes are neglected, the $ \text{LDL}^{T}  $ decomposition of $ \tilde{\K} $ in the three dimensional case gives 
\begin{equation}
	\text{D}=\mydiag\left(\left[\begin{array}{ccccccccccccccc}
		0 & \frac{4}{3}\delta & \delta & \delta & 0 & 0 & 0 & \delta & \delta & 0 & 0 & 0 & 0 & 0 & 0\end{array}\right]\right),
\end{equation}
where $ \delta=p^{\frac{3\gamma+\alpha-2}{(\gamma-1)(\alpha+\gamma)}}\rho^{-\frac{\gamma\left(3\gamma+\alpha-2\right)}{(\gamma-1)(\alpha+\gamma)}} $, affirming that $ \tilde{\K}\succeq 0 $ provided that pressure and density remain positive.

Assuming $ \tilde{\K}\succeq 0 $, an entropy stable entropy-split or Hadamard-form discretization of the \red{compressible} Navier-Stokes equations can be designed by adding the terms
\begin{equation}\label{eq:added term}
	\sum_{i=1}^{d}\sum_{j=1}^{d}\overline{\D}_{{x}_{i}k}\tilde{\K}_{ij,k}\overline{\D}_{{x}_{j}k}\bm{w}_{k}+\text{SAT}_{k}^{V}
\end{equation}
to the RHS of \cref{eq:discrete entropy split} or \cref{eq:2pt disc}, respectively, where $ \text{SAT}_{k}^{V} $ denotes the viscous SATs \cite{yan2018interior,worku2021simultaneous},
\begin{equation} \label{eq:sats visc}
	\text{SAT}_{k}^{V}= - \overline{\H}_{k}^{-1}\sum_{\gamma\in\Gamma_{k}^{I}}\left[\begin{array}{cc}
		\overline{\R}_{\gamma k}^{T} & \overline{\D}_{\gamma k}^{T}\end{array}\right]\left[\begin{array}{cc}
		\overline{\T}_{\gamma k}^{(1)} & \overline{\T}_{\gamma k}^{(3)}\\
		\overline{\T}_{\gamma k}^{(2)} & [{\bm{0}}]
	\end{array}\right]\left[\begin{array}{c}
		\overline{\R}_{\gamma k}\bm{w}_{k}-\overline{\R}_{\gamma v}\bm{w}_{v}\\
		\overline{\D}_{\gamma k}\bm{w}_{k}+\overline{\D}_{\gamma v}\bm{w}_{v}
	\end{array}\right],
\end{equation}
$ \overline{\D}_{\gamma k} $ and $ \overline{\D}_{\gamma v} $ are facet derivative matrices given by 
\begin{align}\label{eq:Dgk}
	\overline{\D}_{\gamma k}=\sum_{i=1}^{d}\sum_{j=1}^{d}\overline{\N}_{{x}_{i}\gamma}\overline{\R}_{\gamma k}\tilde{\K}_{ij,k}\overline{\D}_{{x}_{j}k}, 
	\qquad \overline{\D}_{\gamma v}= - \sum_{i=1}^{d}\sum_{j=1}^{d}\overline{\N}_{{x}_{i}\gamma}\overline{\R}_{\gamma v}\tilde{\K}_{ij,v}\overline{\D}_{{x}_{j}v},
\end{align}
and $ \overline{\T}_{\gamma k}^{(j)} \in \IR{n_c n_{pk} \times n_c n_{pk} }$, $ j=\{1,\dots,3\} $, are symmetric penalty coefficient matrices. The viscous SATs, \cref{eq:sats visc}, couple immediate neighboring elements only, and generalize many known discontinuous Galerkin fluxes for elliptic problems (see, \eg, \cite{arnold2002unified}) such as the second method of Bassi and Rebay (BR2), \cite{bassi1997highbr2}, the symmetric interior penalty Galerkin method (SIPG) \cite{douglas1976interior}, the non-symmetric interior penalty Galerkin method (NIPG) \cite{riviere1999improved}, the Baumann-Oden method (BO) \cite{baumann1999discontinuous}, the Carpenter-Nordstr{\"o}m-Gottlieb method (CNG) \cite{carpenter1999stable}, and the compact discontinuous Galerkin (CDG) method \cite{peraire2008compact}. The viscous SATs can also be extended to include methods that couple second neighbor elements such as the first method of Bassi and Rebay (BR1) \cite{bassi1997highnavier} and the local discontinuous Galerkin (LDG) method \cite{cockburn1998local}. In this work, we do not consider SATs that couple second neighbor elements, as the BR1 and LDG methods are equivalent to the BR2 and CDG methods, respectively, when used with diagonal-$ \E $ multidimensional SBP operators \cite{worku2021simultaneous}.

Our task now is to find penalty coefficients such that 
\begin{equation} \label{eq:visc entropy}
	\sum_{i=1}^{d}\sum_{j=1}^{d}\bm{w}_k^T \overline{\H}_k \overline{\D}_{{x}_{i}k}\tilde{\K}_{ij,k}\overline{\D}_{{x}_{j}k}\bm{w}_{k} + \bm{w}_k^T \overline{\H}_k\text{SAT}_{k}^{V} \le 0,
\end{equation}
so that \cref{eq:added term} dissipates entropy when added to \cref{eq:discrete entropy split} or \cref{eq:2pt disc}. The main difficulty in using existing results presented in \cite{yan2018interior,worku2021simultaneous} is the positive semidefiniteness of the diffusivity tensor, $ \tilde{\K} \succeq 0$, as opposed to being symmetric positive definite, which is the assumption used to develop viscous SATs in the mentioned studies. The goal of the next subsection is to extend the results in  \cite{yan2018interior,worku2021simultaneous} to systems of equations with a symmetric positive semidefinite diffusivity tensor.

\subsection{Entropy stable viscous SATs}
To study the entropy stability of the viscous terms, we first need to write the sum of \cref{eq:visc entropy}  over all elements in terms of facet contributions. This requires that we split the viscous tensor as $ \tilde{\K} = \tilde{\K}^{\frac{1}{2}}  \tilde{\K}^{\frac{1}{2}}$, which is possible due to the following lemma, the proof of which can be found in \cite{horn2013matrix}.
\begin{lemma}\cite[Theorem 7.2.6]{horn2013matrix}\label{lem:horns lemma}
	Let $ \A\in\IR{n\times n} $ be symmetric positive semidefinite, and let $ k\in\left\{ 2,3,\dots\right\}  $, then 
	\begin{enumerate}
		\item there is a unique symmetric positive semidefinite matrix $ \B  $ such that $ \B^{k}=\A $,
		\item $ \B $ commutes with any matrix that commutes with $ \A $.
	\end{enumerate}
\end{lemma}

As a consequence of \cref{lem:horns lemma} and the assumption $ \tilde{\K}\succeq 0$, there exists a unique symmetric positive semidefinite matrix, $ \tilde{\K}^{\frac{1}{2}} \succeq 0 $, such that $ (\tilde{\K}^{\frac{1}{2}})^{2}  = \tilde{\K} $. Furthermore, defining 
\begin{equation}
	\overline{\overline{\H}}_{k} =  \I_{d}\otimes\overline{\H}_{k}=\I_{d}\otimes\H_{k}\otimes\I_{nc},
\end{equation}
and using the fact that $ \H_{k} $ is diagonal and $ \tilde{\K}_{ij,k} $ is block diagonal, it is possible to show that $ \tilde{\K}_k $ commutes with $ \overline{\overline{\H}}_{k}$, \ie, $ \overline{\overline{\H}}_{k} \tilde{\K}_k = \tilde{\K}_k  \overline{\overline{\H}}_{k} $. Hence, the second statement in \cref{lem:horns lemma} guarantees that $ \overline{\overline{\H}}_{k} \tilde{\K}_k^{\frac{1}{2}} = \tilde{\K}_k^{\frac{1}{2}}   \overline{\overline{\H}}_{k} $, which is a critical property needed to write \cref{eq:visc entropy} as a sum of surface contributions. 

\begin{lemma} \label{lem:surf viscous}
	Assuming that $ \tilde{\K} $ is symmetric positive semidefinite, the sum of the viscous terms, \cref{eq:visc entropy}, over all elements can be written as a sum of surface contributions, \ie, 
	\begin{equation} \label{eq:surf viscous}
		\begin{aligned}
			\sum_{\Omega_{k}\in\fn T_{h}}\sum_{i=1}^{d}\sum_{j=1}^{d}\bm{w}_{k}^{T}\overline{\H}_{k}\overline{\D}_{{{x}_{i}}k}&\tilde{\K}_{ij}\overline{\D}_{{x}_{j}k}\bm{w}_{k} +\sum_{\Omega_{k}\in\fn T_{h}}\bm{w}_{k}^{T}\overline{\H}_{k}\text{SAT}_{k}^{V} 
			\\&= -\sum_{\gamma\in\Gamma^{I}}\begin{bmatrix}
				\overline{\R}_{\gamma k}\bm{w}_{k}\\
				\overline{\R}_{\gamma v}\bm{w}_{v}\\
				\P_{k}\bm{w}_{k}\\
				\P_{v}\bm{w}_{v}
			\end{bmatrix}^{T}\begin{bmatrix}
				\overline{\T}_{\gamma k}^{(1)} & -\overline{\T}_{\gamma k}^{(1)} & \left(\overline{\T}_{\gamma k}^{(3)}-\overline{\B}_{\gamma}\right)\C_{\gamma k} & \overline{\T}_{\gamma k}^{(3)}\C_{\gamma v}\\
				-\overline{\T}_{\gamma v}^{(1)} & \overline{\T}_{\gamma v}^{(1)} & \overline{\T}_{\gamma v}^{(3)}\C_{\gamma k} & \left(\overline{\T}_{\gamma v}^{(3)}-\overline{\B}_{\gamma}\right)\C_{\gamma v}\\
				\C_{\gamma k}^{T}\overline{\T}_{\gamma k}^{(2)} & -\C_{\gamma k}^{T}\overline{\T}_{\gamma k}^{(2)} & \alpha_{\gamma k}\overline{\overline{\H}}_{k} & \bm{0}\\
				-\C_{\gamma v}^{T}\overline{\T}_{\gamma v}^{(2)} & \C_{\gamma v}^{T}\overline{\T}_{\gamma v}^{(2)} & \bm{0} & \alpha_{\gamma v}\overline{\overline{\H}}_{v}
			\end{bmatrix}\begin{bmatrix}
				\overline{\R}_{\gamma k}\bm{w}_{k}\\
				\overline{\R}_{\gamma v}\bm{w}_{v}\\
				\P_{k}\bm{w}_{k}\\
				\P_{v}\bm{w}_{v}
			\end{bmatrix},
		\end{aligned}
	\end{equation} 
	where $ \sum_{\gamma\in\Gamma_{k}}\alpha_{\gamma k}=1 $ and $ \sum_{\gamma\in\Gamma_{v}}\alpha_{\gamma v}=1 $ are facet weight parameters and, in three dimensions, 
	\begin{align*}
		\P_{k}&=\tilde{\K}_{k}^{\frac{1}{2}}\begin{bmatrix}
			\overline{\D}_{xk}\\
			\overline{\D}_{yk}\\
			\overline{\D}_{zk}
		\end{bmatrix}, & 
		\P_{v}&=\tilde{\K}_{v}^{\frac{1}{2}}\begin{bmatrix}
			\overline{\D}_{xv}\\
			\overline{\D}_{yv}\\
			\overline{\D}_{zv}
		\end{bmatrix},&
		\C_{\gamma k}^{T}&=\left[\tilde{\K}_{k}^{\frac{1}{2}}\right]^{T}\begin{bmatrix}
			\overline{\R}_{\gamma k}^{T}\overline{\N}_{\red{x_{1}},\gamma}\\
			\overline{\R}_{\gamma k}^{T}\overline{\N}_{\red{x_{2}},\gamma}\\
			\overline{\R}_{\gamma k}^{T}\overline{\N}_{\red{x_{3}},\gamma}
		\end{bmatrix},&
		\C_{\gamma v}^{T}&=-\left[\tilde{\K}_{v}^{\frac{1}{2}}\right]^{T}\begin{bmatrix}
			\overline{\R}_{\gamma v}^{T}\overline{\N}_{\red{x_{1}},\gamma}\\
			\overline{\R}_{\gamma v}^{T}\overline{\N}_{\red{x_{2}},\gamma}\\
			\overline{\R}_{\gamma v}^{T}\overline{\N}_{\red{x_{3}},\gamma}
		\end{bmatrix}.
	\end{align*}
\end{lemma}
\begin{proof}
	Using the SBP property and definition of the facet derivative, \cref{eq:Dgk}, the viscous volume term can be written as
	\begin{align*}
		\sum_{i=1}^{d}\sum_{j=1}^{d} \bm{w}_{k}^{T}\overline{\H}_{k}\overline{\D}_{{x}_{i}k}\tilde{\K}_{ij,k}\overline{\D}_{{x}_{j}k}\bm{w}_{k}
		&=-\sum_{i=1}^{d}\sum_{j=1}^{d}\bm{w}_{k}^{T}\overline{\D}_{{x}_{i}k}^{T}\overline{\H}_{k}\tilde{\K}_{ij,k}\overline{\D}_{{x}_{j}k}\bm{w}_{k}
		+\sum_{i=1}^{d}\sum_{j=1}^{d}\bm{w}_{k}^{T}\overline{\E}_{{x}_{i}k}\tilde{\K}_{ij,k}\overline{\D}_{{x}_{j}k}\bm{w}_{k}
		\\&=-\sum_{i=1}^{d}\sum_{j=1}^{d}\bm{w}_{k}^{T}\overline{\D}_{{x}_{i}k}^{T}\overline{\H}_{k}\tilde{\K}_{ij,k}\overline{\D}_{{x}_{j}k}\bm{w}_{k}+\sum_{\gamma\in\Gamma_{k}} \sum_{i=1}^{d}\sum_{j=1}^{d} \bm{w}_{k}^{T}\overline{\R}_{\gamma k}^{T}\overline{\B}_{\gamma}\overline{\N}_{{x}_{i}k}\overline{\R}_{\gamma k}\tilde{\K}_{ij,k}\overline{\D}_{{x}_{j}k}\bm{w}_{k}
		\\&=-\sum_{i=1}^{d}\sum_{j=1}^{d}\bm{w}_{k}^{T}\overline{\D}_{{x}_{i}k}^{T}\overline{\H}_{k}\tilde{\K}_{ij,k}\overline{\D}_{{x}_{j}k}\bm{w}_{k}
		+\sum_{\gamma\in\Gamma_{k}}\bm{w}_{k}^{T}\overline{\R}_{\gamma k}^{T}\overline{\B}_{\gamma}\overline{\D}_{\gamma k}\bm{w}_{k}.
	\end{align*}
	Applying \cref{lem:horns lemma}, the first term on the RHS of the last line can be written as 
	\begin{align*}
		-\sum_{i=1}^{d}\sum_{j=1}^{d}\bm{w}_{k}^{T}\overline{\D}_{{x}_{i}k}^{T}\overline{\H}_{k}\tilde{\K}_{ij,k}\overline{\D}_{{x}_{j}k}\bm{w}_{k}
		&= -\begin{bmatrix}
			\overline{\D}_{xk}\bm{w}_{k}\\
			\overline{\D}_{yk}\bm{w}_{k}\\
			\overline{\D}_{zk}\bm{w}_{k}
		\end{bmatrix}^{T}\overline{\overline{\H}}_{k}\tilde{\K}_{k}\begin{bmatrix}
			\overline{\D}_{xk}\bm{w}_{k}\\
			\overline{\D}_{yk}\bm{w}_{k}\\
			\overline{\D}_{zk}\bm{w}_{k}
		\end{bmatrix}
		= - \left[\begin{array}{c}
			\overline{\D}_{xk}\bm{w}_{k}\\
			\overline{\D}_{yk}\bm{w}_{k}\\
			\overline{\D}_{zk}\bm{w}_{k}
		\end{array}\right]^{T}\tilde{\K}_{k}^{\frac{1}{2}}\overline{\overline{\H}}_{k}\tilde{\K}_{k}^{\frac{1}{2}}\begin{bmatrix}
			\overline{\D}_{xk}\bm{w}_{k}\\
			\overline{\D}_{yk}\bm{w}_{k}\\
			\overline{\D}_{zk}\bm{w}_{k}
		\end{bmatrix}
		= - \bm{w}_{k}^{T}\P_{k}^{T}\overline{\overline{\H}}_{k}\P_{k}\bm{w}_{k},
	\end{align*}
	hence we have 
	\begin{equation}\label{eq:vol visc}
		\sum_{i=1}^{d}\sum_{j=1}^{d} \bm{w}_{k}^{T}\overline{\H}_{k}\overline{\D}_{{x}_{i}k}\tilde{\K}_{ij,k}\overline{\D}_{{x}_{j}k}\bm{w}_{k}
		= - \bm{w}_{k}^{T}\P_{k}^{T}\overline{\overline{\H}}_{k}\P_{k}\bm{w}_{k} +\sum_{\gamma\in\Gamma_{k}}\bm{w}_{k}^{T}\overline{\R}_{\gamma k}^{T}\overline{\B}_{\gamma}\overline{\D}_{\gamma k}\bm{w}_{k},
	\end{equation}
	and a similar result holds for the viscous volume term of element $ \Omega_{v} $. We also have that
	\begin{align}
		\bm{w}_{k}^{T}\overline{\R}_{\gamma k}^{T}\overline{\T}_{\gamma k}^{(3)}\overline{\D}_{\gamma k}\bm{w}_{k}
		&=\sum_{i=1}^{d}\sum_{j=1}^{d} \bm{w}_{k}^{T}\overline{\R}_{\gamma k}^{T}\overline{\T}_{\gamma k}^{(3)}\overline{\N}_{{x}_{i}\gamma}\overline{\R}_{\gamma k}\tilde{\K}_{ij,k}\overline{\D}_{{x}_{j}k}\bm{w}_{k}=\bm{w}_{k}^{T}\overline{\R}_{\gamma k}^{T}\overline{\T}_{\gamma k}^{(3)}\begin{bmatrix}
			\overline{\R}_{\gamma k}^{T}\overline{\N}_{\red{x_{1}},\gamma}\\
			\overline{\R}_{\gamma k}^{T}\overline{\N}_{\red{x_{2}},\gamma}\\
			\overline{\R}_{\gamma k}^{T}\overline{\N}_{\red{x_{3}},\gamma}
		\end{bmatrix}^{T}\tilde{\K}_{k}^{\frac{1}{2}}\tilde{\K}_{k}^{\frac{1}{2}}\begin{bmatrix}
			\overline{\D}_{xk}\bm{w}_{k}\\
			\overline{\D}_{yk}\bm{w}_{k}\\
			\overline{\D}_{zk}\bm{w}_{k}
		\end{bmatrix} \nonumber
		\\&=\bm{w}_{k}^{T}\overline{\R}_{\gamma k}^{T}\overline{\T}_{\gamma k}^{(3)}\C_{\gamma k}\P_{k}\bm{w}_{k}\label{eq:T3 surf}
	\end{align}
	Similarly, we can write
	\begin{equation}\label{eq:T2 surf}
		\bm{w}_{k}^{T}\overline{\D}_{\gamma k}^{T}\overline{\T}_{\gamma k}^{(2)}\overline{\R}_{\gamma k}\bm{w}_{k}=\bm{w}_{k}^{T}\P_{k}^{T}\C_{\gamma k}^{T}\overline{\T}_{\gamma k}^{(2)}\overline{\R}_{\gamma k}\bm{w}_{k}.
	\end{equation}
	Finally, substituting \cref{eq:vol visc}, \cref{eq:T3 surf}, and \cref{eq:T2 surf} into the LHS of \cref{eq:surf viscous}, writing the result as sums over interfaces and rearranging gives the RHS of \cref{eq:surf viscous}.
\end{proof}

\cref{lem:surf viscous} allows us to write the viscous terms in a form suitable for entropy stability analyses. In fact, the conditions obtained in \cite{yan2018interior,worku2021simultaneous} for the coefficient matrices to produce an energy stable discretization of scalar diffusion problems also apply for the \red{compressible} Navier-Stokes equations upon modifying the coefficient matrices for systems of equations. Starting from \cref{eq:surf viscous}, the entropy stability analyses for the viscous portion of the \red{compressible} Navier-Stokes equations follow similar steps as in \cite{yan2018interior,worku2021simultaneous}; hence, we present the following theorem without a proof.
\begin{theorem}\label{thm:entropy stability nse}
	Let \cref{ass:mapping} hold, the viscous SAT coefficient matrices, $ \overline{\T}_{\gamma k}^{(j)} $, $ j\in\{1,\dots,3\} $, be symmetric, and the diffusivity tensor be symmetric positive semidefinite, $ \tilde{\K}\succeq 0 $. Then, the entropy-split and Hadamard-form SBP-SAT discretizations of the \red{compressible} Navier-Stokes equations are entropy stable if 
	\begin{align}\label{eq:adjoint inconsistent stable sats}
		\overline{\T}_{\gamma k}^{(3)}+\overline{\T}_{\gamma k}^{(2)}-\overline{\B}_{\gamma} & =\overline{\T}_{\gamma v}^{(3)}+\overline{\T}_{\gamma v}^{(2)}-\overline{\B}_{\gamma}=\bm{0}, & \overline{\T}_{\gamma v}^{(3)}-\overline{\T}_{\gamma k}^{(2)} & =\overline{\T}_{\gamma k}^{(3)}-\overline{\T}_{\gamma v}^{(2)}=\bm{0}, & \overline{\T}_{\gamma k}^{(1)} =\overline{\T}_{\gamma v}^{(1)} & \succeq0.
	\end{align}
	Furthermore, the viscous SATs in \cref{eq:added term} are conservative if	
	\begin{align}\label{eq:visc cons}
		\overline{\T}_{\gamma k}^{(1)}=\overline{\T}_{\gamma v}^{(1)},
		\qquad \overline{\T}_{\gamma k}^{(3)}+\overline{\T}_{\gamma v}^{(3)}&=\overline{\B}_{\gamma}.
	\end{align}
	Alternatively, if instead of the conditions in \cref{eq:adjoint inconsistent stable sats},  $ \overline{\T}_{\gamma k}^{(2)}+\overline{\T}_{\gamma v}^{(2)}=-\overline{\B}_{\gamma} $, $ \overline{\T}_{\gamma k}^{(3)} - \overline{\T}_{\gamma v}^{(2)}=\overline{\B}_{\gamma} $, and \cref{eq:visc cons} hold, then the SBP-SAT discretization of the \red{compressible} Navier-Stokes equations is entropy stable provided that 
	\begin{align}
		\overline{\T}_{\gamma k}^{(1)}- \left(\frac{1}{\alpha_{\gamma k}}\overline{\T}_{\gamma k}^{(2)}\C_{\gamma k}\overline{\overline{\H}}_{k}^{-1}\C_{\gamma k}^{T}\overline{\T}_{\gamma k}^{(2)}+\frac{1}{\alpha_{\gamma v}}\overline{\T}_{\gamma v}^{(2)}\C_{\gamma v}\overline{\overline{\H}}_{v}^{-1}\C_{\gamma v}^{T}\overline{\T}_{\gamma v}^{(2)}\right) \succeq0.
	\end{align}
\end{theorem}
We note that the proof of the conservation conditions given in \cref{eq:visc cons} follows directly from the conservation analysis in \cite{worku2021simultaneous}. Furthermore, the conditions $ \overline{\T}_{\gamma k}^{(2)}+\overline{\T}_{\gamma v}^{(2)}=-\overline{\B}_{\gamma} $ and  $ \overline{\T}_{\gamma k}^{(3)} - \overline{\T}_{\gamma v}^{(2)}=\overline{\B}_{\gamma} $ are satisfied by many types of adjoint consistent viscous SATs, including the BR2 SAT for which 
\begin{align}\label{eq:br2 sats}
	\overline{\T}_{\gamma k}^{(1)} = \frac{1}{4}\overline{\B}_\gamma\left(\frac{1}{\alpha_{\gamma k}}\C_{\gamma k}\overline{\overline{\H}}_{k}^{-1}\C_{\gamma k}^{T}+\frac{1}{\alpha_{\gamma v}}\C_{\gamma v}\overline{\overline{\H}}_{v}^{-1}\C_{\gamma v}^{T}\right) \overline{\B}_\gamma,
	\qquad 
	\overline{\T}_{\gamma k}^{(2)} = -\overline{\T}_{\gamma k}^{(3)} = - \frac{1}{2}\overline{\B}_\gamma.
\end{align} 
Computing the coefficient $ \overline{\T}_{\gamma k}^{(1)}  $ for the BR2 SAT is expensive; hence, the SIPG method has been proposed as an alternative in \cite{yan2018interior}. These SATs are adjoint consistent and lead to functional superconvergence for scalar diffusion problems. A less expensive but adjoint inconsistent alternative is the Baumann-Oden SAT given by 
\begin{alignat}{2}\label{eq:bo sats}
	\overline{\T}_{\gamma k}^{(1)} &=\overline{\T}_{\gamma v}^{(1)} = \bm{0}, \qquad & \overline{\T}_{\gamma k}^{(3)}&=\overline{\T}_{\gamma v}^{(3)}=\overline{\T}_{\gamma k}^{(2)}=\overline{\T}_{\gamma v}^{(2)}=\frac{1}{2}\overline{\B}_{\gamma}.
\end{alignat}
In addition to its efficiency, the BO SAT possesses interesting properties such as smaller spectral radius and condition number for scalar diffusion problems \cite{worku2021simultaneous}. Hence, we will use the BO SAT for the numerical experiments in this work. Furthermore, we set the facet weight parameters to a constant, $ \alpha_{\gamma k}=\alpha_{\gamma v}=1/n_f$, where $ n_f $ is the number of facets of the element. Penalty coefficients matrices corresponding to many other types of viscous SATs are reported in \cite{worku2021simultaneous}. We remind the reader that \cref{thm:entropy stability nse} does not imply entropy stability of the entropy-split or Sj{\"o}green-Yee Hadamard-form discretizations of the \red{compressible} Navier-Stokes equations unless heat fluxes are neglected. For Hadamard-form discretizations based on the generalized entropy function, \cref{eq:generalized entropy}, however, \cref{thm:entropy stability nse} guarantees entropy stability. 

\section{Some implementation details}\label{sec:implementation details}
\subsection{Fully entropy conservative or dissipative scheme}
To extend the entropy conservation or stability of the entropy-split and Hadamard-form semi-discrete schemes presented in the previous sections to the time dimension, we use the explicit fourth-order relaxation Runge-Kutta method \cite{ranocha2020relaxation,ketcheson2019relaxation,ranocha2020fully,dekker1984stability,del2002explicit,calvo2006preservation}. A fully entropy conservative or stable scheme is obtained by solving a scalar nonlinear equation at every timestep such that the entropy remains conserved or bounded in time, \ie, we solve for the parameter $ \gamma_n $ such that
\begin{equation}\label{eq:rrk4 gamma_n}
	\bm{1}^T\H_{g}\bm{s}\left(\bm{u}_n + \gamma_n \Delta t \sum_{i=1}^{s} \bm{b}_{i} \bm{r}_i\right) - \bm{1}^T\H_{g}\bm{s}(\bm{u}_n) - \gamma_n \Delta t \sum_{i=1}^{s}\bm{b}_i \bm{w}_i^{T}\H_{g}\bm{r}_i \le 0,
\end{equation}
where $ \bm{r} $ is the global residual vector from the spatial discretization, $ \bm{u} $ and $ \bm{w} $ are the global conservative and entropy variable vectors, respectively, $ \H_{g} $ is the global norm matrix, $ \bm{s}(\cdot)$ indicates the evaluation of the global entropy function vector at the state in the parentheses, $ \bm{b} $ is a vector of the RK4 quadrature weights of the last stage, the subscript $ n $ indicates the current timestep, $ s $ is the number of stages, and $ \Delta t $ is the timestep. The global vectors are obtained by concatenating the corresponding local (elementwise) vectors. For entropy conservative schemes, the inequality in \cref{eq:rrk4 gamma_n} is replaced by an equality. We solve \cref{eq:rrk4 gamma_n} for $ \gamma_n $ to machine precision using the bisection method, and update the solution as
\begin{equation}\label{eq:rrk_update}
	\bm{u}_{n+1} = \bm{u}_n + \gamma_n \Delta t \sum_{i=1}^{s} \bm{b}_i \bm{r}_i.
\end{equation}

\subsection{\blue{Efficient forms} of the Hadamard-form and entropy-split entropy stable discretizations}
Using the SBP property, it is possible to write the Hadamard-form scheme, \cref{eq:2pt disc}, for the \red{compressible} Navier-Stokes equations in the weak-form \cite{shadpey2020entropy,crean2018entropy}
\begin{equation} \label{eq:weak form 2pt}
	\begin{aligned}
		\der[\bm{u}_{k}]t+\overline{\H}_{k}^{-1}\sum_{i=1}^{d} \left[\overline{\S}_{{x}_{i}k}\circ\F_{i}\left(\bm{u}_{k},\bm{u}_{k}\right)\right]\bm{1} &=-\frac{1}{2}\overline{\H}_{k}^{-1}\sum_{\gamma\in\Gamma_{k}^{I}}\sum_{i=1}^{d}\bm{w}_{k}^{T}\left[\overline{\E}_{{x}_{i}}^{kv}\circ\F_{i}\left(\bm{u}_{k},\bm{u}_{v}\right)\right]\bm{1}
		\\&\quad
		+\sum_{i=1}^{d}\sum_{j=1}^{d}\overline{\D}_{{x}_{i}k}\tilde{\K}_{ij,k}\overline{\D}_{{x}_{j}k}\bm{w}_{k}+\text{SAT}_{k}^{V},
	\end{aligned}
\end{equation}
where the viscous SATs are given in \cref{eq:sats visc}. The weak-form is more efficient for implementation than the strong-form discretization, \cref{eq:2pt disc}, as one of the inviscid SAT terms is added as a volume term once per element instead of $ n_f $ times. Note that writing the viscous terms in a  weak-form does not offer efficiency benefits. Using the SBP property, $ \overline{\Q}_{\bm{x}_i k}+ \overline{\Q}_{\bm{x}_i k}^T = \overline{\E }_{\bm{x}_i k}$, it is straightforward to show that the entropy-split scheme, \cref{eq:discrete entropy split}, for the \red{compressible} Navier-Stokes equations can also be written in \blue{a more efficient, but algebraically equivalent, form} as
\begin{equation}\label{eq:weak form split}
	\begin{aligned}
		\der[\bm{u}_{k}]t-\overline{\H}_{k}^{-1}\sum_{i=1}^{d}\left[\frac{1}{\beta+1}\overline{\Q}_{{x}_{i}k}+\overline{\Q}_{{x}_{i}k}^{T}\right]\bm{f}_k^{I(i)}+\frac{1}{\beta+1}\sum_{i=1}^{d}\tilde{\A}_{i}\overline{\D}_{{x}_{i}k}\bm{w}_{k}
		&=-\frac{1}{2}\overline{\H}_{k}^{-1}\sum_{\gamma\in\Gamma^{I}}\sum_{i=1}^{d}\left[\overline{\E}_{{x}_{i}}^{kv}\circ\F_{i}\left(\bm{u}_{k},\bm{u}_{v}\right)\right]\bm{1}
		\\&\quad +\sum_{i=1}^{d}\sum_{j=1}^{d}\overline{\D}_{{x}_{i}k}\tilde{\K}_{ij,k}\overline{\D}_{{x}_{j}k}\bm{w}_{k}+\text{SAT}_{k}^{V}.
	\end{aligned}
\end{equation}
For the numerical experiments in this work, we use the more efficient forms of the Hadamard-form and entropy-split entropy stable discretizations given in \cref{eq:weak form 2pt} and \cref{eq:weak form split}. 
\subsection{Efficiency considerations}\label{subsec:efficiency considerations}
As outlined in the introduction, one of the goals of this paper is to understand the efficiency benefits of the entropy-split scheme compared to the Hadamard-form discretizations of the Euler and \red{compressible} Navier-Stokes equations on unstructured grids. To this end, we implement both schemes in a single framework such that only the computation of the inviscid flux derivatives are handled in separate subroutines. Furthermore, we compare the times taken to compute the various components of the spatial discretizations separately, \eg, the time taken by the volume and facet terms or the inviscid and viscous SATs is broken down by subroutine. The efficiency comparison excludes the time required to construct the SBP operators or advance the state in time as the focus is on the effects of splitting the inviscid fluxes versus using the Hadamard form. \red{All efficiency comparisons are performed in serial and in double precision on the compute nodes of the Niagara supercomputer at the SciNet HPC consortium. Details of the computer architecture can be found in \cite{ponce2019deploying}. All of the numerical experiments are conducted using the Automated PDE Solver (APS) of the Aerospace Computational Engineering Lab at the University of Toronto. APS is written in C++, and was compiled with the optimization level $ \texttt{-O3} $.}

We have implemented many of the recommendations provided in \cite{ranocha2021efficient} to improve the efficiency of the Hadamard-form scheme. These include using the  weak-form, computing only the upper triangular entries of $ \F_{i} $  (making use of its symmetry and the fact that the diagonal entries of the skew symmetric matrix in \cref{eq:weak form 2pt}, $ \overline{\Q}_{{x}_{i}k}- \frac{1}{2}\overline{\E}_{{x}_{i}k} = \overline{\S}_{{x}_ik}$, are zero), computing the SATs only once per facet, and inlining repeatedly called functions.\ignore{We have not implemented the more invasive optimization techniques such as precomputing certain variables or applying SIMD (Single-Instruction-Multiple-Data) vectorization. Furthermore,} Unlike \cite{ranocha2021efficient}, the current study uses multidimensional SBP operators; hence, we do not have a tensor-product structure that can be exploited for efficiency, which has a significant impact on the performance of the schemes considered. Assuming an equal number of nodes per element, say $ n_p = (p+1)^d $, the number of function calls to compute $ \F_{i} $ is asymptotically $ (p+1)^{d-1} $ times more expensive for multidimensional operators compared to tensor-product operators. Similarly, the number of floating point operations without flux evaluation to compute $ \overline{\S}_{\bm{x}_i k} \circ \F_{i}\left(\bm{u}_{k},\bm{u}_{k}\right) $ is $ (p+1)^{d-1} $ times more expensive when using multidimensional operators. Hence, in three dimensions the difference in computational cost can be in orders of magnitude even for modest values of $ p $. One might think that the number nodes of multidimensional diagonal-$ \E $ operators is much fewer compared to the number of nodes of tensor-product elements of equal degree, but currently available diagonal-$ \E $ operators have a comparable number of nodes as tensor-product elements, at least for relatively low degrees. For example, in 3D, the $ p=3 $ multidimensional SBP diagonal-$ \E $ operator has $ 69 $ nodes which is larger than the $ 64 $ nodes on the same degree LGL tensor-product element, and the $ p=4 $ operator has $ 99 $ nodes while the corresponding tensor-product element has $ 125 $. It might be possible to construct SBP diagonal-$ \E $ operators with substantially fewer nodes; we hope to investigate this in the future. Other types of SBP operators, \eg, SBP-$ \Omega $, have fewer nodes, but it is not straightforward to construct entropy conservative entropy-split schemes on non-diagonal $ \E $ operators. 

\subsection{Robustness}
The stability and accuracy of the entropy-split scheme depends on the value of the arbitrary split parameter, $ \beta $, as reported in many studies, \eg, see \cite{yee2000entropy,sandham2002entropy,sjogreen2019entropy,yee2023recent}.  Since $ \beta $ is problem dependent, the entropy-split scheme is, by definition, less robust than the Hadamard-form entropy stable scheme with two-point fluxes based on the generalized entropy function, \cref{eq:generalized entropy}. In general, $ \beta  $ values between $ 1.5 <\beta < 2.5 $ are recommended \cite{sjogreen2019entropy}, although lower or higher values have also been used frequently in the literature. \red{Unless stated otherwise, we set $ \beta=2.5$ for all our numerical experiments, which is a value used in many studies, \eg, \cite{sjogreen2021construction,yee2023recent}, and that we found to work well for the test cases considered.} A thorough study of the $ \beta $-sensitivity of the benchmark problems is out of the scope of the present article, \red{but some studies in this regard can be found in \cite{yee2000entropy}}.

As discussed earlier, another downside of the entropy-split scheme is that it is not provably entropy stable for the \red{compressible} Navier-Stokes equations unless heat fluxes are neglected. However, the method might provide improved stability properties at least for certain types of problems. Mapping the applicability of the entropy-split schemes to types of common fluid dynamics problems is left for future studies. Finally, even though we have developed a hybrid scheme that can be used to locally enforce conservation for problems with discontinuities, it is not clear whether there are other detrimental consequence of the loss of conservation for other types of problems. In conclusion, more investigation is required to fully understand and apply the entropy-split scheme for a wide variety of problems. That being said, however, we show, via numerical examples, that the entropy-split scheme has the potential to be almost as robust as Hadamard-form multidimensional SBP-SAT discretizations, but with substantially less computational cost. 

\section{Numerical results}\label{sec:results}
We apply the entropy-split scheme to a number of benchmark test cases and investigate its properties. In figures and tables, the following acronyms are used: ``ES" for the entropy-split scheme, ``EH" for the hybrid method, and ``IR" for the Ismail-Roe Hadamard-form discretization \cite{ismail2009affordable}. Unless specified otherwise, time integration is handled using the fourth-order relaxation Runge-Kutta method \cite{ranocha2020relaxation}. 

\subsection{The Sod shock tube problem}
We consider the Sod shock tube problem \cite{sod1978survey}, which is governed by the Euler equations, on $ \Omega=[0,1] $ with the following initial conditions,
\begin{alignat*}{5}
	\rho_{L} & =1, \quad& u_{L} & =0, \quad& p_{L} & =1, \quad& \text{if }x &\le0.5,\\
	\rho_{R} & =0.125, \quad & u_{R} & =0, \quad& p_{R} & =0.1, \quad& \text{if }x&>0.5.
\end{alignat*}
The problem is run until $ t_f = 0.2 $ and the Courant-Friedrichs-Lewy (CFL) number is set to 0.3. The timestep is taken to be the minimum of the timesteps computed for each element as 
\begin{equation}
	\Delta t  = \text{CFL} \frac{{h_{\min}}}{|\lambda|_{\max}},
\end{equation}
where $  \abs{\lambda}_{\max}  = \max(\violet{[\sqrt{\bm{V}\cdot\bm{V}} + \sqrt{\gamma p/\rho}]_{1}, \dots, [\sqrt{\bm{V}\cdot\bm{V}} + \sqrt{\gamma p/\rho}]_{n_{pk}}})$ and $  h_{\min} =\min_{a, b \in S_k} \norm{a - b}_2$. \red{The selected timestep is used to compute the relaxation parameter, $ \gamma_{n} $, and march the solution in time according to \cref{eq:rrk_update}.} This test case is well-suited to study the effect of the loss of conservation of the entropy-split scheme by comparing the location of the moving contact and shock discontinuities against the analytical solution, which is calculated using an iterative procedure (\eg, see \cite{pulliam2014fundamental}). Since the exact solution is accessible, we apply the hybrid scheme to a maximum of three closest elements to the contact and shock discontinuities at all times. The matrix-form \blue{interface} dissipation operator developed in \cref{subsec:interface dissipation} is applied for both the entropy-split and hybrid schemes. 
\begin{figure}[!t]
	\centering
	\begin{subfigure}{0.5\textwidth}
		\centering
		\includegraphics[scale=0.31]{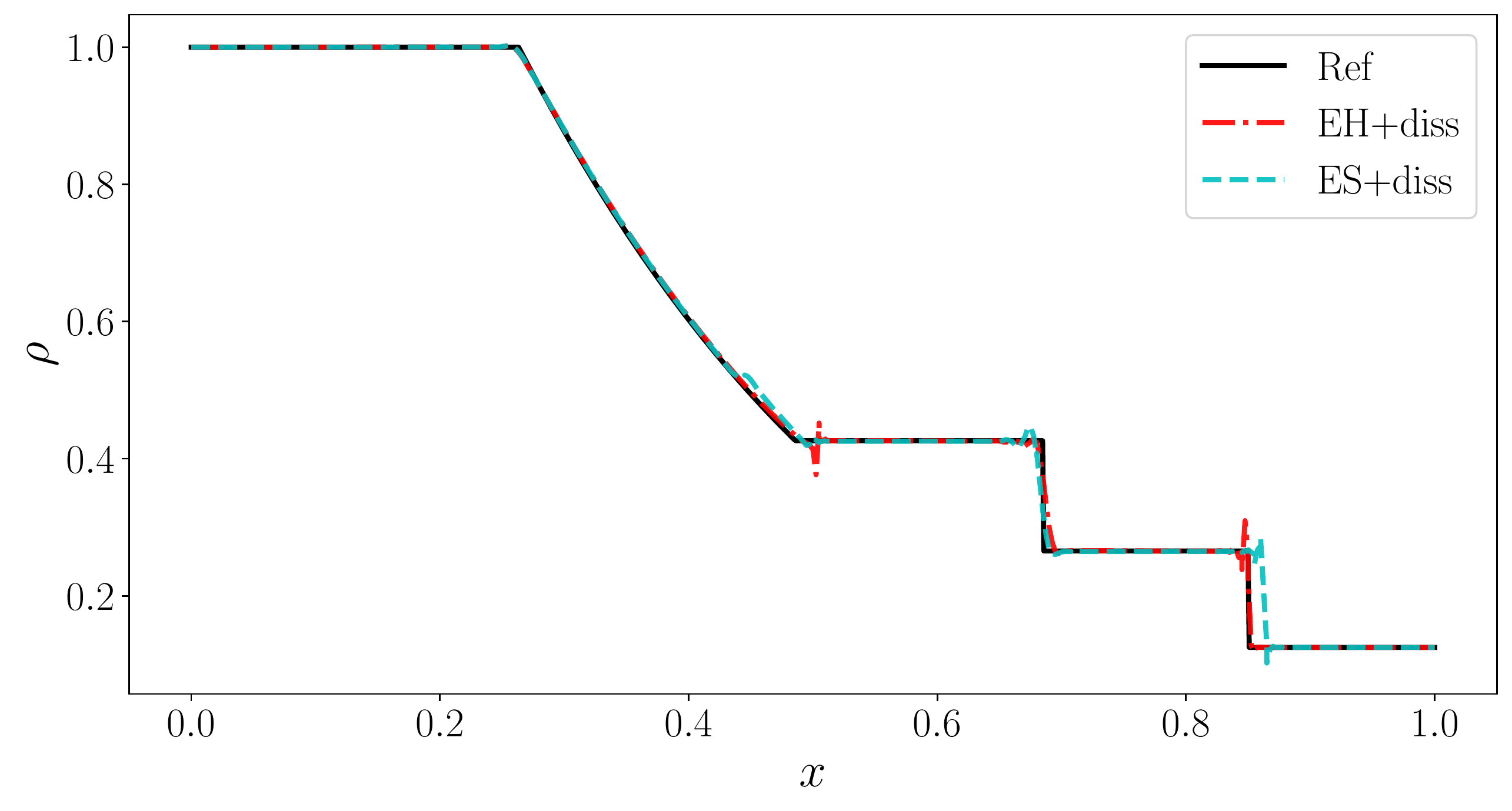}
		\caption{\label{fig:sod_coarse} $ p=2 $, $ n_e = 200 $}
	\end{subfigure}\hfill
	\begin{subfigure}{0.5\textwidth}
		\centering
		\includegraphics[scale=0.31]{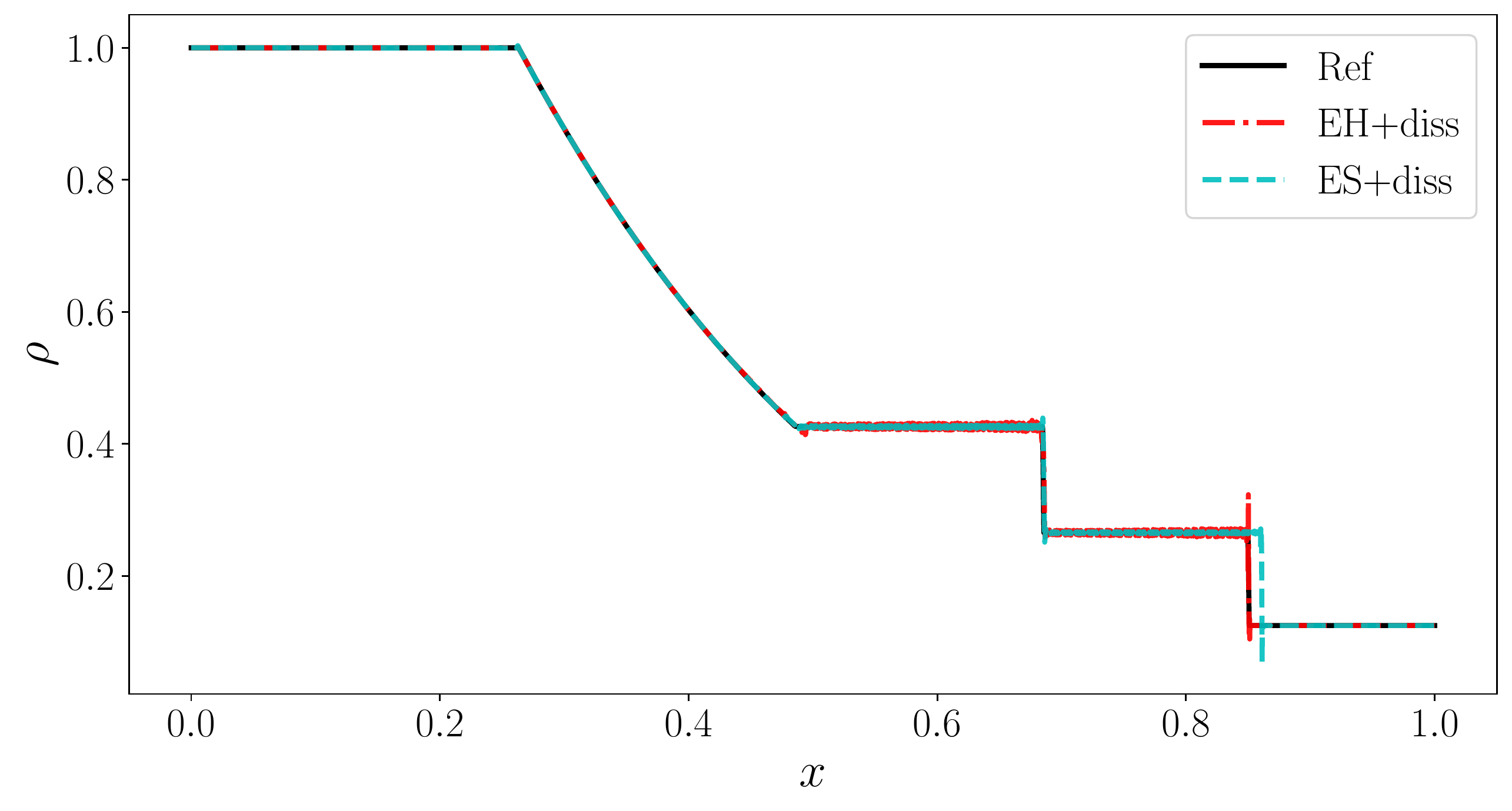}
		\caption{\label{fig:sod_fine} $ p=5 $, $ n_e = 800 $}
	\end{subfigure} 
	\caption{\label{fig:sod} Density profile for the Sod shock tube problem at $ t=0.2 $ using the locally conservative hybrid scheme and the nonconservative entropy-split scheme.}
\end{figure}

\begin{figure}[!t]
	\centering
	\begin{subfigure}{0.5\textwidth}
		\centering
		\includegraphics[scale=0.52]{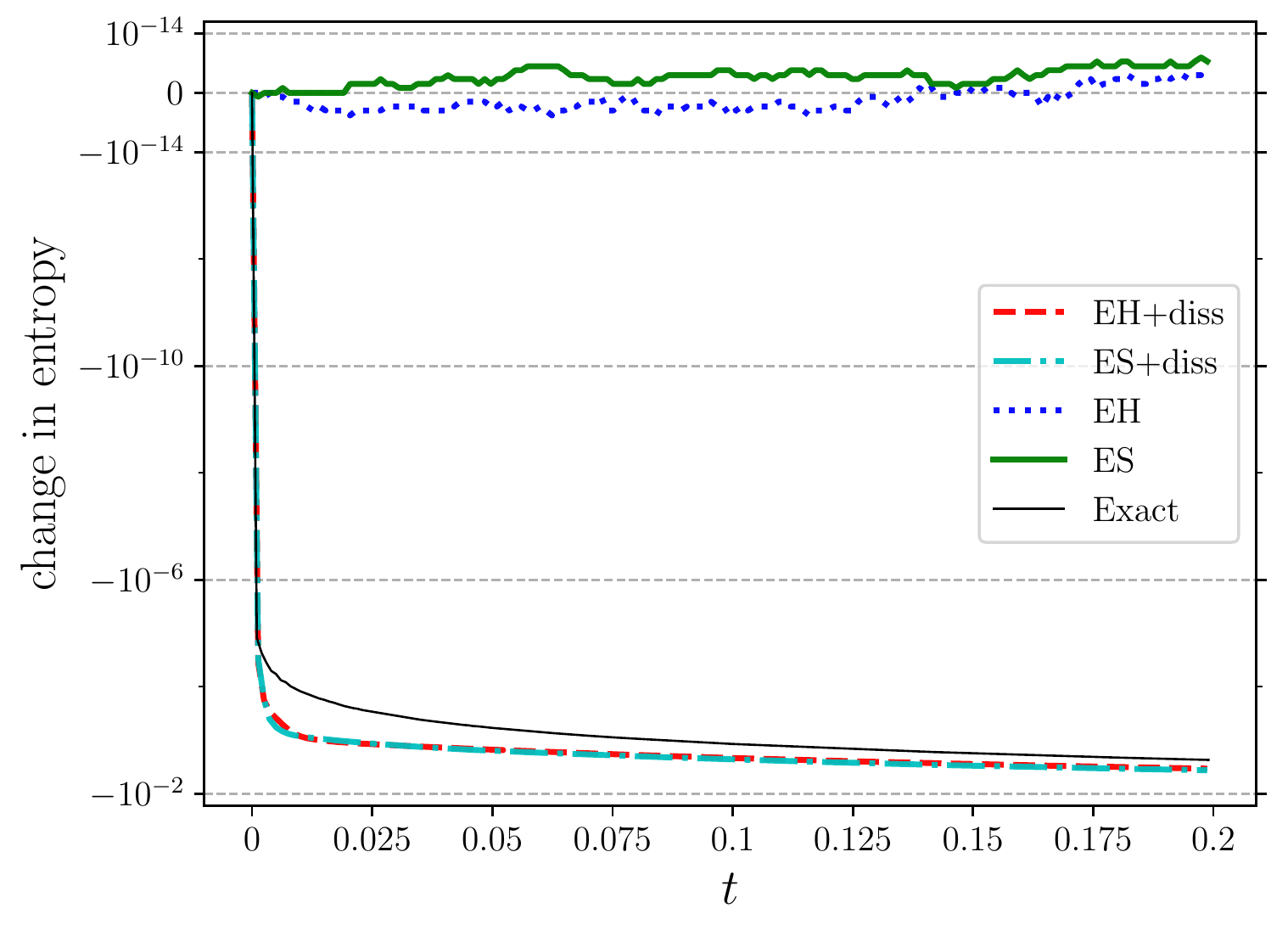}
		\caption{\label{fig:sod_entropy} $ p=1 $, $ n_e = 200 $}
	\end{subfigure}\hfill
	\begin{subfigure}{0.5\textwidth}
		\centering
		\includegraphics[scale=0.52]{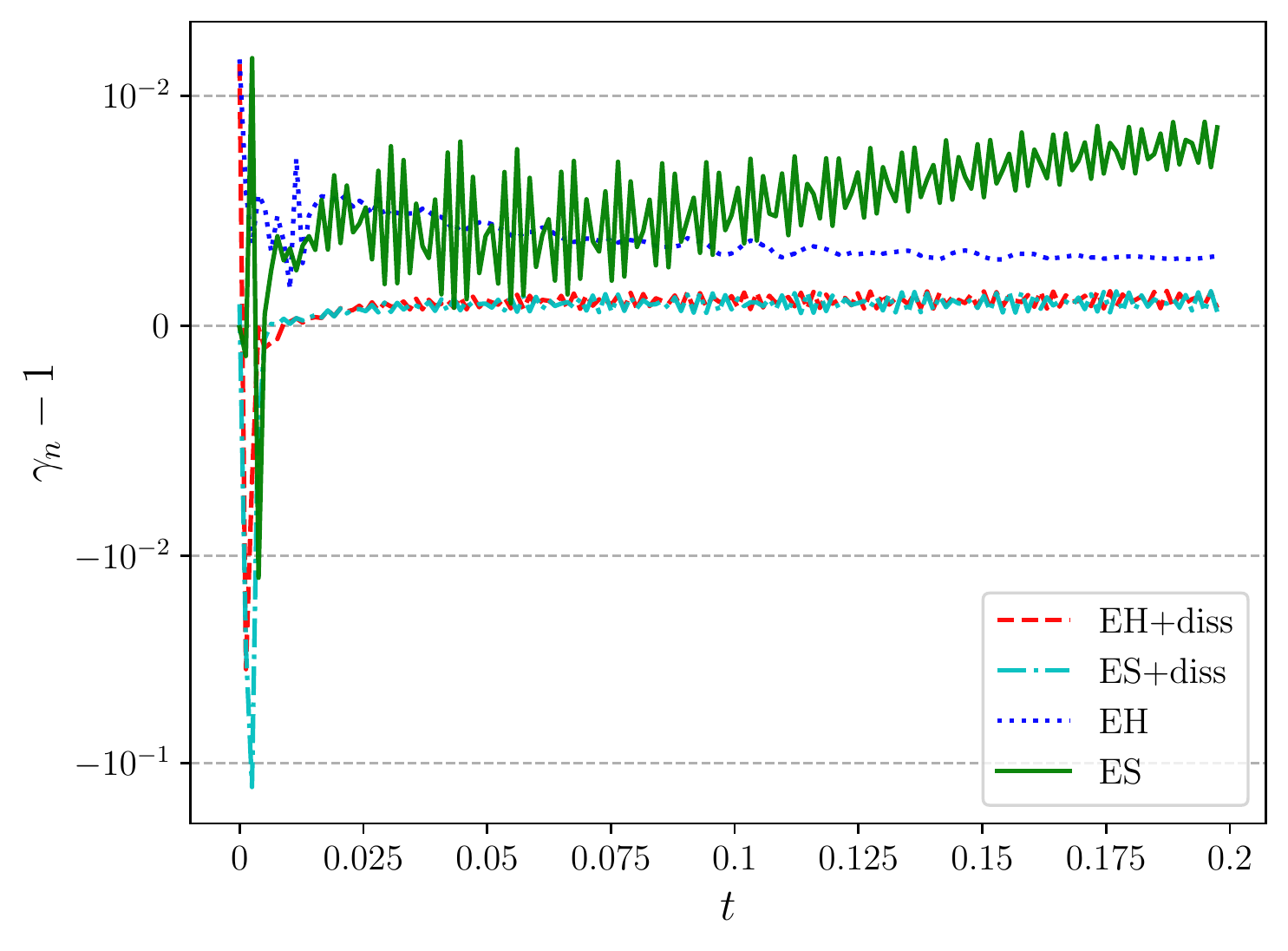}
		\caption{\label{fig:sod_gamma} $ p=1 $, $ n_e = 200 $}
	\end{subfigure} 
	\caption{\label{fig:entropy_change} Change in entropy for the Sod shock tube problem at $ t=0.2 $ using the hybrid and entropy-split schemes with and without dissipation.}
\end{figure}

\cref{fig:sod} shows that the hybrid scheme captures the correct shock speed while the entropy-split method does not. It is evident from \cref{fig:sod_fine} that the entropy-split solution does not converge to the reference solution by refining the mesh or increasing the order of discretization. The numerical solutions oscillate as no limiting or explicit shock capturing methods are employed. \cref{fig:sod_entropy} shows that both the entropy-split and hybrid schemes conserve entropy for problems with shock waves in the absence of \blue{interface} dissipation. It can also be concluded from the figure that the amount of dissipation added using the \blue{interface} dissipation operator, \cref{eq:dissipation operator}, is close to the exact entropy dissipation for this problem. The change in entropy is computed as
\begin{equation}
	\Delta s_t = \frac{s_t - s_0}{\abs{s_0}},
\end{equation}
where $ s_t $ and $ s_0 $ are the entropy values at time $ t $ and $ t=0 $ computed as $ s = \bm{1}^T\H_{g}\bm{s} $. \cref{fig:sod_gamma} shows that the $ \gamma_n $ parameter of the RRK method is close to unity for both schemes with or without dissipation.

\subsection{The Shu-Osher shock tube problem}
To further demonstrate that the hybrid scheme captures discontinuities accurately, we solve the Shu-Osher sine wave problem \cite{shu1989efficient}, which simulates a normal shock wave propagating into a medium with small sinusoidal density perturbations. The problem is governed by the Euler equations, and the initial conditions are 
\begin{alignat*}{5}
	\rho_{L} & =3.857143, \quad& u_{L} & =2.629369, \quad& p_{L} & =10.33333, \quad& \text{if }x &\le-4,\\
	\rho_{R} & =1+0.2\sin(5x), \quad & u_{R} & =0, \quad& p_{R} & =1, \quad& \text{if }x&>-4.
\end{alignat*}
The domain is $ \Omega=[-5,5] $, the final time is $ t_{f}=1.8 $, and the CFL value is set to $ 0.1 $. \red{The timestep computation and solution updates are handled in the same manner as for the Sod shock tube problem.} For this problem, the Sj{\"o}green-Yee two point flux is applied on all elements between the left and right most shock waves, \blue{whose locations are determined from numerical solutions obtained using conservative schemes}. At $ t=0 $, for instance, the conservative scheme is applied on a single or two elements at $ x=-4 $, while at $ t=1.8 $ it is applied on elements between $ x=-2.8 $ and $ x=2.5 $, \blue{covering the region of the domain in which discontinuities are present}. The idea is to apply the two-point flux scheme only at critical portions of the domain where discrete conservation is indispensable. Since the exact solution for this problem is not known, we use the $ p=3 $ Ismail-Roe Hadamard-form discretization with $ n_e = 10000 $ to compute the reference solution. A matrix-type dissipation operator of \cite{ismail2009affordable} is applied with the IR discretization, while the entropy-split and hybrid schemes use the dissipation operator presented in \cref{subsec:interface dissipation}.

\begin{figure}[!t]
	\centering
	\begin{subfigure}{0.5\textwidth}
		\centering
		\includegraphics[scale=0.31]{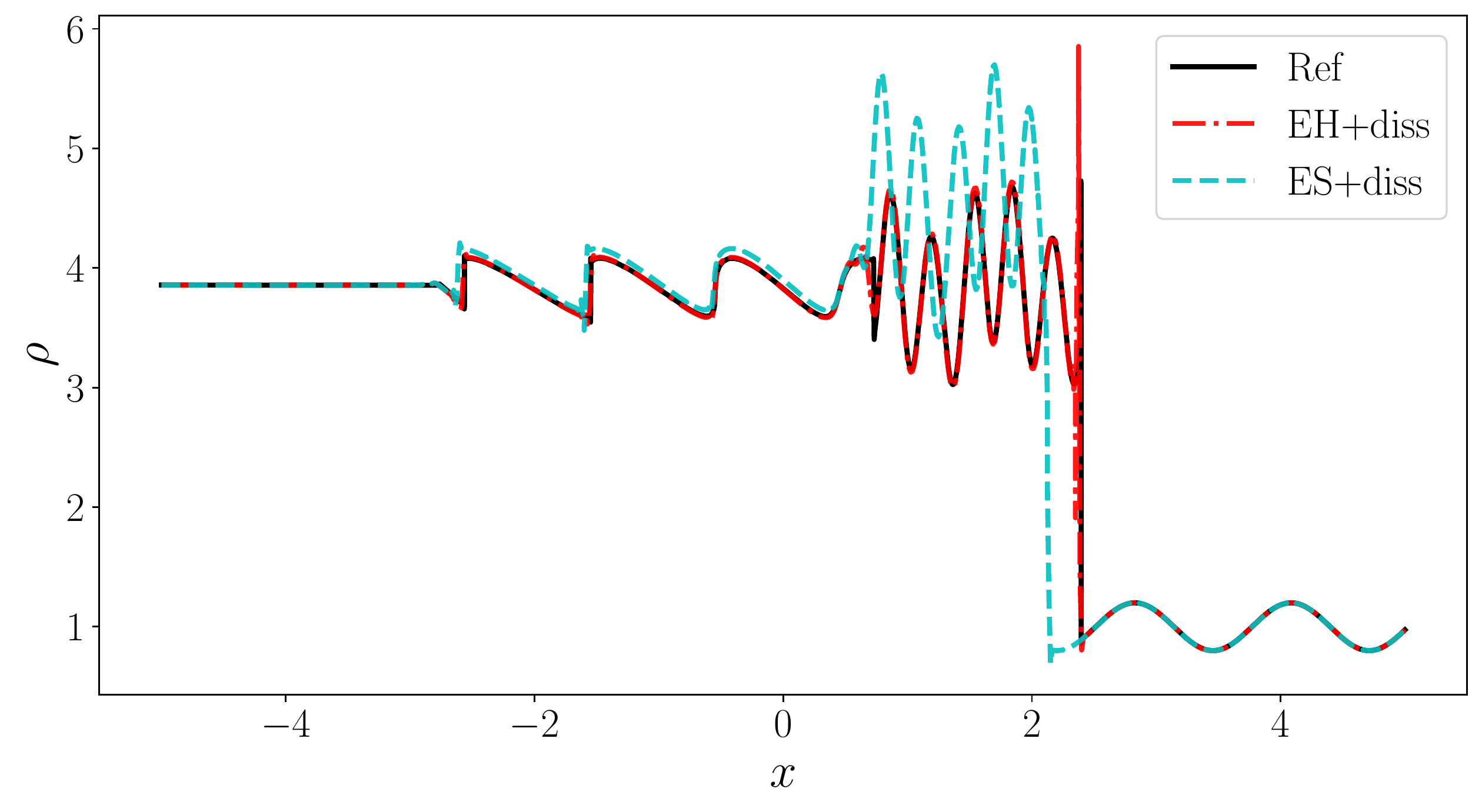}
		\caption{\label{fig:shu_osher_fine} $ p=2 $, $ n_e = 400 $}
	\end{subfigure}\hfill
	\begin{subfigure}{0.5\textwidth}
		\centering
		\includegraphics[scale=0.31]{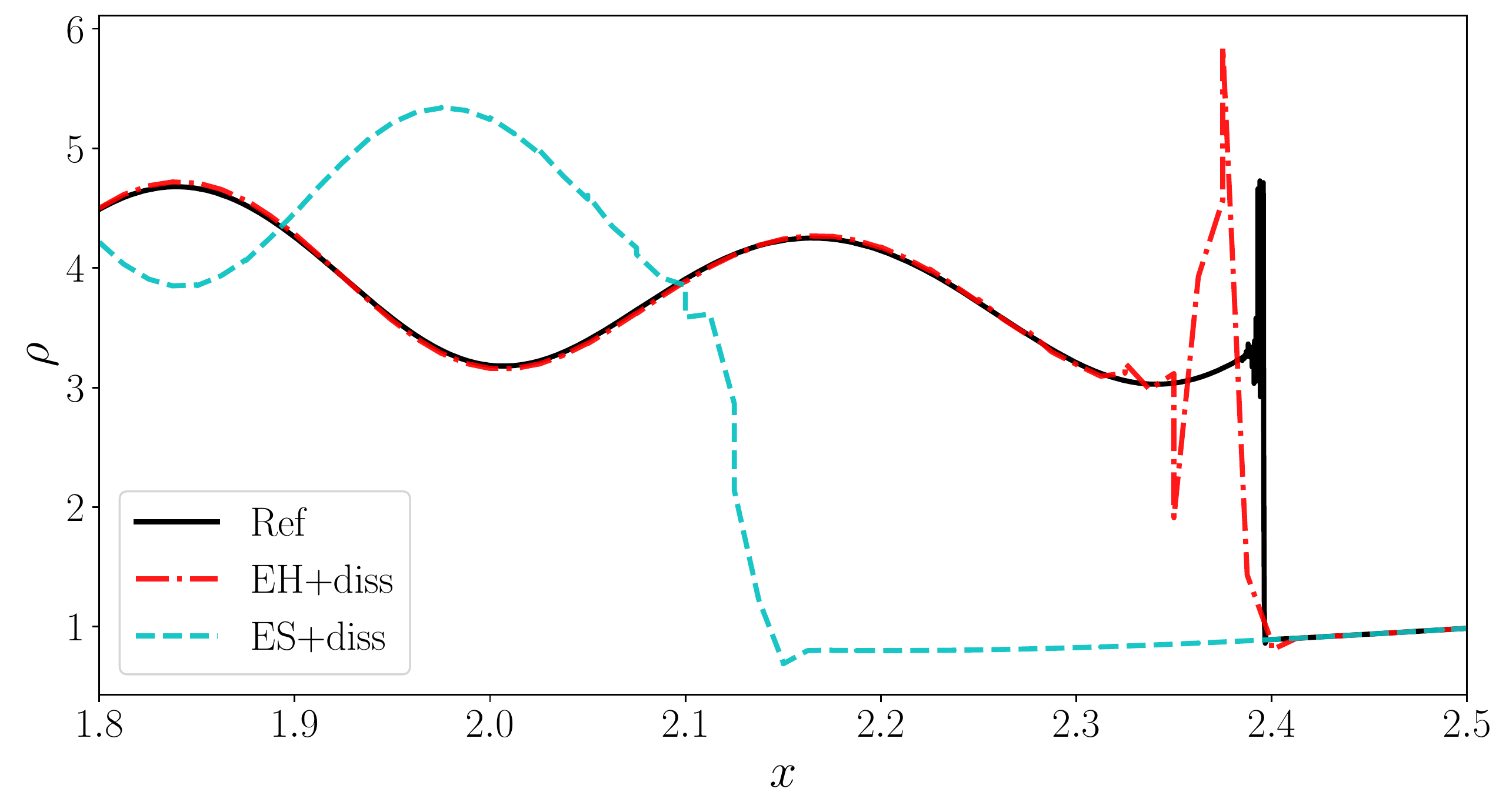}
		\caption{\label{fig:shu_osher_zoomed} $ p=2 $, $ n_e =400 $, closeup}
	\end{subfigure} 
	\caption{\label{fig:shu_osher} Density profile for the Shu-Osher shock tube problem at $ t=1.8 $ using the locally conservative hybrid scheme and the nonconservative entropy-split scheme.}
\end{figure}

\begin{figure}[!t]
	\centering
	\begin{subfigure}{0.5\textwidth}
		\centering
		\includegraphics[scale=0.31]{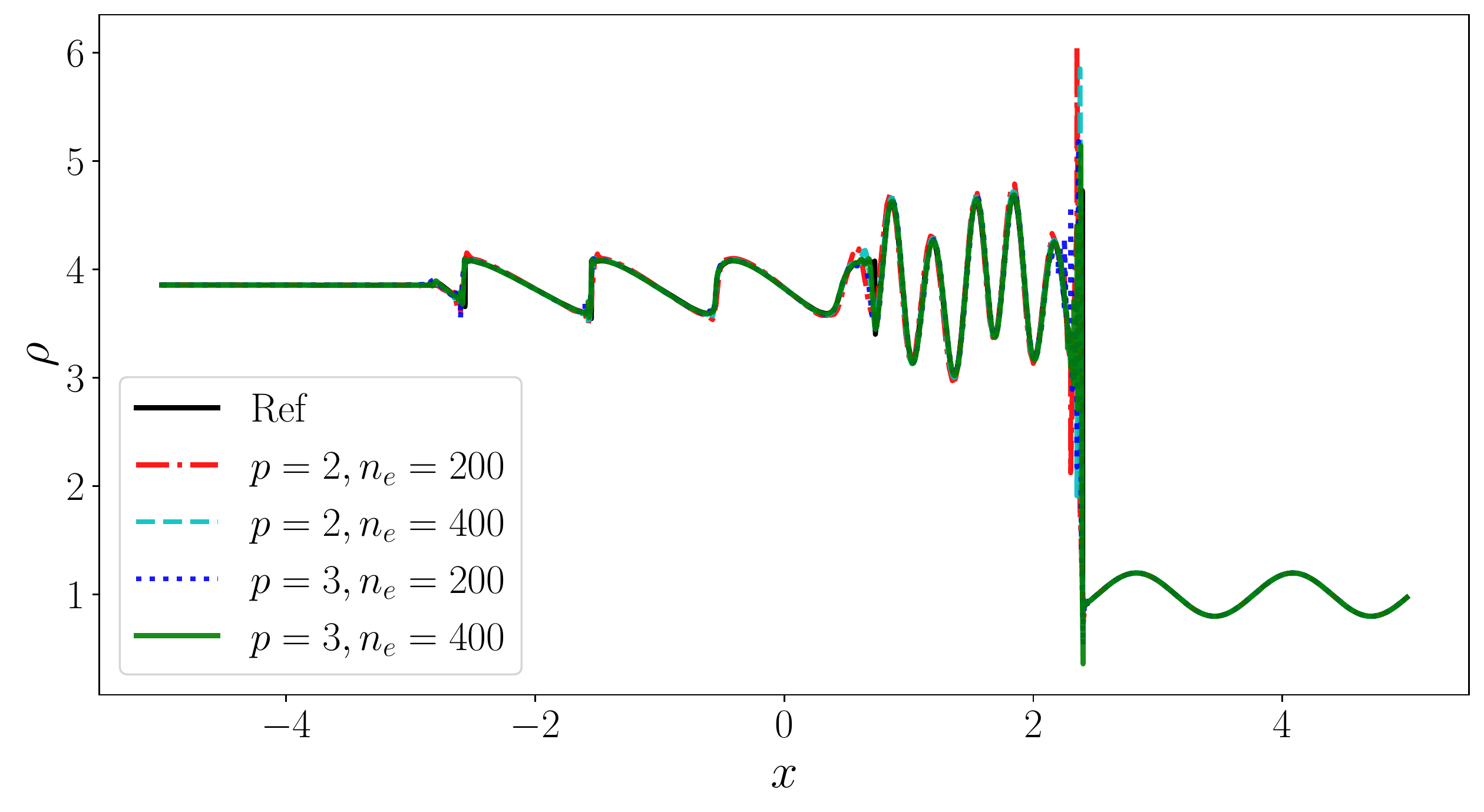}
		\caption{\label{fig:shu_osher_hybrid} hybrid scheme}
	\end{subfigure}\hfill
	\begin{subfigure}{0.5\textwidth}
		\centering
		\includegraphics[scale=0.31]{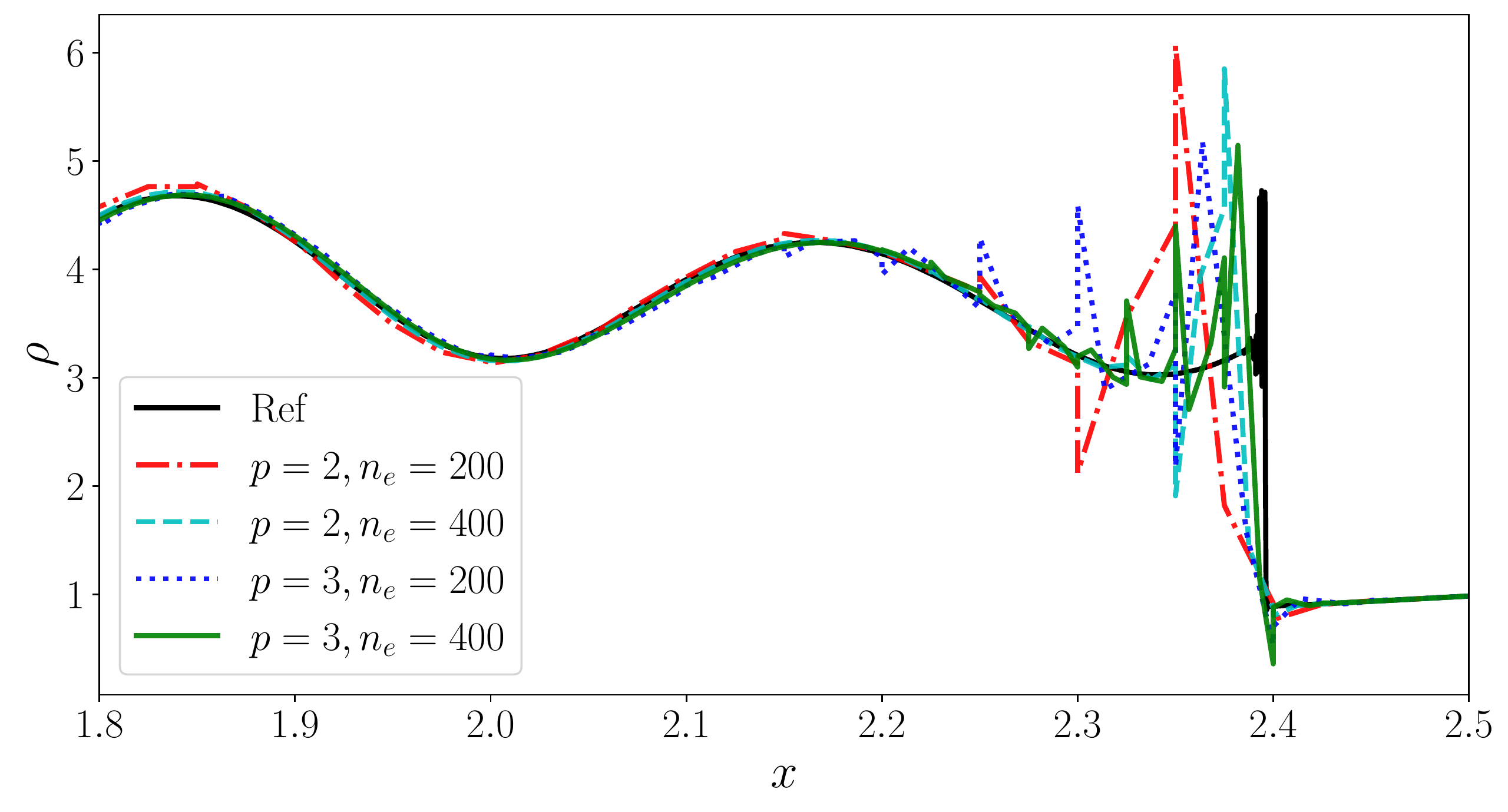}
		\caption{\label{fig:shu_osher_hybrid_zoomed} hybrid scheme, closeup}
	\end{subfigure} 
	\caption{\label{fig:shu_osher_convergence} Density profile for the Shu-Osher shock tube problem at $ t=1.8 $ using the locally conservative hybrid scheme.}
\end{figure}

From the density profile depicted in \cref{fig:shu_osher}, it is clear that the hybrid scheme captures the small and large shock waves accurately. The entropy-split scheme fails to do so as expected. Increasing the mesh resolution or polynomial degree improves the solution accuracy of the hybrid scheme as illustrated in \cref{fig:shu_osher_convergence}. Similar to the Sod shock tube problem, the numerical solution is oscillatory as no limiting or shock capturing method is applied.

\subsection{The isentropic vortex problem}
The isentropic vortex problem is a smooth problem governed by the Euler equations. Both the 2D and 3D versions of the problem are used to study accuracy, primary conservation, entropy conservation, and efficiency of the proposed scheme. The analytical solution for the 2D case on the domain $ \Omega = \left[0,20\right]\times\left[-5,5\right] $ is given as \cite{erlebacher1997interaction, shadpey2020entropy},
\begin{equation}
	\begin{aligned}
		\rho & =\left(1-\frac{\alpha^{2}\left(\gamma-1\right)}{16\gamma\pi}\exp\left(2\left(1-r^{2}\right)\right)\right)^{\frac{1}{\gamma-1}}, & p & =\rho^{\gamma},\\
		u & =1 - \frac{\alpha}{2\pi}\left(\red{x_{2}}-y_{c}\right)\exp\left(1-r^{2}\right), &
		v & =\frac{\alpha}{2\pi}\left(\red{x_{1}}-\left(x_{c}+t\right)\right)\exp\left(1-r^{2}\right),
	\end{aligned}
\end{equation}
where $r^{2}=\left(\red{x_{1}}-x_{c}+t\right)^{2}+\left(\red{x_{2}}-y_{c}\right)^{2}$,
$\alpha=3$ is the vortex strength, $\left(x_c,y_c\right)=\left(5,0\right)$ is the center of the vortex at $ t=0 $,
and $\gamma=7/5$. The grids are curved using the perturbation \cite{chan2019efficient}
\begin{equation} \label{eq:curve 2d vortex}
	\begin{aligned}\red{x_{1}} & =\hat{\red{x}}_{1}+\frac{1}{8}\cos\left(\frac{\pi}{20}\left(\hat{x}_{1}-10\right)\right)\cos\left(\frac{3\pi}{10}\hat{\red{x}}_{2}\right), &
		\red{x_{2}} & =\hat{\red{x}}_2+\frac{1}{8}\sin\left(\frac{\pi}{5}\left(\red{x_{1}}-10\right)\right)\cos\left(\frac{\pi}{10}\hat{\red{x}}_{2}\right),
	\end{aligned}
\end{equation}
where $ (\hat{\cdot}) $ denotes the coordinate of a node before perturbation. In 3D, we use the analytical solution \cite{williams2013nodal,shadpey2020entropy},
\begin{equation}
	\begin{aligned}
		\rho & =\left(1-\frac{2}{25}\left(\gamma-1\right)\exp\left(1-\left(\red{x_{2}}-t\right)^{2}-\red{x}_{1}^{2}\right)\right)^{\frac{1}{\gamma-1}}, & e & =\frac{\rho^{\gamma}}{\gamma(\gamma-1)}+\frac{\rho}{2}\left(\red{V_{1}}^{2}+\red{V_{2}}^{2}+\red{V_{3}}^{2}\right), & \\
		\red{V_{1}} & =-\frac{2}{5}\left(\red{x_{2}}-t\right)\exp\left(\frac{1}{2}\left[1-\left(\red{x_{1}}-t\right)^{2}-\red{x}_{1}^{2}\right]\right), & 
		\red{V_{2}} & =1+\frac{2}{5}\red{x_{1}}\exp\left(\frac{1}{2}\left[1-\left(\red{x_{2}}-t\right)^{2}-\red{x}_{1}^{2}\right]\right), & 
		\red{V_{3}} & =0.
	\end{aligned}
\end{equation}
The domain is $\Omega = \left[-10,10\right]^{3}$, and the grids are curved according
to the perturbation 
\begin{equation}\label{eq:curve 3d vortex}
	{x}_{i}=\hat{{x}}_{i}+0.05\sin\left(\frac{\pi\hat{{x}}_{i}}{2}\right),\quad\forall i\in\{1,2,3\}.
\end{equation} 
\subsubsection{Accuracy}
The high-order accuracy of the entropy-split scheme is validated by running the 2D and 3D isentropic vortex problems on a sequence of curved meshes. In 2D, curved triangular elements are used which are obtained by subdividing rectangular elements into two and applying the perturbation given in \cref{eq:curve 2d vortex}. Similarly, in 3D, curved tetrahedral elements are obtained by subdividing hexahedral elements into six tetrahedra and applying the perturbation given in \cref{eq:curve 3d vortex}. We followed the procedure outlined in \cite{shadpey2020entropy} to compute the $ L^2 $ error, \ie, interpolate the solution from the SBP nodes to a quadrature rule of degree $ 3p + 1 $, integrate the square of the solution error, sum the result over all elements, and take the square root of the sum. Given solution errors on two meshes, the convergence rates are computed as $ \log(\text{error}_2/\text{error}_1 )/ \log(h_2/ h_1)$, where $ h_i = {(\text{vol.}/n_{e,i})} ^{1/d}$ and $ \text{vol.} $ is either the area or volume of the domain.

The solution convergence results for the 2D problem at $ t_f = 2 $ and the 3D case at $ t_f = 1 $ are presented in \cref{tab:accuracy vortex 2d,tab:accuracy vortex 3d}, respectively. CFL values of $ 0.1 $ and $ 0.25 $ are used for the 2D and 3D cases, respectively. For comparison, we have also included the corresponding errors and convergence rates of the Ismail-Roe Hadamard-form discretization. It can be seen from the tables that convergence rates between $ p $ and $ p+1 $ are observed for all cases. In the 2D case, the entropy conservative schemes exhibit an even-odd convergence rate behavior, where the even degree operators converge at rates greater than $ p+0.5 $ while the odd degree operators converge close to a suboptimal rate of $ p $. This behavior has been reported previously, \eg, see \cite{carpenter2015entropy,parsani2016entropy,carpenter2016towards,chan2018discretely}. Adding dissipation improves the solution convergence rates for all cases and decreases the solution errors significantly. The solution errors of the entropy-split and Ismail-Roe Hadamard-form schemes are very close, with the latter having slightly lower errors for almost all cases.

\begin{table*} [!t]
	\small
	\caption{\label{tab:accuracy vortex 2d} Grid convergence study for the 2D isentropic vortex problem on curvilinear grids at $ t_f=2 $.}
	\centering
	\setlength{\tabcolsep}{1.1em}
	\renewcommand*{\arraystretch}{1.2}
	\begin{tabular}{crcccccccc}
		\toprule
		\multirow{3}{*}{$p$}& \multirow{3}{*}{$ n_e $ }   &     \multicolumn{4}{c}{Entropy-split}&    \multicolumn{4}{c}{IR Hadamard-form} \\ 
		\cmidrule(lr){3-6} \cmidrule(lr){7-10}
		& & \multicolumn{2}{c}{conservative} &  \multicolumn{2}{c}{dissipative} &  \multicolumn{2}{c}{conservative} &  \multicolumn{2}{c}{dissipative} \\
		\cmidrule(lr){3-4} \cmidrule(lr){5-6} \cmidrule(lr){7-8} \cmidrule(lr){9-10}
		& & {$ L^2 $ error}  &{rate} & {$ L^2 $ error}  & {rate} &  {$ L^2 $ error}  & {rate} &{$ L^2 $ error}  & {rate} \\
		\midrule
		\multirow{4}{*}{$ 1 $}&    $ 60^2\times 2 $ 		& 2.7515e-01   &     --    & 1.4652e-01   &     --        & 2.8223e-01  &  --         & 1.5185e-01  &    --     \\  
		&    $ 80^2\times 2 $       & 2.1770e-01   & 0.81    & 8.7858e-02   & 1.77         & 2.2491e-01 & 0.78      & 9.0908e-02  & 1.78    \\
		&    $ 100^2\times 2 $      & 1.7857e-01  & 0.88    & 5.8532e-02   & 1.82         & 1.8452e-01 & 0.88     & 6.0450e-02  & 1.82    \\
		&    $ 120^2\times 2 $      & 1.5189e-01 & 0.88    & 4.1900e-02   & 1.83         &1.5668e-01 & 0.89     & 4.3219e-02 & 1.84     \\
		\midrule
		\multirow{4}{*}{$ 2 $}&    $ 60^2\times 2 $ 		& 1.8342e-02  &     --    & 1.0346e-02   &     --        & 1.7201e-02  &  --         & 1.0074e-02  &    --     \\  
		&    $ 80^2\times 2 $       & 8.6064e-03   & 2.63    & 4.7957e-03   & 2.67         & 7.9423e-03  & 2.68     & 4.5671e-03  & 2.74    \\
		&    $ 100^2\times 2 $      & 4.6965e-03  & 2.71    & 2.5913e-03  & 2.75        	& 4.3051e-03 & 2.74     & 2.4300e-03   & 2.82    \\
		&    $ 120^2\times 2 $      & 2.8214e-03   & 2.79     & 1.5550e-03   & 2.80        & 2.5832e-03  & 2.80      & 1.4428e-03  & 2.85     \\
		\midrule
		\multirow{4}{*}{$ 3 $}&    $ 60^2\times 2 $ 		& 1.8747e-03   &     --    & 8.0180e-04   &     --        & 1.9087e-03  &  --         & 7.5947e-04 &    --     \\  
		&    $ 80^2\times 2 $       & 7.3938e-04  & 3.23    & 2.8758e-04   & 3.56       & 7.6724e-04  & 3.16      & 2.7114e-04   & 3.58   \\
		&    $ 100^2\times 2 $      & 3.6119e-04 & 3.21    & 1.3091e-04   & 3.52         & 3.7775e-04   & 3.17     & 1.2265e-04  & 3.55    \\
		&    $ 120^2\times 2 $      & 1.9582e-04  & 3.35    & 6.8669e-05   & 3.53       & 2.0658e-04 & 3.31     & 6.3942e-05  & 3.57    \\
		\midrule
		\multirow{4}{*}{$ 4 $}&    $ 60^2\times 2 $ 		& 1.6301e-04    &     --    & 6.9849e-05  &     --        & 1.4134e-04  &  --         & 6.7207e-05    &    --     \\  
		&    $ 80^2\times 2 $       & 4.4338e-05   & 4.52   & 1.7623e-05   & 4.78         & 3.8150e-05 & 4.55      & 1.6995e-05  & 4.77   \\
		&    $ 100^2\times 2 $      & 1.6038e-05   & 4.55    & 6.0291e-06   & 4.80        & 1.3538e-05 & 4.64      & 5.8206e-06  & 4.80   \\
		&    $ 120^2\times 2 $      & 6.8290e-06  & 4.68    & 2.5207e-06  & 4.78        & 5.7351e-06  & 4.71      & 2.4313e-06   & 4.78    \\
		\bottomrule
	\end{tabular}
\end{table*}

\begin{table*} [!t]
	\small
	\caption{\label{tab:accuracy vortex 3d} Grid convergence study for the 3D isentropic vortex problem on curvilinear grids at $ t_f=1$.}
	\centering
	\setlength{\tabcolsep}{1.1em}
	\renewcommand*{\arraystretch}{1.2}
	\begin{tabular}{crcccccccc}
		\toprule
		\multirow{3}{*}{$p$}& \multirow{3}{*}{$ n_e $ }   &     \multicolumn{4}{c}{Entropy-split}&    \multicolumn{4}{c}{IR Hadamard-form} \\ 
		\cmidrule(lr){3-6} \cmidrule(lr){7-10}
		& & \multicolumn{2}{c}{conservative} &  \multicolumn{2}{c}{dissipative} &  \multicolumn{2}{c}{conservative} &  \multicolumn{2}{c}{dissipative} \\
		\cmidrule(lr){3-4} \cmidrule(lr){5-6} \cmidrule(lr){7-8} \cmidrule(lr){9-10}
		& & {$ L^2 $ error}  &{rate} & {$ L^2 $ error}  & {rate} &  {$ L^2 $ error}  & {rate} &{$ L^2 $ error}  & {rate} \\
		\midrule
		\multirow{4}{*}{$ 1 $}&     $ 15^3\times 6 $  		& 1.1982e+00   &     --    & 1.2204e+00  &     --        & 1.2120e+00  &  --         & 1.2298e+00  &    --     \\  
		&   $ 20^3\times 6 $      & 7.5842e-01   & 1.58    & 7.4754e-01   & 1.70         & 7.5699e-01 & 1.63     & 7.4465e-01  & 1.74    \\
		&    $ 25^3\times 6 $      & 5.2297e-01  & 1.66    & 4.9881e-01   & 1.81         & 5.2184e-01& 1.66      & 4.9923e-01  & 1.79    \\
		&    $ 30^3\times 6 $      & 3.7637e-01 & 1.80    & 3.5070e-01   & 1.93         & 3.7653e-01 & 1.79     & 3.5056e-01, & 1.93     \\
		\midrule
		\multirow{4}{*}{$ 2 $}&     $ 15^3\times 6 $  		& 3.7061e-01   &     --    & 3.3001e-01   &     --        & 3.6917e-01  &  --         & 3.2805e-01  &    --     \\  
		&    $ 20^3\times 6 $       & 1.7847e-01   & 2.54    & 1.5410e-01   & 2.64         & 1.7721e-01  & 2.55     & 1.5354e-01  & 2.63    \\
		&    $ 25^3\times 6 $      & 9.4926e-02  & 2.82    & 8.1664e-02   & 2.84         & 9.3071e-02 & 2.88     & 8.0866e-02   & 2.87    \\
		&    $ 30^3\times 6 $     & 5.9186e-02   &  2.59     & 4.9042e-02   & 2.79        & 5.7696e-02  & 2.62      & 4.8443e-02  & 2.81     \\
		\midrule
		\multirow{4}{*}{$ 3 $}&     $ 15^3\times 6 $  		& 1.2425e-01   &     --    & 9.5332e-02   &     --        & 1.1558e-01 &  --         & 9.3329e-02 &    --     \\  
		&    $ 20^3\times 6 $       & 4.6054e-02  & 3.45    & 3.7759e-02   & 3.21        & 4.4596e-02  & 3.31      & 3.6779e-02  & 3.23    \\
		&    $ 25^3\times 6 $      & 2.2157e-02 & 3.27    & 1.6972e-02   & 3.58        & 2.1349e-02 & 3.30      & 1.6885e-02  & 3.48    \\
		&    $ 30^3\times 6 $      & 1.1566e-02  & 3.56    & 8.5225e-03   & 3.77        & 1.0965e-02  & 3.65     & 8.3907e-03  & 3.83     \\
		\midrule
		\multirow{4}{*}{$ 4 $}&     $ 15^3\times 6 $  		&3.8613e-02    &     --    & 2.7985e-02        &     --        & 3.5419e-02  &  --         & 2.7259e-02    &    --     \\  
		&    $ 20^3\times 6 $       & 1.2417e-02   & 3.94    & 8.4335e-03   & 4.16        & 1.1120e-02 & 4.02      & 7.9801e-03  & 4.27    \\
		&    $ 25^3\times 6 $      & 4.7588e-03   & 4.29    & 2.9824e-03   & 4.65        &  4.1733e-03 & 4.39      & 2.7874e-03  & 4.71   \\
		&    $ 30^3\times 6 $      & 2.0922e-03   & 4.50    & 1.2713e-03    &  4.67        &1.8952e-03  & 4.32     & 1.2137e-03            & 4.56         \\
		\bottomrule
	\end{tabular}
\end{table*}

\subsubsection{Primary conservation}
We investigate the convergence of the primary conservation error with mesh refinement for the entropy-split scheme, where the conservation error is computed as
\begin{equation*}
	\text{error} =\frac{ \abs{\left[\sum_{\Omega_{k}\in\fn T_{h}}\overline{\I}^{T}\overline{\H}_{k}\bm{u}_{k}\right]_{t=t_f} - \left[\sum_{\Omega_{k}\in\fn T_{h}}\overline{\I}^{T}\overline{\H}_{k}\bm{u}_{k}\right]_{t=0}}}{ \abs{ \left[\sum_{\Omega_{k}\in\fn T_{h}}\overline{\I}^{T}\overline{\H}_{k}\bm{u}_{k}\right]_{t=0}}}.
\end{equation*}
The entropy-split scheme is not conservative, but the conservation errors in mass, momentum, and energy decrease quickly with mesh refinement for smooth solutions. The convergence of mass and energy conservation errors for the 2D and 3D entropy-split discretizations of the isentropic vortex problems is reported in \cref{tab:conservation}. We observe superconvergence rates close to or greater than $ 2p $ for both quantities in almost all cases, which is in accordance with \cref{rem:conservation}. This is an encouraging property of the entropy-split scheme as it allows the loss of conservation to be mitigated with modest mesh refinement compared to what would have been needed if the convergence rate was limited to $ p+1 $. Although not reported, we have observed similar superconvergence behavior for the conservation of momentum in each direction.

\begin{table*} [!t]
	\small
	\caption{\label{tab:conservation} Mass and energy conservation grid convergence study for the isentropic vortex problems.}
	\centering
	\setlength{\tabcolsep}{0.85em}
	\renewcommand*{\arraystretch}{1.2}
	\begin{tabular}{crccccrcccc}
		\toprule
		\multirow{3}{*}{$p$}&   \multicolumn{5}{c}{2D ES scheme}&    \multicolumn{5}{c}{3D ES scheme} \\ 
		\cmidrule(lr){2-6} \cmidrule(lr){7-11}
		& \multirow{2}{*}{$ n_e $ }& \multicolumn{2}{c}{mass} &  \multicolumn{2}{c}{energy} &  \multirow{2}{*}{$ n_e $ }& \multicolumn{2}{c}{mass} &  \multicolumn{2}{c}{energy} \\
		\cmidrule(lr){3-4} \cmidrule(lr){5-6} \cmidrule(lr){8-9} \cmidrule(lr){10-11}
		& & {error}  &{rate} & {error}  & {rate} &  & {error}  & {rate} &{error}  & {rate} \\
		\midrule
		\multirow{4}{*}{$ 1 $}&    $ 60^2\times 2 $ 		& 6.0705e-06   &     --    & 1.1264e-05   &     --        &    $ 15^3\times 6 $ & 3.8546e-06  &  --      & 8.3150e-06  &    --     \\  
		&    $ 80^2\times 2 $       & 1.7304e-06   & 4.36    & 1.5949e-06   & 6.80       &    $ 20^3\times 6 $& 1.5919e-06 & 3.07      & 3.0643e-06  & 3.47    \\
		&    $ 100^2\times 2 $      & 2.7436e-07  & 8.25    & 1.3462e-06   & 0.76        &    $ 25^3\times 6 $& 1.0567e-06 & 1.84     & 2.0954e-06  & 1.70    \\
		&    $ 120^2\times 2 $      & 1.6641e-07  & 2.74     & 1.9304e-06   & $-1.98^{\dagger}$      &    $ 30^3\times 6 $& 7.1835e-07 & 2.12     & 1.3534e-06 & 2.40    \\
		\midrule
		\multirow{4}{*}{$ 2 $}&    $ 60^2\times 2 $ 		& 1.3948e-07   &     --    & 3.1346e-07    &     --        &    $ 15^3\times 6 $ & 1.1523e-06  &  --         & 2.1950e-06  &    --     \\  
		&    $ 80^2\times 2 $       & 2.9296e-08   & 5.42    & 7.2323e-08   & 5.10         &   $ 20^3\times 6 $& 4.8005e-07  & 3.04     & 9.1724e-07 & 3.03    \\
		&    $ 100^2\times 2 $      & 9.1568e-09   & 5.21    & 2.4679e-08  & 4.82         &    $ 25^3\times 6 $& 2.0123e-07 & 3.90    & 3.8547e-07   & 3.89    \\
		&    $ 120^2\times 2 $      &3.4885e-09   & 5.29   & 1.0153e-08   & 4.87         &    $ 30^3\times 6 $& 8.7548e-08  & 4.56      & 1.6570e-07 & 4.63    \\
		\midrule
		\multirow{4}{*}{$ 3 $}&    $ 60^2\times 2 $ 		& 2.8481e-09   &     --    & 6.0358e-09   &     --        &    $ 15^3\times 6 $  & 2.9194e-07  &  --         & 4.8042e-07&    --     \\  
		&    $ 80^2\times 2 $       & 4.9214e-10  & 6.10    & 1.0529e-09     & 6.07        &    $ 20^3\times 6 $ & 6.4213e-08  & 5.26      & 4.8042e-07  & 4.93    \\
		&    $ 100^2\times 2 $      & 1.3194e-10  & 5.90    & 2.8771e-10   & 5.81       &    $ 25^3\times 6 $ & 1.5399e-08 & 6.40      & 4.8042e-07  & 6.76    \\
		&    $ 120^2\times 2 $      & 4.1898e-11  & 6.29    & 9.3279e-11    & 6.18        &    $ 30^3\times 6 $ & 5.1313e-09  & 6.03    & 4.8042e-07  & 5.96     \\
		\midrule
		\multirow{4}{*}{$ 4 $}&    $ 60^2\times 2 $ 		& 3.3298e-11    &     --    & 6.9952e-11   &     --       &    $ 15^3\times 6 $  & 5.8265e-08  &  --         & 4.8042e-07    &    --     \\  
		&    $ 80^2\times 2 $       & 2.6775e-12   & 8.76   & 5.5775e-12   & 8.79        &    $ 20^3\times 6 $& 6.9086e-09 & 7.41      & 4.8042e-07  & 7.49    \\
		&    $ 100^2\times 2 $      & 3.8784e-13  & 8.66    & 7.9809e-13    & 8.71        &    $ 25^3\times 6 $ & 1.2298e-09 &  7.73      & 2.0520e-09  & 7.81   \\
		&    $ 120^2\times 2 $      & 8.4176e-14   & 8.38    & 1.5320e-13   & 9.05       &    $ 30^3\times 6 $& 2.4543e-10 &  8.84      & 4.0164e-10   & 8.95     \\
		\bottomrule
	\end{tabular}
	\begin{tablenotes}
		\item[] $^{\dagger}$Upon further mesh refinement, the convergence rate of the error in energy conservation is found to be approximately $ 2$.
	\end{tablenotes}
\end{table*}
\subsubsection{Long-time entropy conservation and error behavior}
\begin{figure}[!t]
	\centering
	\begin{subfigure}{0.48\textwidth}
		\centering
		\includegraphics[scale=0.30]{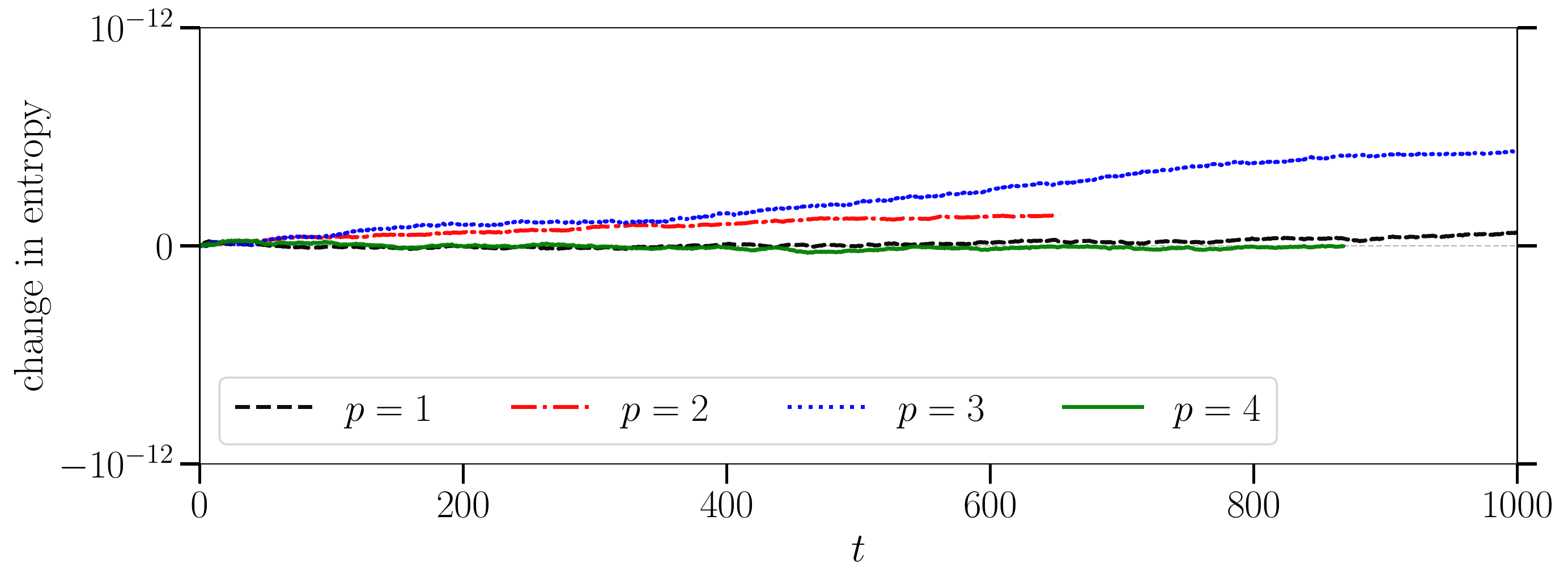}
		\caption{\label{fig:long_entropy} $\beta = 2.5$ }
	\end{subfigure}\hfill
	\begin{subfigure}{0.48\textwidth}
		\centering
		\includegraphics[scale=0.30]{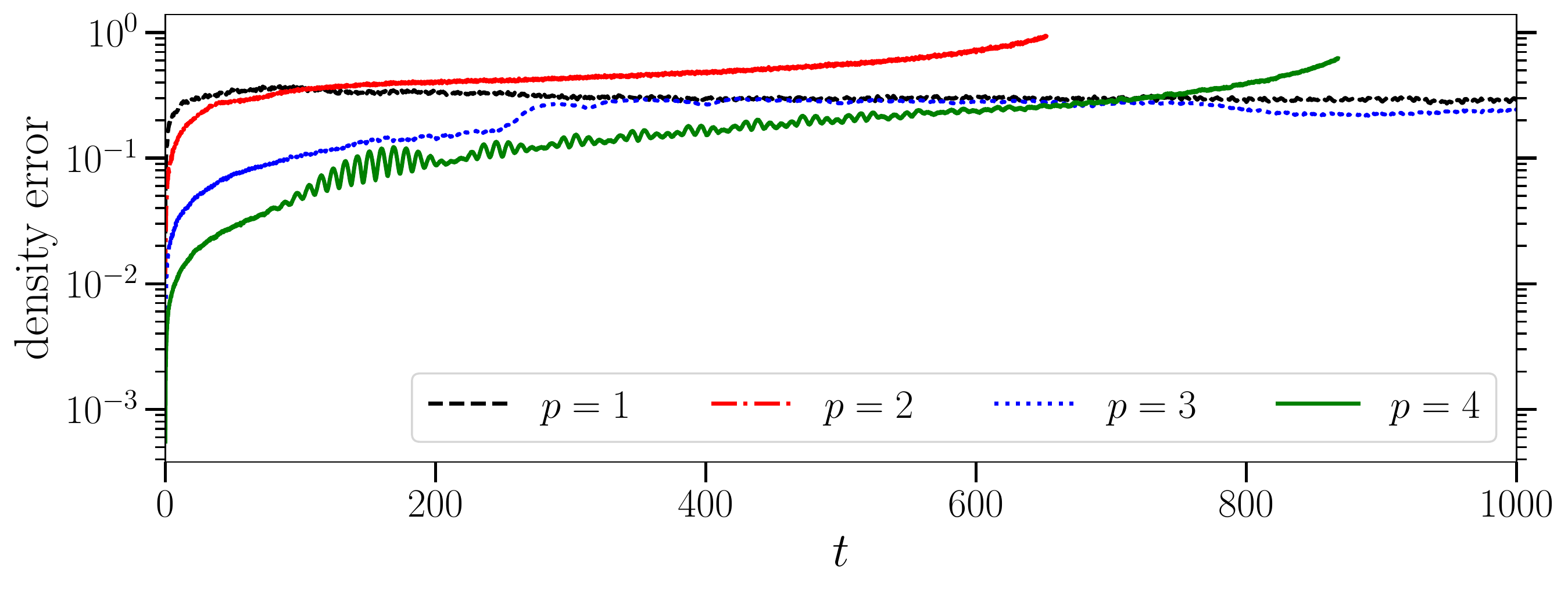}
		\caption{\label{fig:long_error}$\beta = 2.5 $}
	\end{subfigure} 
	\\
	\begin{subfigure}{0.48\textwidth}
		\centering
		\includegraphics[scale=0.30]{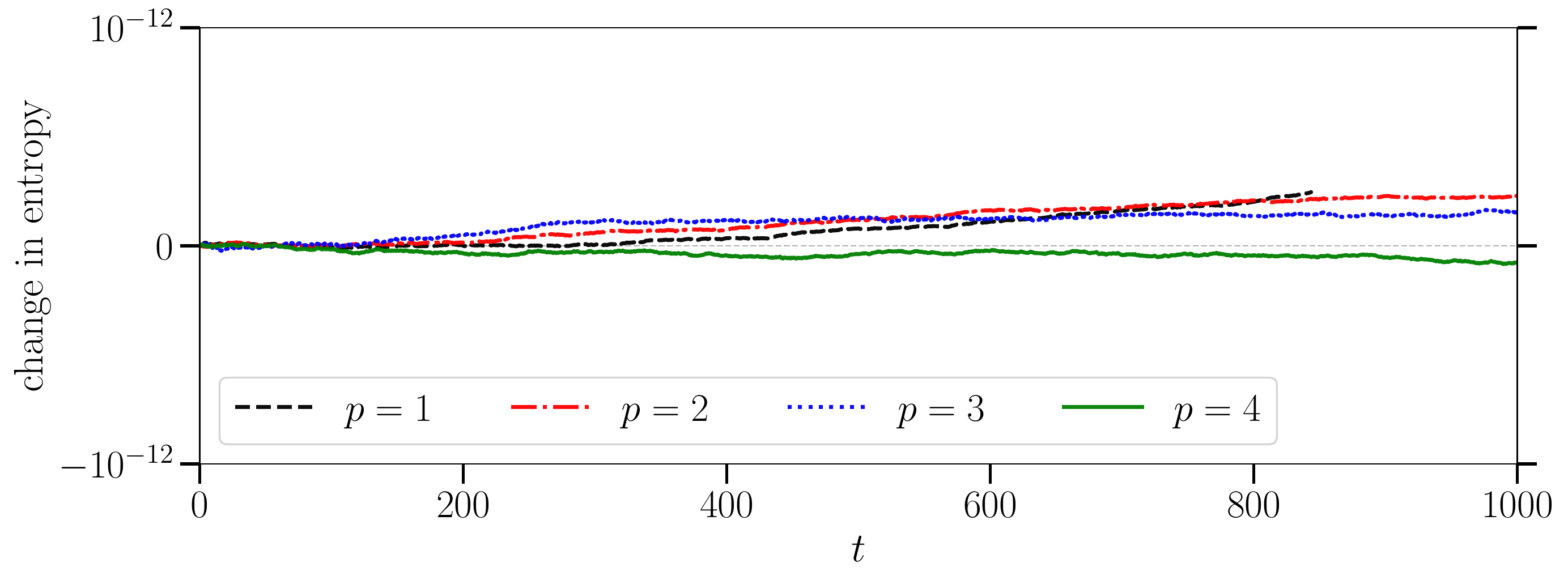}
		\caption{\label{fig:long_entropy_beta2} $\beta = 2 $ }
	\end{subfigure}\hfill
	\begin{subfigure}{0.48\textwidth}
		\centering
		\includegraphics[scale=0.30]{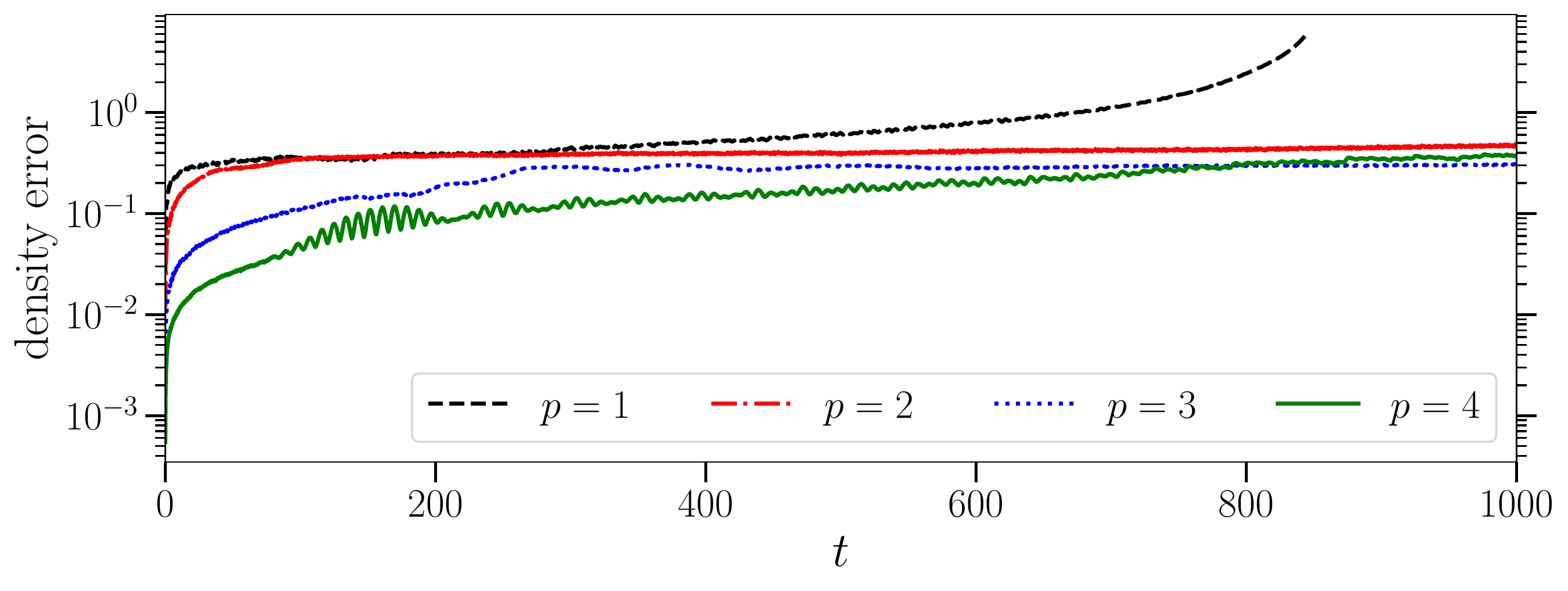}
		\caption{\label{fig:long_error_beta2} $\beta = 2 $}
	\end{subfigure} 
	\caption{\label{fig:long_time_error} Long time entropy conservation and error behavior of the entropy-split method applied for the two dimensional isentropic vortex problem using $ n_e = 800 $ curved elements and $ \text{CFL} = 0.5 $.}
\end{figure}

To study the long-time entropy conservation and error behavior of the entropy-split discretization, we run the 2D isentropic vortex problem until $ t_f = 1000 $, which corresponds to $ 50 $ revolutions of the vortex. For this case, the CFL is set to $ 0.5 $ and a coarse mesh with a total of $ n_e = 800 $ curved elements is used. The results are shown in \cref{fig:long_time_error}. The entropy is conserved for all cases until the end of the simulation or termination due to instabilities. For $ \beta = 2.5 $, the odd degree cases run until the final time, while the even degree cases crash before the end of the simulation. Reducing the CFL value by an order of magnitude did not resolve the issue. The cause of instability is not clear, but it seems to be associated with the value of $ \beta $. \cref{fig:long_entropy_beta2,fig:long_error_beta2} show that for $ \beta = 2 $, all cases run until the end of the simulation except for the $ p=1 $ case, which suggests that further research is required to choose the values of $ \beta $ systematically. When \blue{interface} dissipation is added, all cases run until the end of the simulation for both $ \beta=2.5 $ and $ \beta=2 $ cases. We have also tested the long-time entropy conservation and error behavior of the Ismail-Roe Hadamard-form discretization. We found that the scheme conserves entropy until the final time, and no instability is encountered for all cases. Furthermore, we note that the error values are comparable to those obtained with the entropy-split scheme.

\subsubsection{Efficiency}
One of the appealing aspects of the entropy-split scheme is its efficiency. As described in \cref{subsec:efficiency considerations}, for implementations with multidimensional SBP operators, the cost of computing the entropy-split discretization of the inviscid flux derivatives is lower than the corresponding cost of the Hadamard-form scheme. The efficiency depends on factors such as the type of operators used, the use of tensor-product structure, the SATs, code optimization, and many other implementation details. We have attempted to make fair comparisons between the methods by running both the entropy-split and Hadamard form implementations within a single code base; only the computation of the inviscid flux derivatives is handled with different subroutines. A number of optimization techniques recommended  in \cite{ranocha2021efficient} for the Hadamard-form scheme are applied (see, \cref{subsec:efficiency considerations}), and similar level of optimization is applied for the entropy-split scheme.  

\begin{figure}[!t]
	\centering
	\includegraphics[scale=0.47]{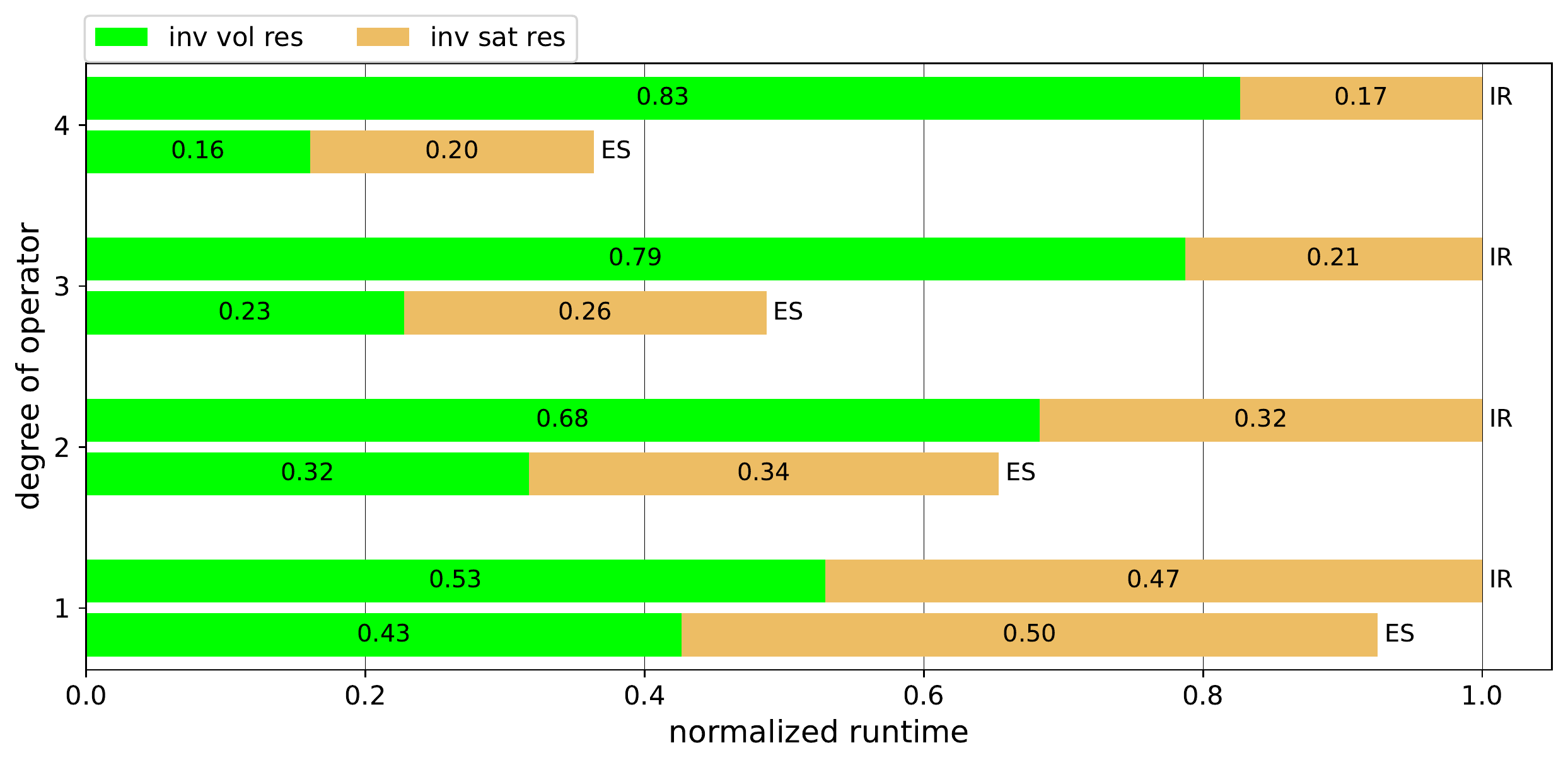}
	\caption{\label{fig:efficiency vortex} Efficiency comparisons of the entropy-split scheme and the Ismail-Roe Hadamard-form discretization of the 3D isentropic vortex problem.}
\end{figure}

The 3D isentropic vortex problem discretized with multidimensional diagonal-$ \E $ operators is used for the efficiency comparison. Only the costs of computing the spatial residuals are compared, and other implementation specific computational costs such as operator construction and time marching are not included. Furthermore, only serial performance is reported. The problem is run for $ 50 $ timesteps with $ n_e = 750 $, and the average performance of each scheme over five runs normalized by the average runtime of the Ismail-Roe Hadamard-form scheme is illustrated in \cref{fig:efficiency vortex}. It can be seen from the figure that as the degree of operators increases the computational cost of the Hadamard-form implementation increases quickly, while the entropy-split scheme scales well for implementations with high-order operators. For the degree three and four cases, the computational cost of the spatial residual terms with the entropy-split scheme is about $ 2 $ and $ 3 $ times less expensive, respectively, compared to the Ismail-Roe Hadamard-form implementation. 

\subsection{Manufactured solution for the \red{compressible} Navier-Stokes equations}
The method of manufactured solutions is used to test the accuracy of the entropy-split discretization of the \red{compressible} Navier-Stokes equations. To this end, we use the manufactured solution \cite{schnucke2020entropy}
\begin{alignat*}{4}
	\rho &=2.0+ 0.1\sin(\pi \left(\Sigma_{i=1}^{d}\bm{x}_i - 0.6t\right)), & \quad  V_i&= 1, &  \quad p& = (\gamma - 1)(\rho ^2 -0.5 \rho \bm{V}^T\bm{V}), \quad i=\{1,\dots, d\},
\end{alignat*}
defined on the periodic domain $ \Omega = [-1,1]^d $. The domain is discretized with SBP diagonal-$ \E $ elements which are curved according to the mapping \cite{chan2019skew}
\begin{align}
	\bm{x}_i &= \hat{\bm{x}}_i + \frac{1}{8}\cos\left(\frac{\pi}{2}\hat{\bm{x}}_i\right)\prod_{j=1,j\neq i}^{d}\sin\left(\pi \hat{\bm{x}}_j\right), \quad i=\{1,\dots, d\}.
\end{align}
The source terms are obtained by substituting the manufactured solution into the \red{compressible} Navier-Stokes equations, \cref{eq:nse}, and solving for the difference of the left and right hand sides. A constant viscosity, $ \mu = 0.01 $, is used, the gas constant and Prandtl number are set to $ R=1 $ and $ Pr=0.71 $, respectively.  The BO SAT coefficients given in \cref{eq:bo sats} are used to couple the viscous terms at element interfaces. No \blue{interface} dissipation is added with either the inviscid or viscous SATs.   

Both the 2D and 3D cases are run for $ 1000 $ timesteps until the final time, $ t_f = 0.001 $. \cref{tab:nse accuracty} shows the $ L^2 $ solution error and convergence rates obtained with the entropy-split and Ismail-Roe Hadamard-form discretizations. We observe convergence rates of approximately $ p+1 $ for all cases. Furthermore, the error values obtained with the two schemes are comparable.

\begin{table*} [!t]
	\small
	\caption{\label{tab:nse accuracty} Grid convergence study for the manufactured solution of the \red{compressible} Navier-Stokes equations at $ t_f = 0.001 $.}
	\centering
	\setlength{\tabcolsep}{0.90em}
	\renewcommand*{\arraystretch}{1.2}
	\begin{tabular}{crccccrcccc}
		\toprule
		\multirow{3}{*}{$p$}&   \multicolumn{5}{c}{2D }&    \multicolumn{5}{c}{3D} \\ 
		\cmidrule(lr){2-6} \cmidrule(lr){7-11}
		& \multirow{2}{*}{$ n_e $ }& \multicolumn{2}{c}{ES scheme} &  \multicolumn{2}{c}{IR scheme} &  \multirow{2}{*}{$ n_e $ }& \multicolumn{2}{c}{ES scheme} &  \multicolumn{2}{c}{IR scheme} \\
		\cmidrule(lr){3-4} \cmidrule(lr){5-6} \cmidrule(lr){8-9} \cmidrule(lr){10-11}
		& & {$ L^2 $ error}  &{rate} & {$ L^2 $ error}  & {rate} &  & {$ L^2 $ error}  & {rate} &{$ L^2 $ error}  & {rate} \\
		\midrule
		\multirow{4}{*}{$ 1 $}&    $ 20^2\times 2 $ 		& 2.8996e-03  &     --    & 2.8992e-03   &     --        &    $ 6^3\times 6 $ & 1.5492e-01  &  --      & 1.5492e-01  &    --     \\  
		&    $ 30^2\times 2 $       &  1.2991e-03   & 1.98    & 1.2988e-03   & 1.98      &    $ 9^3\times 6 $& 7.2321e-02  & 1.87      & 7.2323e-02  & 1.87    \\
		&    $ 40^2\times 2 $      & 7.3772e-04   & 1.96    & 7.3744e-04   & 1.96       &    $ 12^3\times 6 $& 4.1432e-02 & 1.93     & 4.1434e-02  & 1.93    \\ 
		&    $ 50^2\times 2 $      & 4.7759e-04   & 1.94    & 4.7732e-04   & 1.94       &    $ 15^3\times 6 $& 2.6749e-02 & 1.96     & 2.6751e-02 & 1.96    \\
		\midrule
		\multirow{4}{*}{$ 2 $}&    $ 20^2\times 2 $ 		& 9.7812e-05  &     --    & 9.7827e-05    &     --        &    $ 6^3\times 6 $ & 3.4847e-02  &  --         & 3.4845e-02  &    --     \\  
		&    $ 30^2\times 2 $       & 2.8970e-05   & 3.00    &2.8966e-05   &  3.00         &   $ 9^3\times 6 $& 1.1110e-02  & 2.81     & 1.1110e-02  & 2.81    \\
		&    $ 40^2\times 2 $      & 1.2199e-05   & 3.00    & 1.2192e-05   &  3.00        &    $ 12^3\times 6 $& 4.8167e-03 & 2.90    & 4.8179e-03   & 2.90   \\
		&    $ 50^2\times 2 $      & 6.2329e-06   &  3.00   & 6.2266e-06   &  3.00         &    $ 15^3\times 6 $& 2.4989e-03  & 2.94   & 2.5002e-03 & 2.93    \\
		\midrule
		\multirow{4}{*}{$ 3 $}&    $ 20^2\times 2 $ 		& 2.4811e-06   &     --    & 2.4811e-06   &     --        &    $ 6^3\times 6 $  & 7.0398e-03  &  --         & 7.0223e-03 &    --     \\  
		&    $ 30^2\times 2 $       &4.9958e-07  & 3.95   & 4.9848e-07     & 3.96        &    $ 9^3\times 6 $ & 1.5260e-03  & 3.77      & 1.5209e-03  & 3.77    \\
		&    $ 40^2\times 2 $      & 1.6181e-07  & 3.91     & 1.6099e-07   & 3.92        &    $ 12^3\times 6 $ & 5.0177e-04 & 3.86      & 4.9936e-04  & 3.87    \\
		&    $ 50^2\times 2 $      & 6.8144e-08   & 3.87    & 6.7564e-08    & 3.89        &    $ 15^3\times 6 $ & 2.0991e-04  & 3.90    & 2.0852e-04  & 3.91    \\
		\midrule
		\multirow{4}{*}{$ 4 $}&    $ 20^2\times 2 $ 		& 5.8660e-08    &     --    & 5.7522e-08   &     --       &    $ 6^3\times 6 $  & 1.3981e-03 &  --         & 1.3417e-03    &    --     \\  
		&    $ 30^2\times 2 $       & 8.0093e-09   & 4.91   & 7.6089e-09  & 4.98        &    $ 9^3\times 6 $& 2.1934e-04 & 4.56      & 2.0216e-04 & 4.66    \\
		&    $ 40^2\times 2 $      & 1.9966e-09   & 4.82    & 1.8194e-09  & 4.97        &    $ 12^3\times 6 $ & 5.6838e-05 & 4.69      & 5.0252e-05  & 4.83   \\
		&    $ 50^2\times 2 $      & 6.9568e-10    & 4.72    & 6.0473e-10   & 4.93      &    $ 15^3\times 6 $& 1.9956e-05 & 4.69      & 1.6842e-05   & 4.89     \\
		\bottomrule
	\end{tabular}
\end{table*}

\subsection{Viscous Taylor-Green vortex}
The Taylor-Green vortex problem is a challenging test case frequently used to study the accuracy and robustness of high-order numerical methods. For many schemes, the simulation of stretching large scale eddies breaking down into smaller eddies which are eventually dissipated due to viscosity requires a careful application of artificial dissipation or other ad-hoc stabilization techniques. Maintaining a balance between the amount of dissipation required for stability and the amount that can be added without compromising accuracy is challenging. Entropy stable high-order methods offer enhanced stability properties that eliminate the need to use ad-hoc stabilization methods. Our goal is to investigate whether the entropy-split scheme, despite not being provably entropy stable for the \red{compressible} Navier-Stokes equations, remains practically stable for the under-resolved Taylor-Green vortex problem, and whether it provides an accurate solution for a well-resolved case. We also present how the efficiency of the entropy-split scheme for viscous flows compares with that of the Hadamard-form discretization.

The viscous Taylor-Green vortex problem is defined on the domain $ \Omega = \left[-\pi,\pi\right]^{3} $ with the initial conditions \cite{debonis2013solutions}
\begin{alignat*}{3}
	\red{V_{1}}&=u_{0}\sin\left(\frac{\red{x_{1}}}{L}\right)\cos\left(\frac{\red{x_{2}}}{L}\right)\cos\left(\frac{\red{x_{3}}}{L}\right),&\quad p&=p_{0}+\frac{\rho_{0}u_{0}^{2}}{16}\left[\cos\left(\frac{2\red{x_{1}}}{L}\right)+\cos\left(\frac{2\red{x_{2}}}{L}\right)\right]\left[\cos\left(\frac{2\red{x_{3}}}{L}\right)+2\right],&\quad \rho&=\frac{p}{RT_{0}},\\
	\red{V_{2}}&=-u_{0}\cos\left(\frac{\red{x_{1}}}{L}\right)\sin\left(\frac{\red{x_{2}}}{L}\right)\cos\left(\frac{\red{x_{3}}}{L}\right),&\quad \red{V_{3}}&=0,&&
\end{alignat*}
where $ u_{0} =  \rho_{0}= R = L = 1$, $ T_{0} = p_{0}/(R\rho_{0}) $, $ M=u_{0}/\sqrt{\gamma R T_{0}}=0.1 $ is the Mach number, and $ p_{0} = \rho_{0} u_{0}^{2}/(\gamma M^2)$. The Reynolds number is $ Re=\rho_{0} u_{0}L/\mu_{0}=\{100, 200\} $, where $ \mu_{0} $ is a constant viscosity, and the Prandtl number is set to $ Pr=0.71 $. The scaled kinetic energy is computed as
\begin{equation}
	E_{k}=\frac{1}{\rho_{0}\left|\Omega\right|}\int_{\Omega}\rho\frac{\bm{V}\cdot \bm{V}}{2}\dd{\Omega},
\end{equation}
where $ \left|\Omega\right|=8\pi^{3} $, and the kinetic energy dissipation rate is given by 
\begin{equation}
	\varepsilon=-\der[E_{k}]t,
\end{equation}
which is computed using the fifth-order classical SBP derivative operator. 

\begin{figure}[!t]
	\centering
	\includegraphics[scale=0.45]{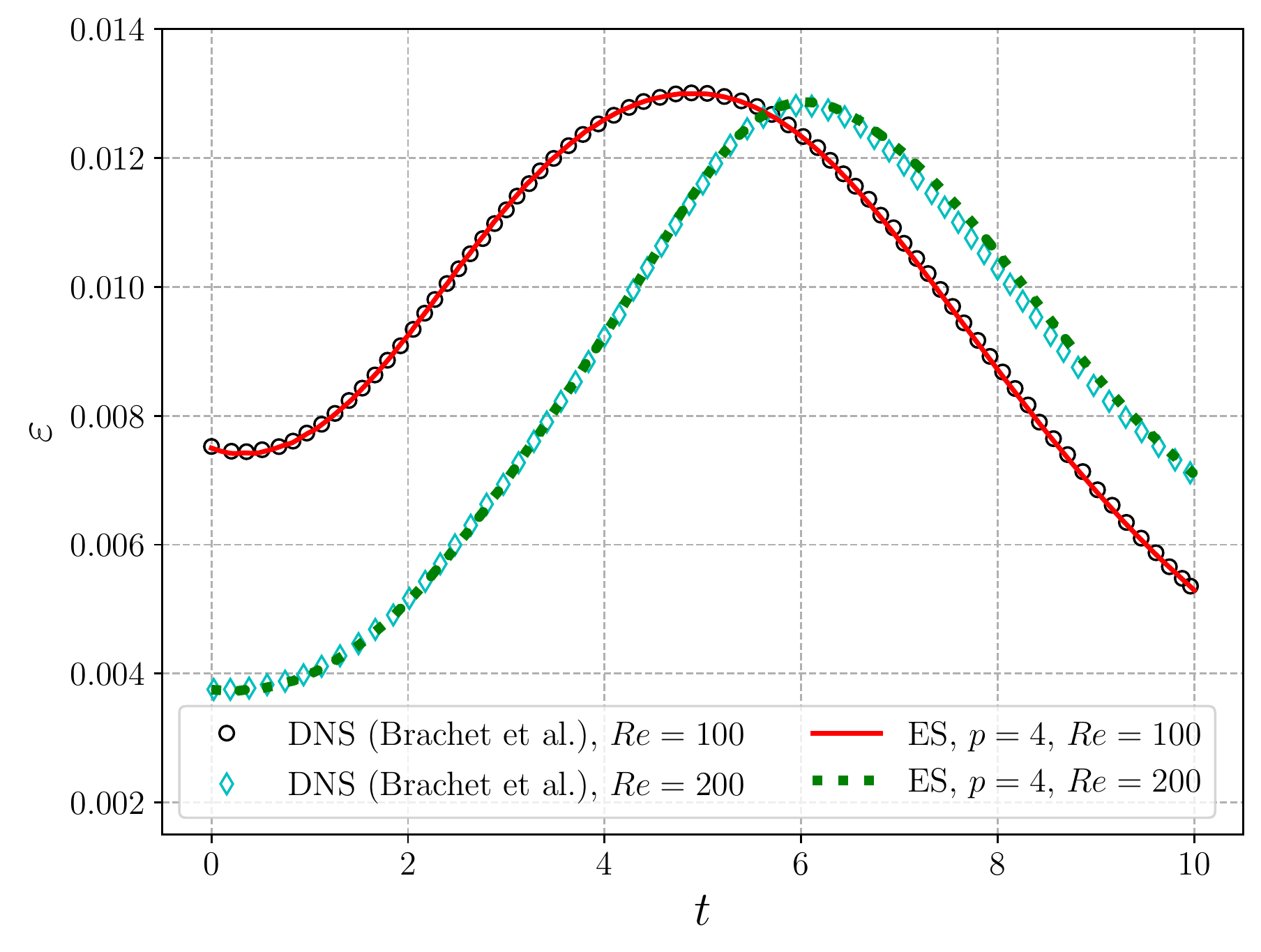}
	\caption{\label{fig:accuracy tgv} Kinetic energy dissipation rate of the viscous Taylor-Green vortex problem solved with the entropy-split scheme using approximately $ 64^3 $ and $ 128^3 $ DOFs for the $ Re=100 $ and $ Re=200 $ cases, respectively.}
\end{figure}

The problem is solved on tetrahedral meshes with approximately $ 64^3 $ and $ 128^3 $ degrees of freedom (DOF) for $ Re=100 $ and $ Re=200 $ cases, respectively. The standard RK4 method is used for time marching. Despite not being provably entropy stable for the \red{compressible} Navier-Stokes equations, the entropy-split scheme remains stable when applied to this challenging test case, including at larger Reynolds numbers, \eg, $ Re=1600 $. \cref{fig:accuracy tgv} shows that the kinetic energy dissipation rates obtained using the entropy-split scheme match the DNS results of Brachet \etal \cite{brachet1983small} reasonably well. For larger values of Reynolds number, the DOFs required to resolve the flow accurately using the SBP diagonal-$ \E $ tetrahedral elements increases quickly. A further investigation on the accuracy of the operators is left for future research.

\begin{figure}[!t]
	\centering
	\includegraphics[scale=0.47]{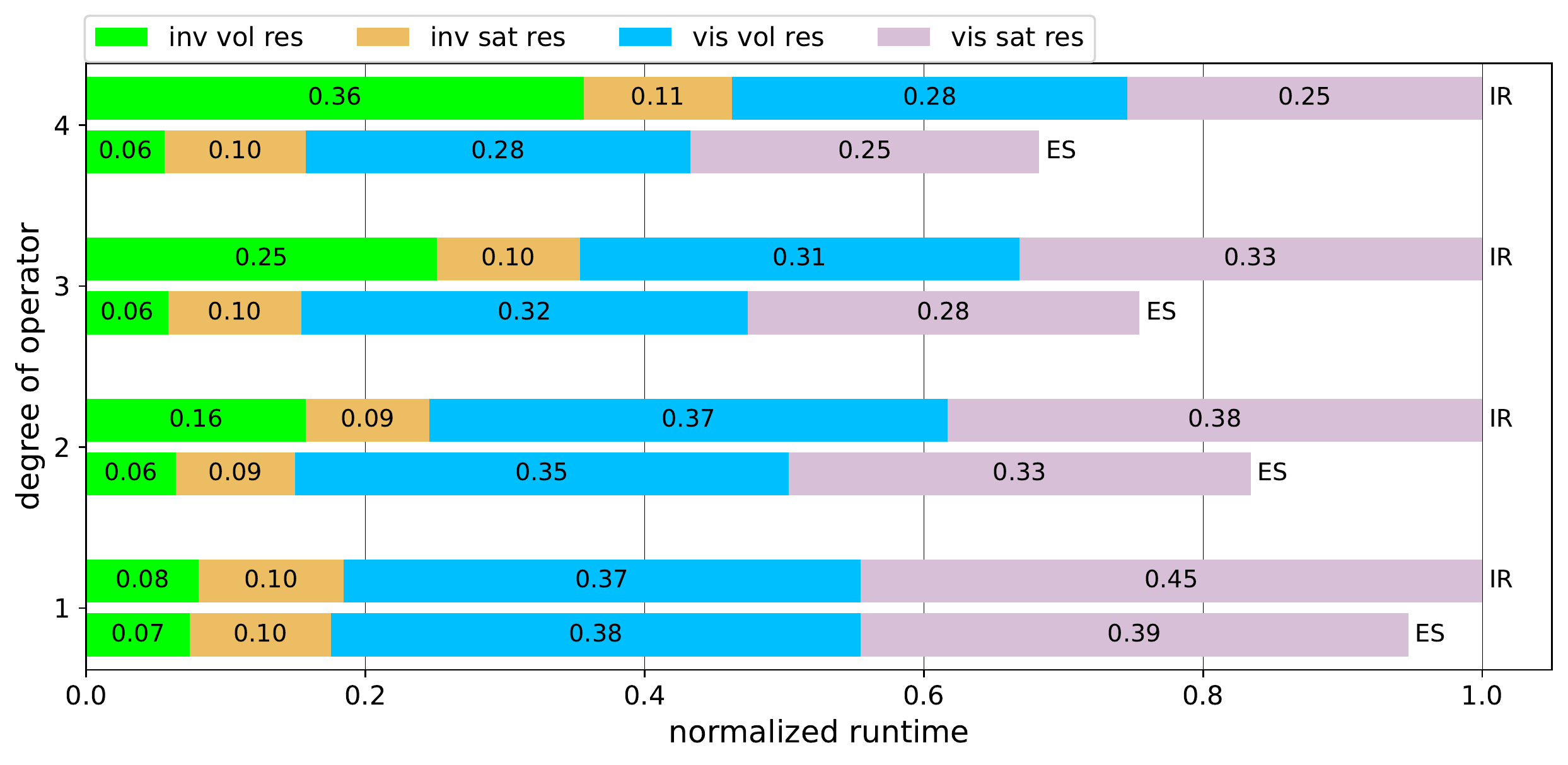}
	\caption{\label{fig:efficiency tgv} Efficiency comparisons of the entropy-split and Ismail-Roe Hadamard-form discretizations of the viscous Taylor-Green vortex problem.}
\end{figure}

An efficiency comparison of the entropy-split and Ismail-Roe Hadamard-form discretizations for the viscous Taylor-Green vortex is depicted in \cref{fig:efficiency tgv}. Again, the normalized average serial performance over five runs, each consisting of 20 timesteps with $ n_e = 750 $ elements, is reported.  While the entropy-split scheme is more efficient than the Ismail-Roe Hadamard-form scheme, the efficiency benefit is not as significant as it is for the Euler equations due to the dominant cost of the viscous terms. However, it is evident from the figure that the proportional cost of the viscous terms is decreasing with increasing operator degree. Hence, for very high degree operators, the efficiency benefit of the entropy-split scheme can be substantial for the \red{compressible} Navier-Stokes equations. 


\section{Conclusions}\label{sec:conclusions}
The entropy-split scheme is extended to element-type implementations using diagonal-norm  diagonal-$ \E $ multidimensional SBP operators and entropy conservative SATs. We show that the scheme is high order accurate and entropy conservative on periodic curvilinear grids for the Euler equations. Its entropy stability does not extend to the \red{compressible} Navier-Stokes equations unless heat fluxes are neglected. Neglecting heat fluxes and applying appropriate viscous SATs, however, we show that the entropy-split scheme is entropy stable for the \red{compressible} Navier-Stokes equations. The viscous SATs are also applicable to existing Hadamard-form discretizations, and they are extensions of the SATs presented in  \cite{yan2018interior,worku2021simultaneous} to systems of equations with a symmetric positive semidefinite diffusivity tensor.

The entropy-split scheme is inherently not conservative in the sense of Lax-Wendroff. To mitigate the effects of this loss of conservation, an entropy stable hybrid scheme consisting of the entropy-split and Sj{\"o}green-Yee Hadamard-form discretizations is proposed. The hybrid scheme is locally conservative, and numerical results show that it converges to the correct weak solution. A matrix-form interface dissipation operator is designed to provide entropy dissipation for the hybrid, entropy-split, and  Sj{\"o}green-Yee Hadamard-form schemes. 

The high-order accuracy of the entropy-split scheme is verified using problems governed by the Euler and Navier-Stokes equations. We found that the solutions to the problems considered converge at rates between $ p $ and $ p+1 $ under mesh refinement. Furthermore, for the entropy-split scheme, the primary conservation errors superconverge at rates close to $ 2p $. For long time simulations of the isentropic vortex problem, we observed that the stability of the entropy-split scheme depends on the value of the arbitrary splitting parameter $ \beta $. However, the stability issues disappear if interface matrix-dissipation is applied. Despite not being entropy stable for the Navier-Stokes equations, no practical instability is observed with the entropy-split scheme when applied to the challenging viscous Taylor-Green vortex problem. \red{In practice, the entropy-split scheme is applied with the heat fluxes included, and despite not being provably entropy stable in this case, the method is likely to be more robust than the standard divergence form discretization for many problems, as the heat fluxes are rarely the source of difficulties in numerical simulations.} In general, however, it can be concluded that the entropy-split scheme is less robust than the Ismail-Roe Hadamard-form scheme due to the problem dependent $ \beta $ parameter, and further research is required to mitigate this robustness issue.

Finally, efficiency comparison tests show that the computational costs of the spatial residual terms of the Euler equations with the entropy-split scheme using the degree three and four SBP diagonal-$ \E $ tetrahedral elements are, respectively, about two and three times smaller than the corresponding costs with the Ismail-Roe Hadamard form scheme. The efficiency benefit for the Navier-Stokes equations is not as large due to the dominant cost of the viscous terms. However, based on the observation that the proportional cost of the viscous terms is decreasing for increasing operator degree, the efficiency benefit of the entropy-split method for the Navier-Stokes equations compared to the Hadamard-form scheme can be substantial for very high order discretizations.

\section*{Declaration of competing interest}
The authors declare that they have no known competing financial interests or personal relationships that could have appeared to influence the work reported in this paper.

\section*{Acknowledgments}
The authors would like to sincerely thank Professor Masayuki Yano and his Aerospace Computational Engineering Lab at the University of Toronto for the use of their software, the Automated PDE Solver (APS), and the helpful implementation discussions. \violet{We are also grateful to the anonymous reviewers for their comments and remarks, which have improved the final manuscript.} Computation were performed on the Niagara supercomputer at the SciNet HPC Consortium \cite{ponce2019deploying}. SciNet is funded by: the Canada Foundation for Innovation; the Government of Ontario; Ontario Research Fund - Research Excellence; and the University of Toronto.

\addcontentsline{toc}{section}{Acknowledgments}

\bibliographystyle{model1-num-names}
{\small
\bibliography{references}
}
\addcontentsline{toc}{section}{\refname}

\end{document}